\documentclass[onefignum,onetabnum]{siamart190516}

\usepackage{amsmath}

\usepackage{amsfonts, psfrag, graphicx, mathtools, stackrel, amssymb,enumitem,indentfirst,float,geometry,color,enumerate,graphicx,subcaption,algorithm,algpseudocode,pifont,hyperref,mathtools,mathrsfs}
\usepackage{pbox}
\usepackage{booktabs, cellspace, hhline}
\allowdisplaybreaks[4]
\usepackage[hyperpageref]{backref}

\usepackage{comment,multicol,xspace}
\usepackage{makecell}

\numberwithin{equation}{section}
\numberwithin{figure}{section}

\usepackage{etoolbox}
\usepackage[title]{appendix}

\usepackage{lipsum}

\newtheorem{thm}{Theorem}[section]
\newtheorem{remark}{Remark}[section]

\newcommand{\commentout}[1]{{}} 

\newcommand{\bfA}{{\bf A}}
\newcommand{\bfa}{{\bf a}}
\newcommand{\bfB}{{\bf B}}
\newcommand{\bfb}{{\bf b}}

\newcommand{\bfc}{{\bf c}}

\newcommand{\bfE}{{\bf E}}
\newcommand{\bfe}{{\bf e}}

\newcommand{\bff}{{\bf f}}

\newcommand{\bfH}{{\bf H}}

\newcommand{\bfJ}{{\bf J}}

\newcommand{\bfn}{{\bf n}}

\newcommand{\bfS}{{\bf S}}
\newcommand{\bft}{{\bf t}}

\newcommand{\bfu}{{\bf u}}

\newcommand{\bfv}{{\bf v}}
\newcommand{\bfw}{{\bf w}}

\newcommand{\bfx}{{\bf x}}

\newcommand{\bfz}{{\bf z}}

\newcommand{\bfphi}{\boldsymbol{\phi}}

\newcommand{\bfxi}{\boldsymbol{\xi}}
\newcommand{\bfeta}{\boldsymbol{\eta}}

\newcommand{\bfPi}{\boldsymbol{\Pi}}

\usepackage{comment,multicol,xspace}
\usepackage[all]{xy} 
\usepackage{tikz-cd}
\usepackage{adjustbox}
\usepackage{multirow}

\makeatletter
\MHInternalSyntaxOn
\def\MT_leftarrow_fill:{%
  \arrowfill@\leftarrow\relbar\relbar}
\def\MT_rightarrow_fill:{%
  \arrowfill@\relbar\relbar\rightarrow}
\newcommand{\xrightleftarrows}[2][]{\mathrel{%
  \raise.55ex\hbox{%
    $\ext@arrow 0359\MT_rightarrow_fill:{\phantom{#1}}{#2}$}%
  \setbox0=\hbox{%
    $\ext@arrow 3095\MT_leftarrow_fill:{#1}{\phantom{#2}}$}%
  \kern-\wd0 \lower.55ex\box0}}
\MHInternalSyntaxOff
\makeatother

\newcommand{\rot}{\operatorname{rot}}
\newcommand{\grad}{\operatorname{grad}}
\newcommand{\curl}{\operatorname{curl}}
\renewcommand{\div}{\operatorname{div}}
\newcommand{\textint}{\operatorname{int}}

\newcommand{\dd}{\,{\rm d}}

\newcommand{\vertiii}[1]{{\left\vert\kern-0.25ex\left\vert\kern-0.25ex\left\vert #1
    \right\vert\kern-0.25ex\right\vert\kern-0.25ex\right\vert}}
    \newcommand{\vertii}[1]{{\left\vert\kern-0.25ex\left\vert #1
    \right\vert\kern-0.25ex\right\vert}}
\newcommand{\verti}[1]{{\left\vert #1
    \right\vert}}

\makeatother

\begin{document}

\headers{$\bfH(\curl)$ IFE Methods}{L. Chen, R. Guo and J. Zou}
\title{A Family of Immersed Finite Element Spaces and Applications to Three Dimensional H(Curl) Interface Problems
}
\author{
Long Chen \thanks{Department of Mathematics, University of California, Irvine, CA 92697, USA. 
The work of this author was funded by NSF DMS-2012465. (lchen7@uci.edu).}
\and Ruchi Guo \thanks{Department of Mathematics, University of California, Irvine, CA 92697, USA. 
The work of this author was supported in parts by HKSAR grant \#15302120. (ruchig@uci.edu). }
\and Jun Zou \thanks{Department of Mathematics, The Chinese University of Hong Kong, Shatin, N.T., Hong Kong. 
The work of this author was substantially supported by Hong Kong RGC General Research Fund 
(projects 14306921 and 14306719). (zou@math.cuhk.edu.hk).}
}


\date{}
\maketitle
\begin{abstract}
Maxwell interface problems are of great importance in many electromagnetic applications. Unfitted mesh methods are especially attractive in 3D computation as they can circumvent generating complex 3D interface-fitted meshes. However, many unfitted mesh methods rely on non-conforming approximation spaces, which may cause a loss of accuracy for solving Maxwell equations, and the widely-used penalty techniques in the literature may not help in recovering the optimal convergence. In this article, we provide a remedy by developing N\'ed\'elec-type immersed finite element spaces with a Petrov-Galerkin scheme that is able to produce optimal-convergent solutions. To establish a systematic framework, we analyze all the $H^1$, $\bfH(\curl)$ and $\bfH(\div)$ IFE spaces and form a discrete de Rham complex. Based on these fundamental results, we further develop a fast solver using a modified Hiptmair-Xu preconditioner which works for both the GMRES and CG methods.
\end{abstract}

\begin{keywords}
Maxwell equations, Interface problems, $\bfH(\text{curl};\Omega)$-elliptic equations, N\'ed\'elec elements, Raviart-Thomas elements, Immersed finite element methods, Petrov-Galerkin formulation, de Rham Complex, preconditioner
\end{keywords}


\section{Introduction}

Efficient and accurate computation of $H(\curl)$ interface problems is of great importance since many applications involve multiple media with different electromagnetic properties. For example, when designing electromagnetic motors and actuators for electric drive vehicles, one must consider the metal-air or metal-metal interface of rotary objects~\cite{1996BrauerRuehl,2014ErcanJaewook}. Another family of critical applications is for non-destructive detection techniques such as EMI sensors, ground-penetrating radar, and Lorentz force eddy current tests~\cite{2006AlbaneseMonk,2002AmmariBaoFleming,2015AmmariChenChenVolkov,2020CHENLIANGZOU} that aim to detect buried metallic objects. 

Let $\bfu$ be an electric or magnetic field in a domain $\Omega$ in $\mathbb R^3$ that covers multiple subdomains representing media with different electromagnetic properties. For simplicity, we assume there are only two subdomains of $\Omega$ denoted as $\Omega^{\pm}$ separated by an interface $\Gamma$.
Then, the mathematical model for the interface conditions of $\bfu$ at $\Gamma$ can be generally written as
\begin{subequations}
 \label{inter_jc}
\begin{align}
[\bfu \times \bfn]_{\Gamma} &:=  \bfu^+\times \bfn -  \bfu^-\times \bfn = 0,  \label{inter_jc_1} \\
[(\alpha \curl\, \bfu)\times \bfn]_{\Gamma} &:=  \alpha^+\curl\, \bfu^+\times \bfn -   \alpha^- \curl\, \bfu^-\times \bfn   = 0,
\label{inter_jc_2} \\
[\beta \bfu\cdot \bfn]_{\Gamma} &:=  \beta^+ \bfu^+\cdot \bfn -   \beta^- \bfu^-\cdot \bfn   = 0,
\label{inter_jc_3}
\end{align}
\end{subequations}
where $\bfn$ denotes the normal vector to $\Gamma$ from $\Omega^-$ to $\Omega^+$, and $\alpha=\alpha^{\pm}$ and $\beta=\beta^{\pm}$ in $\Omega^{\pm}$ are assumed to be positive piecewise constant functions with different electromagnetic meanings in various situations. In fact, the interface conditions in \eqref{inter_jc} can be interpreted from Hodge star operators, see the discussion in Section~\ref{subsec:space} and particularly Diagram~\ref{Hodge_H_local}. In this work, we shall develop IFE spaces to capture those conditions.

The interface condition in \eqref{inter_jc} appears in many formulation of Maxwell equations. For example, the following time-dependent $\curl$-$\curl$ equation can be used to model an electric field $\bfu = \bfE$:
\begin{equation}
\label{Max_mode1}
\epsilon \bfu_{tt} + \sigma \bfu_t + \curl (\mu^{-1}\curl\bfu) = \bfJ_t
\end{equation}
subject to some boundary and initial conditions, where $\epsilon$, $\sigma$, and $\mu$ are the electric permeability, conductivity, and magnetic permeability of the media respectively, and $\bfJ$ is the applied current density. Then solving \eqref{Max_mode1} at each time step by an implicit method yields the equation
\begin{equation}
\label{Max_model3}
\curl (\alpha \curl \bfu ) + \beta \bfu = \bff.
\end{equation}
In this model, the jump conditions in \eqref{inter_jc} have $\alpha=\mu^{-1}$ and $\beta=\epsilon \Delta t^{-2} + \sigma \Delta t^{-1}$. We refer readers to~\cite{2015CaiCao,2000ChenDuZou,1999CiarletZou,2012HiptmairLiZou} for details of this formulation. One can also consider a Helmholtz-type equation~\cite{2009ChenXiaoZhang,2007ChenWangZheng,2011HIPTMAIRMOIOLAPERUGIA,2003Monk} where the parameters in \eqref{inter_jc} become $\alpha=\mu^{-1}_r$ and $\beta = \kappa^2\epsilon_r$ with $\mu_r$ and $\epsilon_r$ being the relative magnetic and electric permeability. 


It is widely known that traditional finite element (FE) methods must be applied on geometrically fitted meshes to solve interface problems, but generating a high-quality fitted mesh is non-trivial, especially for the considered three-dimensional (3D) space. For complex geometry, a global refinement is needed to generate high-quality meshes, which is expensive. Instead, locally modifying the mesh near the interface is more efficient. For instance, the approach in~\cite{2009ChenXiaoZhang} first generates a background unfitted mesh and further triangulates those interface elements to fit the interface. However, the local refinement cannot guarantee shape regularity or even the maximal angle condition~\cite{2005BuffaCostabelDauge}. In particular, for arbitrarily shaped interface surfaces, the so-called slivers having edges nearly coplanar are difficult to be completely removed from tetrahedral meshes~\cite{2000Edelsbrunner,2004PerssonStrang}, which may deteriorate the accuracy. Another approach is to directly use polyhedral meshes to fit the interface which avoids further triangulation~\cite{2020VeigaDassiMascotto,2018CaoChen,2017ChenWeiWen}. For all these approaches, small edges of anisotropic meshes may make the conditioning of the resulted linear system even worse and pose a severer restriction on the step size for time-dependent problems.

In contrast, some special numerical methods are designed for solving interface problems on unfitted meshes which completely avoids the mesh generation issue. Typical examples include penalty-type methods~\cite{2015BurmanClaus,2016BurmanGuzmanSanchezSarkisHcontrast,2020LiuZhangZhangZheng}, non-matching mesh methods~\cite{2016CasagrandeHiptmairOstrowski,2016CasagrandeWinkelmannHiptmairOstrowski,2000ChenDuZou,2008HuShuZou} and immersed finite element (IFE) methods to be discussed in this article. Indeed, all these methods have been successfully applied to various interface problems arising from fluid, elasticity, wave propagation, and so on, see~\cite{2009BeckerBurmanHansbo,2007GroReusken,2020Guo,2020GuoLinLinZhuang,2014HansboLarsonZahedi} and the reference therein. We also refer readers to FDTD methods~\cite{2004ZhaoWei} based on finite difference formulation for Maxwell interface equations. 

The methodology of IFE methods is to encode the jump conditions in the construction of approximation spaces, i.e., the IFE spaces, which then have optimal approximation capabilities for functions satisfying the related jump conditions. One of its distinguishing features from others is that they inherit the degrees of freedom (DoFs), which is particularly advantageous for the considered $\bfH(\curl)$ problems. 

For almost all the aforementioned unfitted-mesh FEMs, one of the prices paid is the loss of the conformity of the approximation spaces. For example, the interface-penalty methods and non-matching mesh methods employ spaces that are discontinuous across the interface; while the IFE spaces are discontinuous across interface edges/faces. For many situations, the discontinuity can be well-handled by suitable penalties such as the Nitsche's penalty~\cite{2015BurmanClaus,2015LinLinZhang}, but it becomes one of the main obstacles for Maxwell-type equations due to the underlying $\bfH(\curl)$ space. Here let us briefly explain the issue. Solutions of \eqref{Max_model3} generally admit low regularity especially near the interface~\cite{2000CostabelDauge,1999MartinMoniqueSerge,2004CostabelDaugeNicaise}. Suppose that $\bfu\in\bfH^s(\curl;\Omega)$, $s>0$, and an approximation space is non-conforming on a face $F$, then the penalty $h^{-1}\int_F[\bfu_h\times \bfn]\cdot[\bfv_h\times \bfn] \dd s$ is generally needed to ensure the stabilization. Then, in the error estimation, one needs the inequality:
\begin{equation}
\label{approx_F}
h^{-1/2}_K \| \bfu - \pi_F \bfu \|_{L^2(F)} \lesssim  h^{s-1}_K\| \bfu \|_{\bfH^s(\curl;K)},
\end{equation} 
where $\pi_F$ is a projection to some polynomial space on $F$. Note that even for a moderate regularity $s=1$, \eqref{approx_F} immediately leads to a failure of convergence. This issue was first realized in~\cite{2001BenBuffaMaday} for a mortar FEM on non-matching meshes. It was also observed by the authors in~\cite{2016CasagrandeHiptmairOstrowski,2016CasagrandeWinkelmannHiptmairOstrowski} for the Nitsche's penalty method (see Theorem 2 in~\cite{2016CasagrandeHiptmairOstrowski}) and recently numerically studied in~\cite{2020RoppertSchoderToth} with a realistic example. 
From the results in the literature, this issue seems to be essential rather than caused by the limitation of analysis techniques.

Some approaches in the literature can overcome this issue. One is to search for a conforming subspace of the underlying non-conforming space. This is indeed true for standard discontinuous Galerkin (DG) spaces on fitted meshes for Maxwell equations~\cite{2005HoustonPerugiaSchneebeli,2004HoustonPerugiaSchotzau} where the analysis is completely different from the standard DG framework. The authors in~\cite{2008HuShuZou} employ a similar idea to improve the estimate of the mortar FEM, which relies on a ``nested mesh" assumption. But since the unfitted mesh methods modify their spaces near the interface, the conforming subspaces may not exist. Without resorting to conformity, another approach is to use the appropriate scaling in the stabilization~\cite{2020VeigaDassiMascotto,2021CaoChenGuo}. In fact, the scaling factor is highly related to the trace spaces of the underlying Sobolev space ~\cite{2008BrennerCuiLiSung}.
For example, if the solution has $\bfH^1(\curl;\Omega)$ regularity, the scaling factor should be $h$ instead of $h^{-1}$ in \eqref{approx_F}, which can indeed achieve optimal convergence. But for many unfitted mesh methods, such a scaling factor may not yield a stable scheme.



\begin{equation}
   \label{IFE_deRham}
\begin{tikzcd}
  \mathbb{R} \arrow[r,"\hookrightarrow",shorten >= 1ex] 
   & H^2(\beta;\mathcal{T}_h) \arrow[d, "I^n_h"] \arrow[r, "\nabla"] 
   & \bfH^1(\curl,\alpha,\beta;\Omega) \arrow[d, "I^e_h"] \arrow[r, "\curl"] 
   & \bfH^1(\div,\alpha;\mathcal{T}_h) \arrow[d, "I^f_h"] \arrow[r, "\div"] 
   & H^1(\mathcal{T}_h) \arrow[d, "\Pi^0"] \arrow[r] & 0 
   \\
  \mathbb{R} \arrow[r, "\hookrightarrow",shorten >= 1ex] 
   & S^n_h   \arrow[d, "(I^n_h)^{-1}", shift left = 0.5ex]  \arrow[r, "\nabla", shorten <= 1ex, shorten >= 1ex] 
   & \bfS^e_h \arrow[d, "(I^e_h)^{-1}", shift left = 0.5ex] \arrow[r, "\curl", shorten <= 1ex, shorten >= 1ex] 
   & \bfS^{f}_h \arrow[d, "(I^f_h)^{-1}", shift left = 0.5ex] \arrow[r, "\div", shorten <= 1ex, shorten >= 1ex] 
   & Q_h \arrow[d, "I", shift left = 0.5ex] \arrow[r, shorten <= 1ex] & 0
    \\
    \mathbb{R} \arrow[r, "\hookrightarrow",shorten >= 1ex] 
   & \widetilde{S}^n_h \arrow[u, "I^n_h", shift left = 0.5ex]  \arrow[r, "\nabla", shorten <= 1ex, shorten >= 1ex] 
   & \widetilde\bfS^e_h \arrow[u, "I^e_h", shift left = 0.5ex] \arrow[r, "\curl", shorten <= 1ex, shorten >= 1ex] 
   & \widetilde\bfS^{f}_h \arrow[u, "I^f_h", shift left = 0.5ex] \arrow[r, "\div", shorten <= 1ex, shorten >= 1ex] 
   & Q_h \arrow[u, "I", shift left = 0.5ex] \arrow[r, shorten <= 1ex] & 0
\end{tikzcd}
\end{equation}

Therefore, it would be highly desirable if the penalty can be completely avoided, which is the key motivation for this work. With this motivation, here we employ a Petrov-Galerkin (PG) formulation that uses IFE trial spaces but keeps the standard conforming FE test spaces. Then, there will be no penalty introduced from integration by parts, which nicely addresses the issue of \eqref{approx_F}. In order for the resulted matrix being square, the trial and test spaces should have the same number of DoFs. Thus, we shall design the IFE spaces to be isomorphic to the standard FE spaces through the usual DoFs. The isomorphism is shown by the mapping between the middle and bottom sequences in the Diagram \eqref{IFE_deRham}, where the definition of the spaces and operators is given in Section~\ref{subsec:complex}.  The original idea of PG formulation can be found in the fundamental work~\cite{1994BabuskaCalozOsborn} of Babu\v{s}ka et al. and is then used in~\cite{2010ChuGrahamHou} for Multiscale FEM and~\cite{2013HouSongWangZhao} for IFE methods on $H^1$ interface problems. This is a generalization of the previous work in~\cite{2020GuoLinZou} for $H(\curl)$ interface problems from 2D to 3D. However, the generalization is non-trivial due to the more complicated 3D geometry and the well-known essential difference between 2D and 3D $\bfH(\curl)$ spaces.

There are multiple novelties in this work. We construct both the 3D $\bfH(\curl)$ and $\bfH(\div)$ IFE spaces according to their corresponding jump conditions which have not appeared in any literature. With the $H^1$ IFE space in the literature, we establish the complete de Rham complex as shown by the middle sequence in Diagram \eqref{IFE_deRham} that is identical to the FE counterparts, including the exactness and commutative properties. It is well known that solving the resulting linear system from the $\bfH(\curl)$ system is challenging, and the exact sequence is critical for developing fast solvers as it can accurately describe the kernel space of $\curl$. 
Based on the de Rham complex and multigrid methods, we develop an IFE version of HX preconditioner~\cite{2007HiptmairXu}. In addition, the non-SPD matrix and the potentially small-cutting elements add another layer of difficulty for many unfitted mesh methods. We design a special smoother by computing the exact inverse of the non-SPD submatrix near the interface. The idea is closely related to domain decomposition methods~\cite{1998XuZou}. We demonstrate by numerical experiments that the resulting fast solver works remarkably well for general minimum residual methods (GMRES) and even Conjugate-Gradient (CG) methods which typically do not converge for non-SPD systems. It is worthwhile to mention that CG-like methods have been developed for more general semi-symmetric matrices~\cite{1986MakinsonShah}. We also refer readers to multigrid methods for interface problems in~\cite{2008XuZhu,2021HuWang}.


The article consists of 5 additional sections. In the next section, we describe the geometry and set up some notations. In Section~\ref{sec:IFEspaces}, we establish the IFE spaces and their de Rham complex. In Section~\ref{sec:formulation}, we describe the PG-IFE formulation, the fast solver, and a numerical test for the \textit{inf-sup} condition. In Section~\ref{sec:approx}, we show the optimal approximation capabilities. In Section~\ref{sec:numer}, we present some numerical experiments.


\section{Preliminary}

In this section, we first describe an unfitted background mesh and introduce some spaces and notation which will be frequently used in this work.

\subsection{Meshes}
\begin{figure}[h]
\centering
\begin{subfigure}{.2\textwidth}
    \includegraphics[width=1.05in]{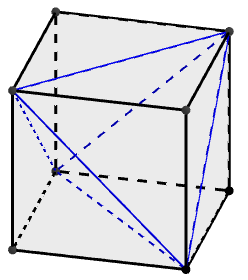}
     \label{fig:cubetet1} 
\end{subfigure}
~~~
 \begin{subfigure}{.2\textwidth}
    \includegraphics[width=1.05in]{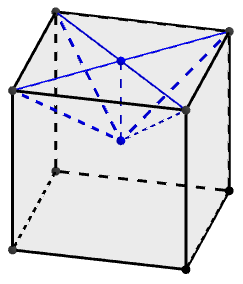}
     \label{fig:cubetet3} 
\end{subfigure}
\begin{subfigure}{.18\textwidth}
    \includegraphics[width=0.8in]{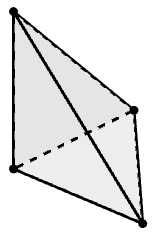}
     \label{fig:righttet} 
\end{subfigure}
\begin{subfigure}{.2\textwidth}
     \includegraphics[width=1.2in]{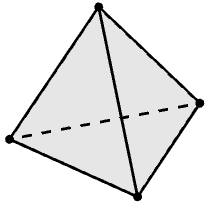}
     \label{fig:regulartet} 
\end{subfigure}
     \caption{Triangulation of a cube. The left two plots show two triangulation approaches resulting in 5 tetrahedra (Type I) and 24 tetrahedra (Type II). The right two plots show two types of tetrahedra from the triangulation: a trirectangular tetrahedron and a regular tetrahedron.}
  \label{fig:cubetet} 
\end{figure}

Note that electromagnetic wave generally propagates over the whole space, and a widely used approach is to cover the interesting media by a box as a modeling and computational domain. So in this work, we shall assume our domain is cubic, and generate a background Cartesian mesh in which each cubic element is then partitioned into several tetrahedra. For example, we show two types of triangulation in the left two plots of Figure~\ref{fig:cubetet} where the triangulation in Type I does not introduce extra nodes but Type II introduces an extra node on each face and one inside the element. (We only show 4 tetrahedra in the second plot of Figure~\ref{fig:cubetet}\subref{fig:cubetet3} for better illustration).
All these tetrahedra form an interface-unfitted shape regular mesh denoted by $\mathcal{T}_h$ which will be used for discretization. Let $\mathcal{N}_h$, $\mathcal{E}_h$, and $\mathcal{F}_h$ be the collection of nodes, edges, and faces of $\mathcal{T}_h$, respectively. For each element $K$, we let $\mathcal{N}_K$, $\mathcal{E}_K$, and $\mathcal{F}_K$ be the collection of nodes, edges, and faces of $K$ with $|\mathcal{N}_K|=4$, $|\mathcal{E}_K|=6$, and $|\mathcal{F}_K|=4$.

Given the signed-distance level-set function of the interface denoted by $\varphi^{\Gamma}$, we can compute its linear interpolation $\varphi^{\Gamma}_h$ on the mesh $\mathcal{T}_h$; namely, $\varphi^{\Gamma}_h$ interpolates $\varphi^{\Gamma}$ at the mesh nodes and is piecewise linear on each element. Then, the numerical interface is defined as
\begin{equation}
\label{NumInterf}
\Gamma_h = \{ \bfx~:~ \varphi^{\Gamma}_h(\bfx) =0 \}.
\end{equation}
We note that, as $\varphi^{\Gamma}_h$ is piecewise linear, $\Gamma_h$ is also piecewise linear and forms a polyhedron as a linear approximation to the exact interface. Linear interpolation of signed-distance functions has been widely used for discretizing interface numerically~\cite{1999Sethian,2001OsherFedkiw}. For a smooth interface, the geometric error is in the order of $\mathcal{O}(h^2)$, see Lemma~\ref{lemma_interface_flat} below, which is sufficient for a linear method~\cite{2010LiMelenkWohlmuthZou}. Then, we let $\Gamma_h$ cut $\Omega$ into $\Omega^{\pm}_h$. An element $K$ is called an interface element if 
\begin{equation}
\label{interfElem}
\min_{\bfx\in\mathcal{N}_K}\varphi^{\Gamma}_h(\bfx) \max_{\bfx\in\mathcal{N}_K}\varphi^{\Gamma}_h(\bfx) < 0.
\end{equation}
As $\varphi^{\Gamma}_h$ interpolates $\varphi^{\Gamma}$ at the mesh nodes, \eqref{interfElem} also implies that $K$ intersects the exact interface. Let $\mathcal{T}^i_h$ be the collection of interface elements.

\begin{figure}[h]
\centering
\begin{minipage}{.42\textwidth}
  \centering
  \includegraphics[width=2in]{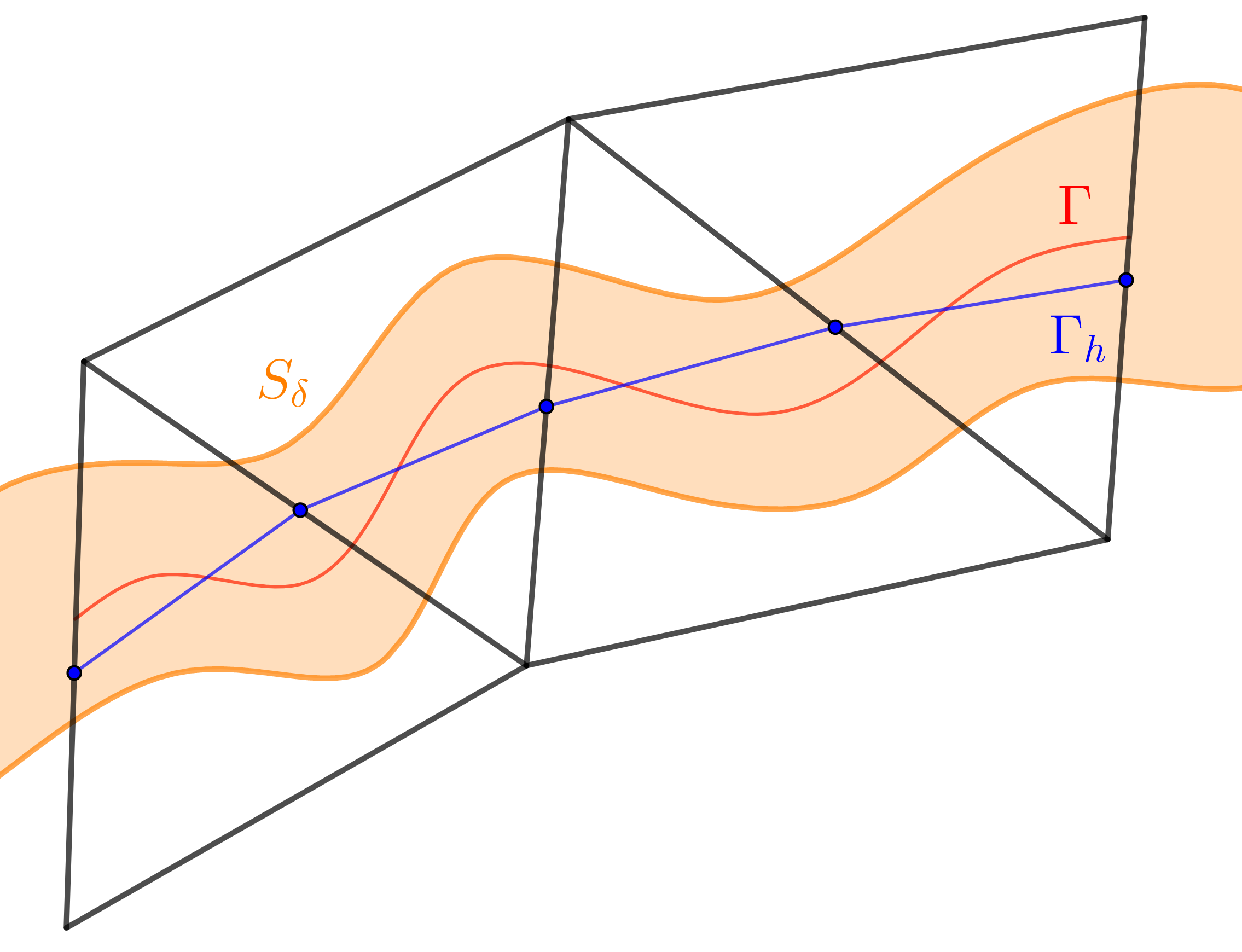}
  \caption{The mismatching region and $\delta$ Strip.}
  \label{fig:deltasrip}
\end{minipage}
\hspace{1cm}
\begin{minipage}{.4\textwidth}
  \centering
  \ \includegraphics[width=1in]{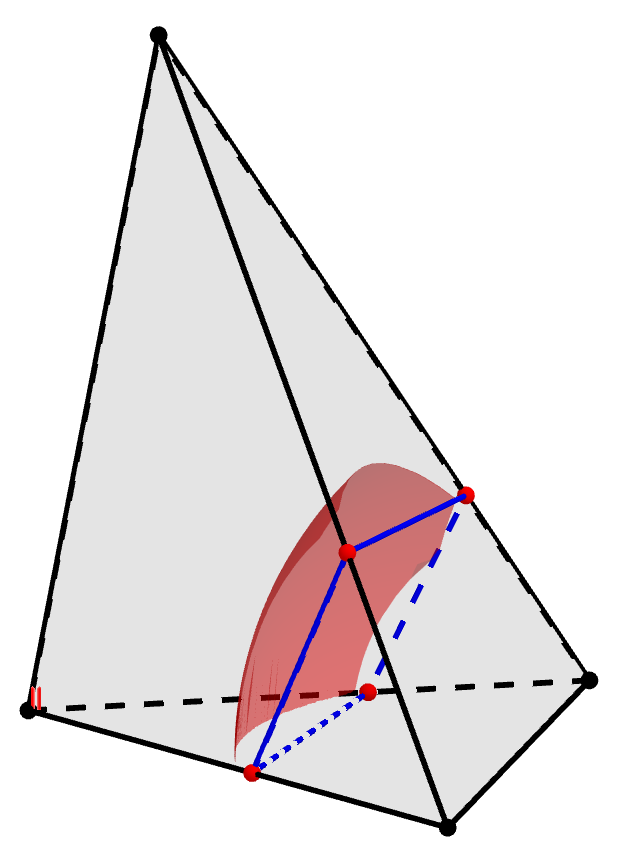}
  \caption{The intersection points of $\Gamma^K_h$ are always coplanar while those of $\Gamma$ may not.}
  \label{fig:NumInterf}
  \end{minipage}
\end{figure}

For IFE methods, there are two advantages to using $\varphi^{\Gamma}_h$. First, as $\Gamma^K_h $ is a plane, it is not hard to see that there are only two possible cases for $\Gamma^K_h$ cutting a tetrahedron: three or four intersection points, see Figure~\ref{fig:regulartetInter}. Second, it is possible that the four intersection points of the original interface surface are not coplanar; but the cutting points of $\Gamma^K_h$ are always coplanar, as shown in Figure~\ref{fig:NumInterf}. Of course, the cutting points of $\Gamma$ may not be those of $\Gamma_h$, see Figure~\ref{fig:deltasrip} for a 2D illustration where $\Gamma^K_h$ does not pass through the cutting points of $\Gamma$. With these two properties, the IFE shape functions always exist, see the discussion in Section~\ref{sec:IFEspaces}.

In addition, given a $K\in\mathcal{T}^i_h$, define the interface patch $\Gamma^K_h = \Gamma_h\cap K$. Of course, there is a mismatching region sandwiched by $\Gamma_h$ and $\Gamma$, see the subregion enclosed by the red curve and blue polyline in Figure~\ref{fig:deltasrip} for a 2D illustration. We define it as $K_{\textint}:=(K^+\cap K^-_h)\cup(K^-\cap K^+_h)$ and $\Omega_{\textint}:= \cup_{K\in\mathcal{T}^i_h} K_{\textint}$. 

We present the following geometric error estimate of which the proof is in Appendix~\ref{lemma_interface_flat}. 
\begin{lemma}
\label{lemma_interface_flat}
Suppose $\Gamma$ is smooth enough such that $\varphi^{\Gamma}$ is also smooth near the interface, and let the mesh size be sufficiently small. Then, for each interface element $K$ and each $\bfx\in\Gamma\cap K$, there holds
\begin{subequations}
\label{lemma_interface_flat_eq0}
\begin{align}
    & \| \bfx - \bfx^{\bot} \| \le  c_{\Gamma,1} h^2_K, \label{lemma_interface_flat_eq01}   \\
    &  \| \bfn - \bar{\bfn}_K \| \le  c_{\Gamma,2} h_K, \label{lemma_interface_flat_eq02}   
\end{align}
\end{subequations}
where $\bfx^{\bot}$ is the projection of $\bfx$ onto $\Gamma^K_h$, and $\bar{\bfn}_K$ is the normal vector to $\Gamma^K_h$ with the same orientation as the normal vector to interface, i.e., $\bfn$. Here, the constants $c^1_{\Gamma}$ and $c^2_{\Gamma}$ only depend on the interface curvature.
\end{lemma}


\subsection{Spaces}
\label{subsec:space}

Given a subdomain $D\subseteq \Omega$, for a non-negative number $s$, we let $H^s(D)$ be the standard Sobolev space. In addition, we introduce the following vector function spaces:
\begin{subequations}
\label{sob_space}
\begin{align}
    & \bfH^s(\curl;D) = \{ \bfu\in \bfH^s(D): \curl \, \bfu\in \bfH^s(D) \},  \\
    & \bfH^s(\div;D) = \{ \bfu\in \bfH^s(D): \div \, \bfu\in H^s(D) \}.
\end{align}
\end{subequations}
Given a subdomain $D\subseteq \Omega$, if $D\cap\Gamma \neq \emptyset$, we let $D^{\pm}=\Omega^{\pm}\cap D$, and further denote $H^s(\cup D^{\pm})$, $\bfH^s(\curl;\cup D^{\pm})$ and $\bfH^s(\div;\cup D^{\pm})$ consisting of functions that belong to the corresponding Sobolev spaces on each $D^{\pm}$.
Given a face $F$ in the mesh, we let $\rot_F$ be the 2D rotationary operator on $F$, 
and define the space
\begin{equation}
\label{2Drot_space}
\bfH^s(\rot;F) = \{\bfu\in \bfH^s(F): \rot_F \, \bfu\in H^s(F) \}.
\end{equation}

As the IFE version of the de Rham complex is to be established, we need to additionally revisit the $H^1$ and $\bfH(\div)$ interface problems and understand their connection. 
The $H^1$-interface problem reads as
\begin{equation} 
\label{model_H1}
\begin{aligned}
      -\nabla\cdot(\beta\nabla u)&=f \quad \text{ in }  \Omega^-  \cup \Omega^+, 
\end{aligned}
\end{equation}
with $f\in L^2(\Omega)$, subject to $u=0$ on $\partial \Omega$ and the jump conditions:
\begin{subequations}
\label{H1_jump}
\begin{align}
       [u]_{\Gamma} & := u^+ - u^- = 0, \\
  [\beta \nabla u \cdot \bfn]_{\Gamma} &:=  \beta^{+}\nabla u^{+}\cdot \mathbf{ n}  - \beta^{-}\nabla u^{-}\cdot \mathbf{ n} =0,
\end{align}
\end{subequations}
where the parameter $\beta$ is the same as the one in \eqref{inter_jc_3} that has the similar physical meaning related to conductivity. The $\bfH(\div)$-interface problem is given by
\begin{equation}
\label{model_Hidv}
- \nabla \div(\bfu) + \alpha \bfu = \bff,
\end{equation}
with $\bff \in \bfH(\curl;\Omega)$, subject to $\bfu\cdot\bfn=0$ on $\partial \Omega$ and the jump conditions
\begin{subequations}
\label{Hdiv_jump}
\begin{align}
       [\bfu\cdot\bfn]_{\Gamma} & := \bfu^+\cdot\bfn - \bfu^-\cdot\bfn = 0, \\
  [\alpha  \bfu \times \bfn]_{\Gamma} &:=  \alpha^{+} \bfu^{+} \times \mathbf{ n}  - \alpha^{-} \bfu^{-} \times \mathbf{ n} = \mathbf{ 0}, \\
  [\div(\bfu)]_{\Gamma} : &= \div(\bfu^+) - \div(\bfu^-) = 0.
\end{align}
\end{subequations}
This system comes from the first-order least-squares formulation of the elliptic interface problem~\cite{1994CaiLazarovManteuffelMcCormick} or preconditioning for the mixed finite element using a gradient formulation~\cite{1997ARNOLDRICHARDWINTHER}. The related finite element approximation for interface problems can be found in~\cite{2010HiptmairLiZou}.

Now, we encode the jump conditions for the $H^1$ case \eqref{H1_jump}, the $\bfH(\curl)$ case \eqref{inter_jc}, and the $\bfH(\div)$ case \eqref{Hdiv_jump} into the definition of the following special Sobolev spaces:
\begin{subequations}
\label{sob_space_inter}
\begin{align}
     H^{2}(\beta;D) &:= \{ u\in H^{2}(\cup D^{\pm})\cap H^1(D) :  \beta \nabla u \in \bfH(\div;D) \} ,  \label{sob_space_inter1}\\
    \bfH^1(\alpha,\beta,\curl;D) & := \{ \bfu\in \bfH^1(\curl;\cup D^{\pm})\cap \bfH(\curl;D): \beta\bfu \in \bfH(\div;D),~~ \alpha\curl\bfu\in \bfH(\curl;D)\},  \label{sob_space_inter2} \\
    \bfH^1(\alpha,\div;D) &:= \{ \bfu\in \bfH^1(\div;\cup D^{\pm})\cap \bfH(\div;D): \alpha \bfu \in \bfH(\curl;D), ~~  \div \, \bfu\in H^1(D) \} . \label{sob_space_inter3}
\end{align}
\end{subequations}
In fact, the other interpretation of these interface conditions is to regard $\alpha$ and $\beta$ as Hodge star operators mapping between fields in different Sobolev spaces, which is illustrated by the following Diagram \eqref{Hodge_H_local}.
The construction of IFE spaces is just to mimic it locally on each interface element.
\begin{equation}
\label{Hodge_H_local}
\left.\begin{array}{ccccccc}
H^{2}(\beta;\Omega ) &  \xrightarrow[]{~~\grad~~} & \bfH^1(\curl,\alpha,\beta; \Omega) & \xrightarrow[]{~~\curl~~} &  \bfH^1(\div,\alpha; \Omega) & \xrightarrow[]{~~\div~~} & H^1(\Omega) \\
 \bigg\downarrow I & &   \bigg\downarrow \beta &  &  \bigg\downarrow \alpha & &   \bigg\downarrow I  \\
 L^2(\Omega) & \xleftarrow[]{~~~ \div ~~~}  & \bfH(\div;\Omega) &   \xleftarrow[]{~~~ \curl ~~~} & \bfH(\curl;\Omega)  &  \xleftarrow[]{~~~ \grad ~~~} & H^1(\Omega) 
\end{array}\right. 
\end{equation}



\section{Immersed Finite Element Spaces}
\label{sec:IFEspaces}

In this section, we develop N\'ed\'elec-type and Raviart-Thomas-type IFE spaces that admit the edge and face DoFs, respectively. 
We shall use $\alpha_h$ and $\beta_h$ for the discontinuous coefficients partitioned by $\Gamma_h$ instead of $\Gamma$. In the following discussion, let $\mathbb{P}_k(D)$ be the $k$-th degree polynomial space on a domain $D$. For readers' sake, we recall the lowest order local N\'ed\'elec~\cite{1980Nedelec,1986Nedelec} and Raviart-Thomas spaces~\cite{1977RaviartThomas}: 
\begin{subequations}
\label{NeRT}
\begin{align}
    \mathcal{ND}_0(K):= \{ \bfa\times\bfx + \bfb: \bfa\in\mathbb{R}^3,~~ \bfb\in\mathbb{R}^3 \},  ~~~~ \mathcal{RT}_0(K):= \{ c~\bfx + \bfa: c\in \mathbb{R}, ~~  \bfa\in\mathbb{R}^3 \}.
\end{align}
\end{subequations}


\subsection{The local $H^1$, $\bfH(\curl)$ and $\bfH(\div)$ IFE Spaces}

Let us briefly recall the local IFE spaces~\cite{2022CaoChenGuo,2005KafafyLinLinWang}. Consider two piecewise constant vector spaces:
\begin{subequations}
\label{AB_space1}
\begin{align}
    &  \bfA_h(K) = \{ \bfc: \bfc^{\pm}= \bfc|_{K^{\pm}_h} \in \left[\mathbb{P}_0(K^{\pm}_h)\right]^3, ~ \bfc\in \bfH(\div;K), ~ \alpha_h\bfc \in \bfH(\curl;K) \}, \\
    &   \bfB_h(K) = \{ \bfc: \bfc^{\pm}= \bfc|_{K^{\pm}_h} \in \left[\mathbb{P}_0(K^{\pm}_h)\right]^3, ~ \bfc\in \bfH(\curl;K), ~ \beta_h\bfc \in \bfH(\div;K) \}.
\end{align}
\end{subequations}
We can derive the specific format for functions in \eqref{AB_space1}. Given each interface element with the approximate plane $\Gamma^K_h$, let $\bar{\bfn}_K$, $\bar{\bft}_{K,1}$ and $\bar{\bft}_{K,2}$ be its unit normal vector and two orthogonal unit tangential vectors. Further denote the transformation matrix by $T_K=[\bar{\bfn}_K, \bar{\bft}_{K,1},\bar{\bft}_{K,2}]$. Define the transformation matrices: 
\begin{equation}
\label{lem_ABeigen_eq1}
A_K = T_K \left[\begin{array}{ccc} 1 & 0 & 0 \\ 0 & \tilde{\alpha} & 0 \\ 0 & 0 & \tilde{\alpha} \end{array}\right] T^{\top}_K ~~~\text{and} ~~~ 
B_K = T_K \left[\begin{array}{ccc} \tilde{\beta} & 0 & 0 \\ 0 & 1 & 0 \\ 0 & 0 & 1 \end{array}\right] T^{\top}_K,
\end{equation}
where $\tilde{\alpha}=\alpha^+/\alpha^-$ and $\tilde{\beta}=\beta^+/\beta^-$. Then,
the spaces $\bfA_h(K)$ and $\bfB_h(K)$ can be rewritten as
\begin{equation}
\label{AB_space2}
     \bfA_h(K) = \{  \bfc: \bfc^{\pm} \in \left[\mathbb{P}_0(K^{\pm}_h)\right]^3, ~ \bfc^- = A_K \bfc^+ \},   ~~~  \bfB_h(K) = \{  \bfc: \bfc^{\pm} \in \left[\mathbb{P}_0(K^{\pm}_h)\right]^3, ~ \bfc^- = B_K \bfc^+ \}. 
\end{equation}

\begin{table}[htp]
  \centering 
	\renewcommand{\arraystretch}{2}  
  \begin{tabular}{ c c c c } 
  \toprule
 IFE spaces & $S^n_h(K)$ & $\bfS^e_h(K)$ & $\bfS^f_h(K)$  \medskip  \\ 
 \hline
\makecell{ Sobolev \\ spaces} & $H^1(K)$ & $\bfH(\curl;K)$ & $\bfH(\div;K)$ \medskip \\ 
 \hline 

\makecell{ Function \\ format } & \makecell{ $\bfb_h\cdot(\bfx-\bfx_K)+c$ \\ 
 $ \bfb_h \in \bfB_h(\beta_h;K),$\\
 $ c\in \mathbb{P}_0(K)$ } & \makecell{ $\bfa_h\times(\bfx-\bfx_K) + \bfb_h$ \\ 
 $ \bfa_h \in \bfA_h(\alpha_h;K),$\\
 $ \bfb_h \in \bfB_h(\beta_h;K)$ } &  \makecell{ $c(\bfx-\bfx_K)+\bfa_h$\\ $c\in \mathbb{P}_0(K),$\\
 $ \bfa_h \in \bfA_h(\alpha_h;K)$ } \medskip  \\ 
 \hline
 \makecell{Jump \\ conditions} 
 & \makecell{ $[v_h]_{\Gamma^K_h} =0$ \\ $[\beta_h \nabla v_h\cdot \bar{\bfn}]_{\Gamma^K_h} = 0$ } 
 & \makecell{ $[\bfv_h\times \bar{\bfn}]_{\Gamma^K_h} = \mathbf{ 0}$ \\ $[\alpha_h\curl\bfv_h\times \bar{\bfn}]_{\Gamma^K_h} = \mathbf{ 0}$ \\ $[\beta_h \bfv_h\cdot \bar{\bfn}]_{\bfx_K} = { 0}$ }
 & \makecell{ $[\bfv_h\cdot \bar{\bfn}]_{\Gamma^K_h} = { 0}$ \\ $[\alpha_h\bfv_h\times \bar{\bfn}]_{\bfx_K} = \mathbf{ 0}$ \\ $[\div \bfv_h ]_{\Gamma^K_h} = { 0}$ } \medskip \\
 \bottomrule
\end{tabular}
  \caption{IFE spaces, their function format, the jump conditions and the Sobolev spaces to which they belong.}
  \label{tab_IFE_func}
\end{table}
With these two fundamental spaces, we are able to construct the $H^1$, $\bfH(\curl)$, and $\bfH(\div)$ IFE spaces, denoted by $S^n_h(K)$, $\bfS^e_h(K)$, and $\bfS^f_h(K)$, which are reported in Table~\ref{tab_IFE_func} including their corresponding jump condition and general function format. One can easily check that the IFE spaces above belong to the corresponding Sobolev spaces. 
Moreover, the formats are identical to the standard Lagrange, N\'ed\'elec, and Raviart-Thomas elements in which the only difference is the vectors $\bfa$ and $\bfb$ replacing the simple constant vectors in $[\mathbb{P}_0(K)]^3$. Then, it is easy to see the following sequence
\begin{equation}
\label{local_IFE_exatSeq}
\left.\begin{array}{ccccccccccc}
\mathbb{R}& \xrightarrow[]{~~\hookrightarrow~~} &S^{n}_h(K) &  \xrightarrow[]{~~\nabla~~} & \bfS^e_h(K) & \xrightarrow[]{~~\curl~~} & \bfS^{f}_h(K) & \xrightarrow[]{~~\div~~} & \mathbb{P}_0(K) & \xrightarrow & 0.
\end{array}\right. 
\end{equation}

\subsection{Degrees of Freedom}

For the $H^1$ case, it has been proved in~\cite{2005KafafyLinLinWang} that the functions in $S^n_h(K)$ admit nodal DoFs for the two types of tetrahedra shown in the right two plots of Figure~\ref{fig:cubetet}. Namely, there exist $\phi^n_i$, $i=1,2,3,4$, such that $\phi^n_i(\bfz_j)=\delta_{ij}$, where $\bfz_j$, $j=1,2,3,4$, are the nodes of $K$. Then, the global $H^1$ IFE space is defined as
\begin{equation}
\begin{split}
\label{H1_IFE_globSpace}
S^n_h = \{ v_h\in L^2(\Omega)~:~& v_h|_K\in S^n_h(K), ~ \text{if} ~ K\in \mathcal{T}^i_h, ~ \text{and} ~ v_h|_K\in \mathbb{P}_1(K), ~ \text{if} ~ K\in \mathcal{T}^n_h \\
& v_h ~ \text{is continuous at} ~ \bfx\in \mathcal{N}_h \}.
\end{split}
\end{equation}
Although the local $\bfH(\curl)$ and $\bfH(\div)$ IFE spaces are constructed in~\cite{2022CaoChenGuo}, their edge and face DoFs have not been addressed yet, which is our focus in this subsection. 
Thanks to the special geometry of the considered tetrahedra, the discussion will be relatively easier. 


Here, the most difficult one is the edge space $\bfS^e_h(K)$. But, we will see that the edge shape functions can be constructed from the nodal and face shape functions in a manner exactly the same as the standard FE spaces. 
We shall see that the local exact sequence for the IFE spaces in \eqref{local_IFE_exatSeq} can help in simplifying the analysis. 

Let us first consider the unisolvence of $\bfS^f_h(K)$ and begin by deriving a formula for computing $\bfH(\div)$ IFE shape functions. Given an interface element $K$ with $4$ faces $F_i$, $i=1,2,...,4$, we let $v_i=\int_{F_i} \bfv_h\cdot\bfn \dd s$, with $\bfn_{F_i}$ being the normal vector to $F_i$, for any function $\bfv_h=c(\bfx-\bfx_K) + \bfa_h \in \bfS^f_h(K)$ based on Table~\ref{tab_IFE_func}, and let $v_0=v_1+v_2+v_3+v_4$. Determining $c$ is straightforward:
\begin{equation}
\label{comput_c0}
c = \frac{\div(\bfv_h)}{3} = \frac{1}{3|K|} \int_K \div(\bfv_h) \dd \bfx = \frac{1}{3|K|} \sum_{i=1}^4 \int_{F_i} \bfv_h\cdot\bfn_{F_i} \dd s =  \frac{ v_0}{3|K|}.
\end{equation}
Then, by \eqref{AB_space2} $\bfa_h$ should satisfy the equations: $\forall i=1,...,4$,
\begin{equation}
\label{ConstructHdivIFE_eq1}
|F^+_i| \bfa^+_h\cdot\bfn_{F_i} + |F^-_i| \bfa^-_h\cdot\bfn_{F_i} =  |F^+_i| \bfn_{F_i}^{\top}\bfa^+_h + |F^-_i| \bfn_{F_i}^{\top}A_K\bfa^+_h = v_i -  \frac{ v_0}{3|K|}\int_{F_i}(\bfx-\bfx_K)\cdot\bfn_{F_i} \dd s.
\end{equation}
Considering \eqref{ConstructHdivIFE_eq1} on three arbitrary faces $F_1$, $F_2$ and $F_3$, we obtain a linear system for the vector $\bfa^+_h$:
\begin{equation}
\label{ConstructHdivIFE_eq2}
M_K\bfa_h^+ = \bar{\bfv},~~~~ \text{with} ~~ M_K =\Phi^+N_K+ \Phi^-N_KA_K,
\end{equation}
where $\Phi^{\pm} = \text{diag}( |F_1^{\pm}|, |F_2^{\pm}|, |F_3^{\pm}| )$, $N_K = [\bfn_{F_1},\bfn_{F_2},\bfn_{F_3} ]^{\top}$ and $\bar{\bfv} $ is given by the right-hand side of \eqref{ConstructHdivIFE_eq1}. Thus, the solvability of $\bfa^+_h$ is equivalent to the invertibility of $M_K$.

\begin{lemma}
\label{lem_Hdiv_unisol}
Given a trirectangular tetrahedral or regular tetrahedral interface element $K$, the space $\bfS^f_h(K)$ has the DoFs $\int_{F_i} \bfv_h\cdot\bfn_{F_i} \dd s$, $i=1,2,...,4$.
\end{lemma}
\begin{proof}
It remains to show the invertibility of $M_K$. We treat the two types of tetrahedra separatly. 

If $K$ is a trirectangular tetrahedron, we chose $F_i$, $i=1,2,3$, as the three orthogonal faces. Then, $N_K$ is an orthonormal matrix, and we can rewrite
\begin{equation}
\label{Hcurl_lem_unisol_eq1}
M_KN^{\top}_K = \Phi^+ + \Phi^-N_K A_K N^{\top}_K = \Phi + \Phi^- N_K (A_K-I) N^{\top}_K = \Phi\left( I + \Phi^{-1}\Phi^-  N_K (A_K-I) N^{\top}_K \right),
\end{equation} 
where $\Phi = \text{diag}( |F_1|, |F_2|, |F_3| ) = \Phi^- + \Phi^+$. Without loss of generality, we can assume $\tilde{\alpha}=\alpha^+/\alpha^- \le 1$; otherwise we simply consider the matrix $M_KA^{-1}_K N^{\top}_K  = \Phi^+N_K A_K N^{\top}_K + \Phi^- $. Let $\lambda_{\max}(\cdot)$ be the largest magnitude of eigenvalues of a matrix. Note that $\lambda_{\max}(\Phi^{-1}\Phi^- )\le 1$ and $\lambda_{\max}(N_K (A_K-I) N^{\top}_K) = \lambda_{\max}(A_K-I) = 1 -\tilde{\alpha}$. In addition, since $\Phi^{-1}\Phi^-$ is a diagonal matrix and $N_K (A_K-I) N^{\top}_K$ is symmetric, it is well-known that $\Phi^{-1}\Phi^-N_K (A_K-I) N^{\top}_K $ has real eigenvalues, and
\begin{equation}
\label{Hcurl_lem_unisol_eq2}
\lambda_{\max}( \Phi^{-1}\Phi^-N_K (A_K-I) N^{\top}_K ) \le\lambda_{\max}( \Phi^{-1}\Phi^-) \lambda_{\max}(N_K (A_K-I) N^{\top}_K ) \le 1- \tilde{\alpha}.
\end{equation}
Therefore, the eigenvalues of $I + \Phi^{-1}\Phi^-  N_K (A_K-I) N^{\top}_K$ must be bounded below by $1-(1-\tilde{\alpha} ) =  \tilde{\alpha}$, and thus, from \eqref{Hcurl_lem_unisol_eq1}, the matrix $M_K$ is non-singular. 

If $K$ is a regular tetrahedron, such a simple proof is not available. Instead, we have to employ a computer-aided argument to verify that $M_K$ is non-singular. As it is technical, we put it in Appendix~\ref{append:unisol}.
\end{proof}

Lemma~\ref{lem_Hdiv_unisol} guarantees the existence of face shape functions $\bfphi^f_i$, $i=1,2,3,,4$ satisfying
\begin{equation}
\label{Hdiv_IFE_shapefun}
\int_{F_j}\bfphi^f_i\cdot\bfn \dd s = \delta_{ij},  ~~~~ j=1,2,3,4,
\end{equation}
regardless of interface location. Then the global $\bfH(\div)$ IFE space can be constructed as
\begin{equation}
\begin{split}
\label{Hdiv_IFE_globSapce}
\bfS^f_h = \{ \bfv_h\in L^2(\Omega): & \bfv_h|_K\in \bfS^f_h(K), ~\text{if} ~ K\in \mathcal{T}^i_h, ~ \text{and} ~ \bfv_h|_K \in \mathcal{RT}_0(K), ~\text{if} ~ K\in \mathcal{T}^n_h,\\
 & \int_F [\bfv_h\cdot\bfn] \dd s = 0, ~ \forall F\in\mathcal{F}_h \}.
\end{split}
\end{equation}


With the nodal and face shape functions, we are ready to analyze the DoFs of $\bfS^e_h(K)$. 

\begin{lemma}
\label{Hcurl_lem_unisol}
Given a trirectangular tetrahedral or regular tetrahedral interface element $K$, the space $\bfS^e_h(K)$ has the DoFs $\int_{e_i} \bfv_h\cdot\bft_i \dd s$, $i=1,2,...,6$.
\end{lemma}
\begin{proof}
Given an interface element $K$ with $6$ edges $e_i$, $i=1,...,6$, we let $v_i=\int_{e_i} \bfv_h\cdot\bft_i \dd s$ for a function $\bfv_h\in\bfS^e_h(K)$ where $\bft_i$ have some assigned orientation. We need to show that $\bfv_h$ uniquely exists for arbitrary $v_i$. As the dimensions match, we only need to prove the existence which is done by explicit construction. For each face $F_j$ of $K$, we let 
\begin{equation}
\label{Hcurl_lem_unisol_eq_new2}
w_j = \sum_{e_i\subset \partial F_j} c_{j,e_i}v_i, ~~~ j=1,...,4, ~~~ \text{and} ~~~ \bfphi^f = \sum_{j=1}^4w_j \bfphi^f_j,
\end{equation}
where $c_{j,e_i} = 1$ or $-1$ is a sign to correct the orientation such that the edges associated with each face has the counter-clockwise orientation. So, trivially $\sum_{j=1}^4w_j=0$ and thus $\bfphi^f$ is $\div$-free, i.e., $\bfphi^f\in \curl \bfS^e_h(K)$ by the local exact sequence in \eqref{local_IFE_exatSeq}. We then proceed to construct the component in $S^n_h(K)$. Let the four nodes of $K$ be $\bfz_l$, $l=1,...,4$. Picking any node, say $\bfz_4$,  without loss of generality, we let the three neighbor edges be $e_i$ pointing to $\bfz_i$, $i=1,2,3$. Define the orientation coefficients $d_{i} = 1$ if $\bft_i$ points from $\bfz_4$ to another ending point of $\bfe_i$, otherwise $d_i = -1$. We then construct
\begin{equation}
\label{Hcurl_lem_unisol_eq_new3}
\phi^n = \sum_{i=1,2,3} d_{i} \left(v_i - \int_{e_i} (\bfphi^f \times (\bfx - \bfx_K))\cdot\bft_i/2 \dd s \right) \phi^n_i.
\end{equation}
Now, the edge function has the edge integral values $v_i$, $i=1,...,6$, can be constructed as
\begin{equation}
\label{Hcurl_lem_unisol_eq_new4}
\bfphi^e = \bfphi^f \times (\bfx - \bfx_K)/2 + \nabla \phi^n.
\end{equation}
It is easy to see that $\bfphi^e$ matches the DoFs on $e_i$, $i=1,2,3$. The verification of the rest edges can be shown by the DoFs of $\bfphi^f$ and $\phi^n$ with the formula $\int_F\curl\bfphi^e\cdot\bfn \dd s = \int_F \rot_F\bfphi^e \dd s = \int_{\partial F} \bfphi^e\cdot\bft \dd s$ for each face $F$.
\end{proof}

Thanks to the same DoFs, the whole construction procedure exactly mimics the standard FE spaces. Lemma~\ref{Hcurl_lem_unisol} guarantees the existence of N\'ed\'elec IFE shape functions $\bfphi^e_i$, $i=1,2,...,6$, satisfying
\begin{equation}
\label{Hcurl_IFEshapeFun} 
\int_{e_j} \bfphi^e_i\cdot\bft \dd s = \delta_{ij}, ~~~ j = 1,2,...,6.
\end{equation}
Then, the global $\bfH(\curl)$ IFE space can be constructed as
\begin{equation}
\begin{split}
\label{Hcurl_IFE_globSapce}
\bfS^e_h = \{ \bfv_h\in L^2(\Omega): &~ \bfv_h|_K \in \bfS^e_h(K), ~\text{if} ~ K\in \mathcal{T}^i_h, ~ \text{and} ~ \bfv_h|_K \in \mathcal{ND}_0(K), ~\text{if} ~ K\in \mathcal{T}^n_h,\\
 & \int_e (\bfv_h|_{K_1} - \bfv_h|_{K_2})\cdot\bft \dd s = 0, ~ \forall e\in\mathcal{E}_h, ~ \forall K_1,K_2 ~ \text{sharing} ~e \}.
\end{split}
\end{equation}

%
%

\begin{figure}[h]
\centering
\begin{subfigure}{.31\textwidth}
    \includegraphics[width=0.9in]{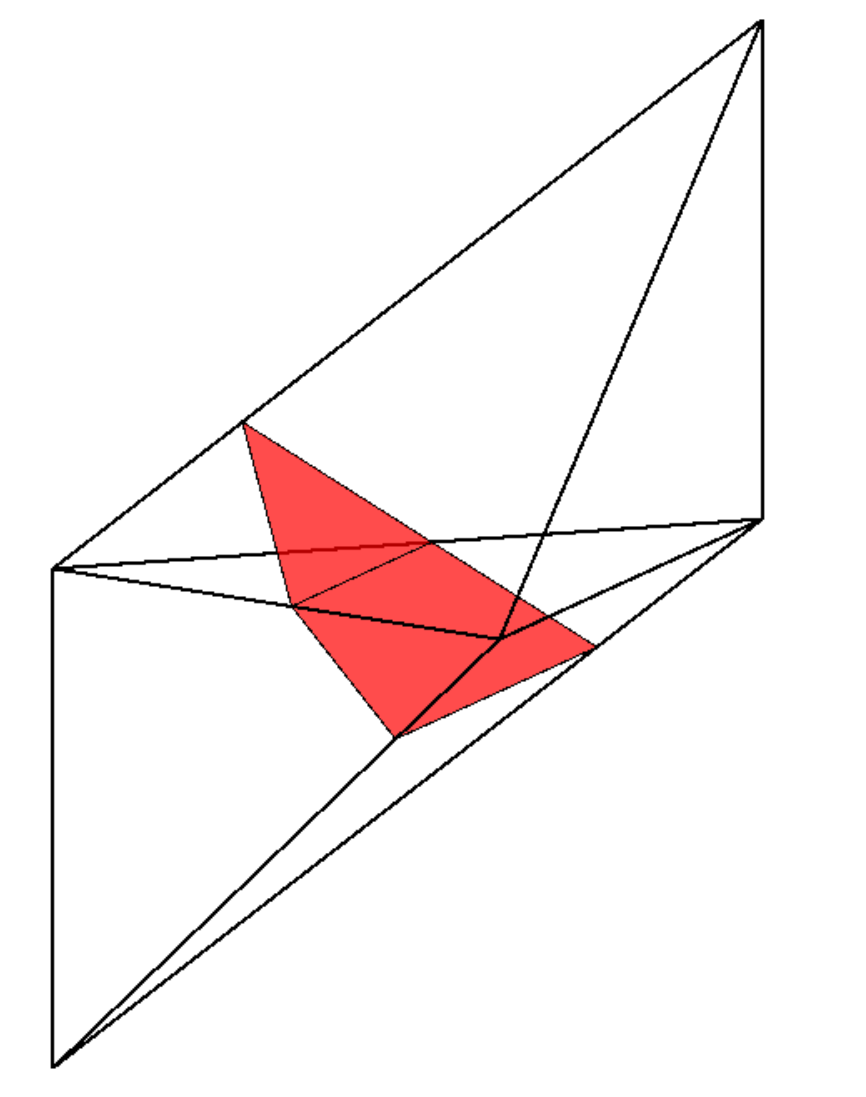}
\end{subfigure}
\begin{subfigure}{.4\textwidth}
     \includegraphics[width=2.2in]{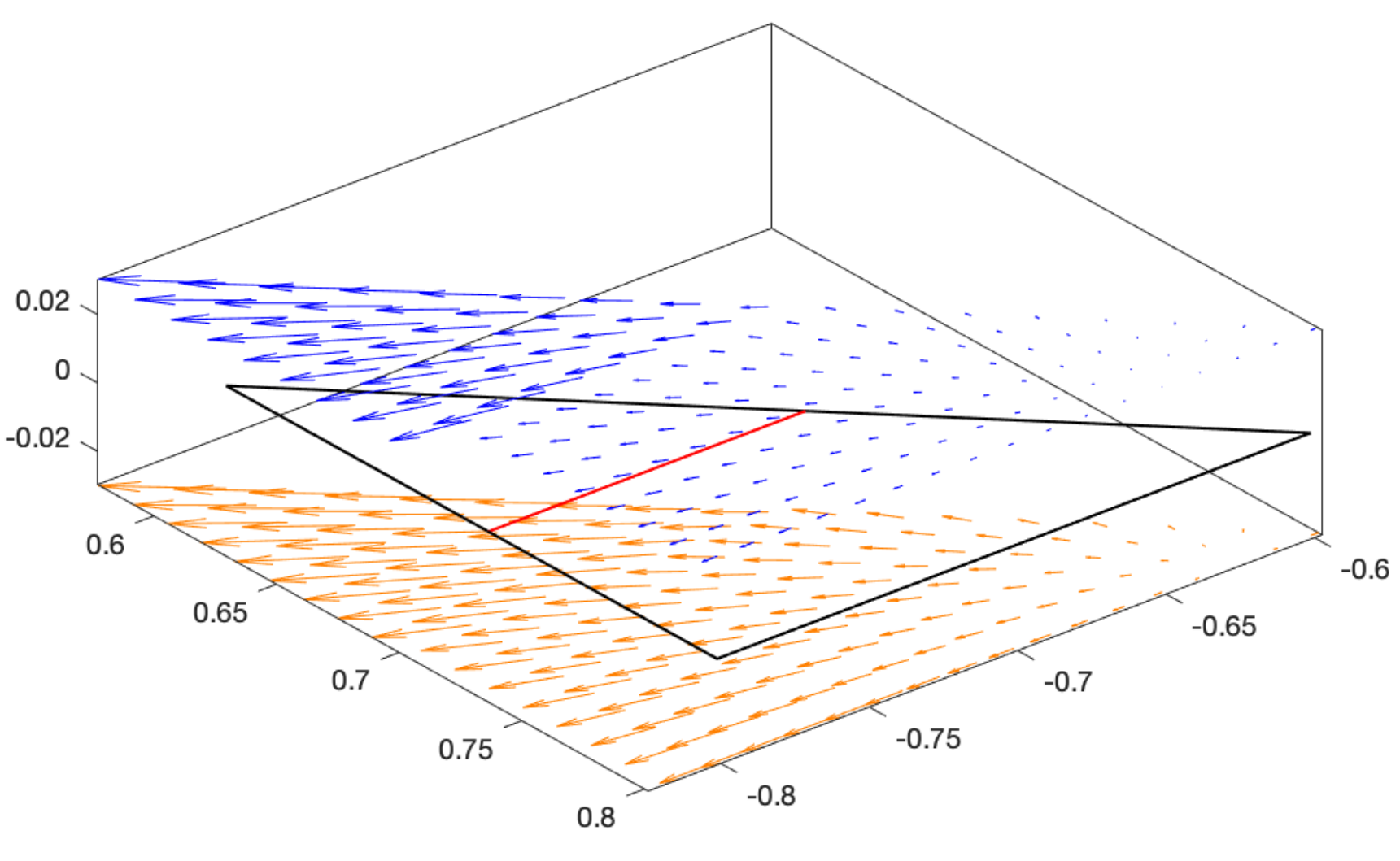}
\end{subfigure}
     \caption{Left: two neighbor elements. Right: the tangential traces of an IFE function on the sharing face from the two neighbor elements where the blue field is from the top and the orange field is from the bottom.}
  \label{fig:2tettrace} 
\end{figure}

\begin{remark}
The integral conditions in both \eqref{Hcurl_IFE_globSapce} and \eqref{Hdiv_IFE_globSapce} simply imply the tangential and normal continuity on non-interface edges and faces, which is the same as the standard FE spaces. 
However, the difference is on interface edges and faces where the integral conditions can not imply the exact continuity. Take the $\bfH(\curl)$ case as example: $\int_e(\bfv_h|_{K_1}-\bfv_h|_{K_2})\cdot\bft \dd s$ can not imply $(\bfv_h|_{K_1}-\bfv_h|_{K_2})\cdot\bft =0$ since $\bfv_h\cdot\bft$ is piecewise constant on $e^{\pm}_h$. 
So the DoFs should be only understood as continuity in certain average sense. We illustrate this kind of discontinuity in Figure~\ref{fig:2tettrace} for an edge IFE function $\bfv_h$ where we show the tangential traces $\bfv_h|_{K_i}\times \bfn$ (right) from the two neighbor elements, $K_i$, $i=1,2$, (left) on their sharing the face $F$ (right). Furthermore, for the $\bfH(\curl)$ case, 
the integral conditions $\int_e (\bfv_h|_{K_1}-\bfv_h|_{K_2})\cdot\bft \dd s=0$, $\forall e\subset \partial F$ can not even yield $\int_F [\bfv_h\times \bfn] \dd s = \mathbf{ 0}$, which is another difference from the standard FE space.
\end{remark}

Finally, the node, edge and face DoFs make the following interpolation operators well-defined:
\begin{subequations}
\label{interpolation}
\begin{align}
    I^n_h& : H^2(\beta;\Omega) \rightarrow S^n_h, ~~~~~~~~~~~~ \text{satisfying} ~~~ I^n_hu(\bfx) = u(\bfx), ~ \forall \bfx\in \mathcal{N}_h, \label{H1_IFE_interp}   \\
    I^e_h&: \bfH^1(\alpha,\beta,\curl;\Omega) \rightarrow \bfS^e_h, ~~~ \text{satisfying} ~ \int_e I^e_h\bfu\cdot\bft \dd s = \int_e \bfu\cdot\bft \dd s, ~ \forall e\in \mathcal{E}_h, \label{Hcurl_IFE_interp} \\
    I^f_h&: \bfH^1(\alpha,\div;\Omega) \rightarrow \bfS^f_h, ~~~~~~~ \text{satisfying} ~ \int_F I^f_h\bfu\cdot\bfn \dd s = \int_F \bfu\cdot\bfn \dd s, ~ \forall F\in \mathcal{F}_h. \label{Hdiv_IFE_interp}
\end{align}
\end{subequations}


\subsection{The de Rham Complex}
\label{subsec:complex}

Now, we recover the well-known and important de Rham complex in Diagram \eqref{IFE_deRham} for all the IFE spaces introduced above. Here, $Q_h$ denotes the piecewise constant space on the mesh $\mathcal{T}_h$, and $\Pi^0$ denotes the $L^2$ projection to $Q_h$. 

\begin{theorem}
\label{thm_deRham_IFE}
The sequence in the middle of \eqref{IFE_deRham} is a complex.
\end{theorem}
\begin{proof}
First of all, it is clear that $\mathbb{R}\subset S^n_h$. Besides, with the local result in \eqref{local_IFE_exatSeq}, we know that the functions mapped by the operators $\nabla$, $\curl$ and $\div$ must satisfy the corresponding jump conditions on each element. Thus, we only need to verify the corresponding integral conditions on edges and faces in \eqref{Hcurl_IFE_globSapce} and \eqref{Hdiv_IFE_globSapce}, respectively. For each edge $e\in\mathcal{E}_h$ with the nodes $\bfz_1$ and $\bfz_2$, we let $K_1$ and $K_2$ be any two elements sharing $e$. Then, given any $v_h\in S^n_h$, the continuity on nodes implies
\begin{equation}
\label{thm_deRham_IFE_eq1}
\int_e \nabla v_h|_{K_1} \cdot \bft \dd s = v_h|_{K_1}(\bfz_2) - v_h|_{K_1}(\bfz_1) = v_h|_{K_2}(\bfz_2) - v_h|_{K_2}(\bfz_1)  = \int_e \nabla v_h|_{K_2} \cdot \bft \dd s,
\end{equation}
where $\bft$ is assumed from $\bfz_1$ to $\bfz_2$. So, we conclude $\nabla v_h \in \bfS^e_h$. Next, for each face $F\in \mathcal{F}_h$, we let $K_1$ and $K_2$ be the two neighbor elements. Given any $\bfv_h\in \bfS^e_h$, we have 
\begin{equation}
\label{thm_deRham_IFE_eq2}
\int_F \curl \bfv_h|_{K_1}\cdot\bfn_F \dd s = \int_F \rot_F \bfv_h|_{K_1} \dd s = \int_{\partial F} \bfv_h|_{K_1}\cdot\bft \dd s = \int_{\partial F}  \bfv_h|_{K_2}\cdot\bft \dd s = \int_F \curl \bfv_h|_{K_2}\cdot\bfn_F \dd s.
\end{equation}
Therefore, we have $\curl \bfv_h \in \bfS^f_h$. Finally, it is trivial that $\div(\bfv_h)\in Q_h$ for any $\bfv_h\in \bfS^f_h$.
\end{proof}

\begin{theorem}
\label{thm_deRham_IFE_commut}
The diagram formed by the sequences in the middle and top of \eqref{IFE_deRham} is commutative. 
\end{theorem}
\begin{proof}
First of all, given every $u\in H^2(\beta;\Omega)$, there holds for each $e\in\mathcal{E}_h$ with the nodes $\bfz_1$ and $\bfz_2$ that
\begin{equation}
\begin{split}
\label{thm_deRham_IFE_commut_eq1}
\int_e I^e_h (\nabla u) \cdot \bft \dd s = \int_e \nabla u \cdot \bft \dd s = u(\bfz_2) - u(\bfz_1) = I^n_h u(\bfz_2) - I^n_h u(\bfz_1) =  \int_e \nabla I^n_h u \cdot \bft \dd s
\end{split}
\end{equation}
which implies $I^e_h \circ \nabla = \nabla \circ I^n_h$. Next, given $\bfu \in \bfH^1(\alpha,\beta,\curl;\Omega)$, there holds
\begin{equation}
\label{thm_deRham_IFE_commut_eq2}
\int_F I^f_h \curl \bfu \cdot \bfn \dd s = \int_{ F} \curl \bfu \cdot \bfn \dd s  = \int_{\partial F} \bfu\cdot\bft \dd s = \int_{\partial F} I^e_h\bfu\cdot\bft \dd s = \int_{ F} \curl I^e_h\bfu\cdot\bfn \dd s
\end{equation}
which shows $I^f_h\circ\curl = \curl \circ I^e_h$. Finally, for every $\bfu\in\bfH^1(\alpha,\div;\Omega)$, integration by parts yields
\begin{equation}
\label{thm_deRham_IFE_commut_eq3}
\int_K \Pi^0_K \div\bfu \dd \bfx = \int_K  \div\bfu \dd \bfx = \int_{\partial K} \bfu\cdot\bfn \dd s = \int_{\partial K} I^f_n\bfu\cdot\bfn \dd s =  \int_{ K} \div I^f_h \bfu \dd s
\end{equation}
which yields $\Pi^0\circ \div = \div \circ I^f_h$. Then, the commutativity has been proved.
\end{proof}

To show the exactness, let us first describe the isomorphism between the IFE spaces and standard FE spaces. Let $\widetilde{S}^n_h$, $\widetilde{\bfS}^e_h$ and $\widetilde{\bfS}^f_h$ be the standard Lagrange, N\'ed\'elec and Raviart-Thomas spaces defined on the same mesh $\mathcal{T}_h$. Note that the interpolations $I^n_h$, $I^e_h$ and $I^f_h$ are all well-defined on these standard spaces, which lead to certain isomorphism; namely their inversions from the IFE spaces to FE spaces are also well-defined. We shall keep the same notation $I^n_h$, $I^e_h$ and $I^f_h$ as the isomorphism, and illustrate this structure by the diagram formed from the middle and bottom sequences in \eqref{IFE_deRham},
where $I$ is the identity operator. 

\begin{lemma}
\label{lem_IFE_commut}
The diagram formed by the middle and bottom sequences in \eqref{IFE_deRham} is commutative.
\end{lemma}
\begin{proof}
The argument is the same as Theorem~\ref{thm_deRham_IFE_commut} by verifying the DoFs.
\end{proof}

\begin{lemma}
\label{lem_IFE_ker}
The operators $I^n_h$, $I^e_h$ and $I^f_h$ and their inversions preserve the kernels of $\nabla$, $\curl$ and $\div$: 
\begin{subequations}
\label{lem_IFE_ker_eq0}
\begin{align}
& \ker(\nabla)\cap \widetilde{S}^n_h \xrightleftarrows[(I^n_h)^{-1}]{~~~I^n_h~~~}     \ker(\nabla)\cap S^n_h, \label{lem_IFE_ker_eq01} \\ 
& \ker(\curl)\cap \widetilde{\bfS}^e_h \xrightleftarrows[(I^e_h)^{-1}]{~~~I^e_h~~~} \ker(\curl)\cap \bfS^e_h,  \label{lem_IFE_ker_eq02}\\
& \ker(\div)\cap \widetilde{\bfS}^f_h \xrightleftarrows[(I^f_h)^{-1}]{~~~I^f_h~~~} \ker(\div)\cap \bfS^f_h .  \label{lem_IFE_ker_eq03}
\end{align}
\end{subequations}
\end{lemma}
\begin{proof}
\eqref{lem_IFE_ker_eq01} is trivial since $\ker(\nabla)\cap \widetilde{S}^n_h= \ker(\nabla)\cap S^n_h = \mathbb{P}_0(K)$ and $I^n_h$ preserves constants. For \eqref{lem_IFE_ker_eq02}, thanks to the local exact sequence, we know $\bfv_h\in \ker(\curl)\cap \widetilde{\bfS}^e_h$ if and only if
$
\int_F \curl\bfv_h\cdot\bfn \dd s  = \int_{\partial F} \bfv_h\cdot\bft \dd s= 0,
$
for each face $F\in\mathcal{F}_K$, which is equivalent to $\int_{\partial F} I^e_h\bfv_h\cdot\bft \dd s = \int_F \curl I^e_h \bfv_h\cdot\bfn \dd s = 0$. Then, Lemma~\ref{lem_Hdiv_unisol} yields the desired result. The argument for \eqref{lem_IFE_ker_eq03} and its inverse direction is similar. 
\end{proof}

\begin{theorem}
\label{thm_exact_IFE}
The sequence on the middle of \eqref{IFE_deRham} is exact.
\end{theorem}
\begin{proof}
The result directly follows from Lemma~\ref{lem_IFE_ker} and the exactness for the standard FE spaces.
\end{proof}




\section{The Petrov-Galerkin IFE Method for The $\bfH(\curl)$ Interface Problem}
\label{sec:formulation}

In this section, we apply the IFE spaces, particularly the $\bfH(\curl)$ IFE space, to the $H(\curl)$ interface problem. 

\subsection{Schemes}
For simplicity of presentation, we here only consider the perfect electric conductor (PEC) or correspondingly the perfect magnetic conductor (PMC) boundary condition, i.e., $\bfu\times \bfn = \mathbf{ 0}$ on $\partial \Omega$. For the quasi-static equation \eqref{Max_mode1},
the proposed PG-IFE scheme is to find $\bfu_h\in \bfS^e_{h,0}$, i.e., the subspace of $\bfS^e_{h}$ with the zero trace, such that
\begin{equation}
\label{Max_model3_wf2}
a_h(\bfu_h,\bfv_h) = (\bff,\bfv)_{\Omega}, ~~ \forall \bfv_h \in \widetilde{\bfS}^e_{h,0}~~~~~ \text{with} ~~ a_h(\bfu_h,\bfv_h) : = (\alpha_h\curl\bfu_h, \curl \bfv_h)_{\Omega} + (\beta_h \bfu_h, \bfv_h)_{\Omega},
\end{equation}
where $(\cdot,\cdot)_{\Omega}$ denotes the standard $L^2$ inner product. Note that \eqref{Max_model3_wf2} follows from the integration by parts of \eqref{Max_model3} which always holds as long as the test function space is $\bfH(\curl)$-conforming.
Namely, even without any penalties, the scheme \eqref{Max_model3_wf2} is still consistent. Thanks to the isomorphism between the IFE and FE spaces, the resulting linear system is square. So the PG formulation completely avoids the non-conformity issue aforementioned in the introduction.


\subsection{Fast Solvers: HX Preconditioner}
\label{subsec:solver}

{The resulting linear system of \eqref{Max_model3_wf2} is denoted by 
\begin{equation}
\label{PG_IFE_system}
\bfB^{PG} \bar{\bfu} = \bar{\bff}
\end{equation}
where $\bfB^{PG}$ is not symmetric positive definite (SPD).} However, the majority portion of $\bfB^{PG}$ is indeed SPD, and the ``problematic" portion is only around the interface. So, naturally, it may help if we employ the exact solver on the subsystem around the interface and use an iterative solver on the rest part of the system, which motivates the following fast solver. Define the collections of elements and edges around the interface:
\begin{equation}
\begin{split}
\label{InterEdgeSet}
 \mathcal{T}^l_h = \{ K \in \mathcal{T}_h: \overline{K}\cap\overline{K'}\neq\emptyset, ~ \text{for some } K'\in \mathcal{T}^{l-1}_h \}, ~~~  \mathcal{E}^l_h = \{ e \in \mathcal{E}_h: e\in\mathcal{E}_K ~\text{for some }K \in \mathcal{T}^l_h \},
\end{split}
\end{equation}
where $ \mathcal{T}^0_h= \mathcal{T}^i_h$ is simply the set of interface elements. It is not hard to see that $\mathcal{T}_h^l$ and $\mathcal{E}^l_h$ can expand to the whole domain from the interface when we increase $l$. But with relatively small $l$, the set $\mathcal{E}^l_h$ is sufficient to cover the ``problematic" DoFs. Assuming the DoFs associated with $\mathcal{E}^l_h$ are indexed at the end, we can write
 \begin{equation}
\label{PG_IFE_system_mat}
\bfB^{PG} = 
\left[\begin{array}{cc} \bfB^{(1,1)}_l & \bfB^{(2,1)}_l  \\ \bfB^{(1,2)}_l & \bfB^{(2,2)}_l \end{array}\right],
\end{equation}
where $\bfB^{(1,1)}_l$ is SPD but $\bfB^{(2,2)}_l$ is non-SPD and also $\bfB^{(1,2)}_l\neq (\bfB^{(2,1)}_l)^{\top}$.

We write $\bfB^{\nabla}$ as the matrix associated with the IFE discretization of the bilinear form
\begin{equation}
\label{solver_eq1}
a^{\nabla}_{h,1}(\bfu_h,\bfv_h) := (\alpha_h\nabla \bfu_h, \nabla \bfv_h)_{\Omega} + (\beta_h \bfu_h, \bfv_h)_{\Omega}, ~~~~ \forall \bfu_h,\bfv_h \in [S^n_{h,0}]^3,
\end{equation}
and write $B^{\nabla}$ as the matrix of the bilinear form
\begin{equation}
\label{solver_eq2}
a^{\nabla}_{h,2}(u_h,v_h) := (\alpha_h\nabla u_h, \nabla v_h)_{\Omega} , ~~~~ \forall u_h,v_h \in S^n_{h,0}.
\end{equation}
{Based on \eqref{PG_IFE_system_mat}, we then construct a special smoother:
\begin{equation}
\label{smoother}
\bfB^{PG}_{l}=
\left[\begin{array}{cc} \text{diag}(\bfB^{(1,1)}_l) & 0 \\0 & \bfB^{(2,2)}_l \end{array}\right],
\end{equation}
In the following discussion, $\bfB^{PG}_l$ and $l$ shall be frequently referred to as a ``block diagonal smoother" and an ``expanding width". 
}
 Then, the IFE version of the HX preconditioner~\cite{2007HiptmairXu} is
\begin{equation}
\label{solver_eq3}
\bfB^{\curl} = [\bfB^{PG}_{l}]^{-1} + P_{\curl} (\bfB^{\nabla})^{-1}  P^{\top}_{\curl} + G (B^{\nabla})^{-1} G^{\top}
\end{equation}
where $G$ is the matrix associated with $\nabla:\widetilde{S}^n_h\rightarrow\widetilde{\bfS}^e_h$, and $P_{\curl}$ is the matrix associated with the operator $\bfPi^{\curl}_h:[\widetilde{S}^n_h]\rightarrow \widetilde{\bfS}^e_h$ such that $\bfPi^{\curl}_h\bfw|_e\cdot\bft = (\bfw_1 + \bfw_2)\cdot\bft$, $\forall e\in \mathcal{E}_h$ with $\bfw_1$ and $\bfw_2$ being the nodal values of $\bfw$ at the nodes of $e$. We emphasize that forming such a matrix highly relies on the isomorphism between IFE and FE spaces through the DoFs. The matrices $\bfB^{\nabla}$ and $B^{\nabla}$ are SPDs and their inverses are approximated by a few algebraic multigrid (AMG) cycles. {The $[\bfB^{PG}_{l}]^{-1}$, particularly the block $(\bfB^{(2,2)}_l)^{-1}$}, is computed by the direct solver (the backslash in MATLAB) as it is of small size, and we shall see that this direct solver only has a very subtle effect on the total computational time but is the key to handle the ``problematic" DoFs near the interface. We adopt the HX preconditioner implemented in $i$FEM~\cite{Chen:2008ifem}.



\subsection{Numerical Test of the \textit{Inf-Sup} Condition}
Define the natural energy norm 
$$
\vertii{\bfw_h}^2_{h,D} = \int_{D}\alpha_h \curl\bfw_h\cdot\curl\bfw_h \dd \bfx +  \int_{D}\beta_h \bfw_h\cdot\bfw_h \dd \bfx,~~~~\forall \bfw_h\in \bfS^e_h~ \text{or} ~\widetilde{\bfS}^e_h.
 $$
 Clearly, there holds the identity $\|\bfv_h\|_{h,\Omega}^2 = a_h(\bfv_h,\bfv_h)$ and the norm equivalence $\|\cdot\|_{\bfH(\curl;\Omega)}\simeq \vertii{\cdot}_{h,\Omega}$. With the optimal approximation capabilities of the $\bfH(\curl)$ IFE spaces discussed later, one can immediately show the optimal convergence for the PG-IFE scheme as long as the following \textit{inf-sup} condition holds: 
\begin{equation}
\label{infsup_cond}
\eta_s := \sup_{\bfv_h\in \widetilde{\bfS}^e_{h,0}}\frac{ a_h(\bfu_h,\bfv_h) }{ \| \bfv_h \|_{h,\Omega} \| \bfu_h \|_{h,\Omega} } \ge C_s, ~~~~ \forall \bfu_h \in \bfS^e_{h,0}.
\end{equation}
However, it is generally quite challenging to show such a \textit{inf-sup} condition. The approach of Babu\v{s}ka in~\cite{1994BabuskaCalozOsborn} assumes a certain structure of the elliptic coefficient. The proof of the 2D PG-IFE method~\cite{2020GuoLinZou} highly relies on the symmetry of the high-order $\curl$-$\curl$ term, which does not hold anymore for the 3D case. In this work, instead of pursuing a rigorous theoretical analysis, we provide a numerical test for the \textit{inf-sup} condition \eqref{infsup_cond}. 

Denote $\bfB^{FE}$ and $\bfB^{IF}$ as the matrices corresponding to the FE and IFE schemes on the same mesh, i.e., the matrices associated with the bilinear forms $a_h(\bfw_h,\bfv_h)$, $\bfw_h,\bfv_h \in \widetilde{\bfS}^e_{h,0}$ and $a_h(\bfw_h,\bfv_h)$, $\bfw_h,\bfv_h \in \bfS^e_{h,0}$, respectively. Further let $N$ be the total number of edge DoFs. Then, \eqref{infsup_cond} is equivalent to 
\begin{equation}
\label{infsup_cond_1}
\sup_{\bar{\bfv}\in \mathbb{R}^N} \frac{ \bar{\bfv}^{\top}\bfB^{PG}\bar{\bfu}  }{ \sqrt{ \bar{\bfu}^{\top}\bfB^{IF}\bar{\bfu} }  \sqrt{ \bar{\bfv}^{\top}\bfB^{FE}\bar{\bfv} } } \ge C_s, ~~~~ \forall \bar{\bfu}\in \mathbb{R}^N.
\end{equation}
In general, 
one may verify \eqref{infsup_cond_1} by constructing a suitable test function $\bar{\bfv} = \mathbf{P}\bar{\bfu}$, where $\mathbf{P}$ is a non-singular matrix, such that
\begin{equation}
\label{infsup_cond_4}
\frac{ \bar{\bfu}^{\top} \mathbf{P}^{\top} \bfB^{PG}\bar{\bfu}  }{ \sqrt{ \bar{\bfu}^{\top}\bfB^{IF}\bar{\bfu} }  \sqrt{ \bar{\bfu}^{\top} \mathbf{ P}^{\top} \bfB^{FE} \mathbf{ P} \bar{\bfu} } } \ge C_s, ~~~~ \forall \bar{\bfu}\in \mathbb{R}^N.
\end{equation}
It is not hard to see that \eqref{infsup_cond_4} implies \eqref{infsup_cond_3} as it takes sup over $\bar{\bfv}\in \mathbb{R}^N$. We shall concentrate on \eqref{infsup_cond_4}.

Note that one natural choice of the test function $\bfv_h = (I^e_h)^{-1}\bfu_h$. Consequently, in \eqref{infsup_cond_4} $\bar{\bfu}=\bar{\bfv}$, i.e., $\mathbf{ P}$ is an identity matrix, and thus an equivalent condition for \eqref{infsup_cond_4} to hold is $\bfB^{PG} + (\bfB^{PG})^{\top}$ being SPD. Unfortunately, numerical results show that $\bfB^{PG} + (\bfB^{PG})^{\top}$ has negative eigenvalues. Our experience suggests that such a local construction approach may not be appropriate. Hence, we consider a global construction approach. Given each IFE function $\bfu_h\in \bfS^e_{h,0}$, we define its projection to the FE space denoted by $\pi_h\bfu_h\in \widetilde{\bfS}^e_{h,0}$ through $a_h(\cdot,\cdot)$, namely
\begin{equation}
\label{IFE2FEproj}
a_h(\bfu_h, \bfv_h) = a_h(\pi_h\bfu_h, \bfv_h), ~~~~ \forall \bfv_h\in \widetilde{\bfS}^e_{h,0}. 
\end{equation}
Note that \eqref{IFE2FEproj} is always well-defined since $a_h$ is SPD on the FE space. 
The corresponding vector solution of $\pi_h\bfu_h$ is $\bar{\bfv} = (\bfB^{FE})^{-1} \bfB^{PG} \bar{\bfu}$, i.e., $\mathbf{ P} = (\bfB^{FE})^{-1} \bfB^{PG}$. Notice that
\begin{equation}
\label{infsup_cond_3}
 \bar{\bfu}^{\top} \mathbf{P}^{\top} \bfB^{PG}\bar{\bfu} = \bar\bfu^{\top} (\bfB^{PG})^{\top} (\bfB^{FE})^{-1} \bfB^{PG} \bar{\bfu} = \bar{\bfu}^{\top} \mathbf{ P}^{\top} \bfB^{FE} \mathbf{ P} \bar{\bfu} .
\end{equation}
Thus, \eqref{infsup_cond_4} becomes the following minimization problem
\begin{equation}
\label{infsup_cond_2}
\eta_s = \min_{\bar{\bfu}\in \mathbb{R}^N} \frac{ \sqrt{ \bar{\bfu}^{\top} (\bfB^{PG})^{\top} (\bfB^{FE})^{-1} \bfB^{PG} \bar{\bfu} } }{ \sqrt{ \bar{\bfu}^{\top}\bfB^{IF}\bar{\bfu} } } .
\end{equation}
Similar to~\cite{2001Bathe,1993ChapelleBathe}, the solution to the minimization problem in \eqref{infsup_cond_2} is the smallest generalized eigenvalue:
\begin{equation}
\label{geigen}
 (\bfB^{PG})^{\top} (\bfB^{FE})^{-1} \bfB^{PG} \bar{\bfw} = \lambda \bfB^{IF} \bar{\bfw}.
\end{equation}
Here, $\lambda$ must be non-negative as the matrices on the left is semi-SPD and the right one is SPD. Thus, we have $\eta_s = \sqrt{\lambda_{\min}}$ which can be used as a numerical estimate of the \textit{inf-sup} constant. We refer readers to Section~\ref{sec:numer} for the numerical results.

\section{Approximation Capabilities}
\label{sec:approx}

In this section, we proceed to show the optimal approximation capabilities of the IFE spaces. Due to the interface conditions, it will be very different from the estimation of the interpolation for the standard FE spaces in~\cite{1999CiarletZou,2003Monk}. Before the technical details, let us describe some major challenges and the fundamental idea, which also illustrates the structure of this section. The first one concerns the mismatching portion sandwiched by $\Gamma_h$ and $\Gamma$ around the interface.
To deal with this issue, we employ the $\delta$-strip argument~\cite{2012HiptmairLiZou,2010LiMelenkWohlmuthZou} that is widely used for interface problems, but special attention should be paid to the low regularity of functions in $\bfH^1(\alpha,\beta,\curl;\Omega)$. This result together with some other fundamental estimates will be presented in Section~\ref{subsec:fundanalysis}. The second one is the non-smoothness of the IFE functions across the interface which makes classical techniques unapplicable here. For this issue, a special quasi-interpolation will be constructed to encode the jump information that is discussed in Section~\ref{subsec:quasiInterp}. In Section~\ref{subsec:faceedgeinterp}, the quasi-interpolation will be used as a bridge to estimate the face and edge interpolations.



\subsection{Some Fundamental Estimates}
\label{subsec:fundanalysis}
We first define the $\delta$-strip $S_{\delta} = \{ \bfx: \text{dist}(\bfx,\Gamma) \le \delta \}$.
Note that Lemma~\ref{lem_geo_gamma} immediately implies that
\begin{equation}
\label{delta_strip_1}
\Omega_{\textint} \subset S_{\delta}, ~~~~ \text{with} ~~ 0<\delta \le c_0 h^2.
\end{equation}
The key is the following estimate that enables us to control the $L^2$-norm on the $\delta$-strip by its width.
\begin{lemma}
\label{lem_delta}
It holds for any $z\in H^1_0(\Omega^{\pm})$ that
\begin{equation}
\label{lem_delta_eq0}
\| z \|_{L^2(S_{\delta})} \le C \sqrt{\delta} \| z \|_{H^1(\Omega)}.
\end{equation}
\end{lemma}

Next, we recall the $\bfH^1(\text{curl};\Omega)$ and $\bfH^1(\text{div};\Omega)$-extension operators, see Theorem 3.4 and Corollary 3.5 in~\cite{2012HiptmairLiZou} and Theorem 3.2 in~\cite{2010HiptmairLiZou}. To simplify the discussion, we shall adopt the general differential operator $\dd = \curl$ or $\div$ widely used in differential forms.
\begin{thm}
\label{thm_ext}
For $\dd = \curl$ or $\div$, there exist two bounded linear operators
\begin{equation}
\label{thm_ext_eq0}
\bfE^{\pm}_{\dd} ~:~ \bfH^1(\dd;\Omega^{\pm})\rightarrow \bfH^1(\dd;\Omega),
\end{equation}
such that for each $\bfu\in\bfH^1(\dd;\Omega^{\pm})$:
\begin{itemize}
  \item[1.]  $\bfE^{\pm}_{\dd}\bfu = \bfu ~~ \text{a.e.} ~\text{in} ~ \Omega^{\pm}$,
  \smallskip
  \item[2.] $\| \bfE^{\pm}_{\dd}\bfu \|_{\bfH^1(\dd;\Omega)} \le C_E \| \bfu \|_{\bfH^1(\dd;\Omega^{\pm})}$ with the constant $C_E$ only depending on $\Omega$.
  \end{itemize}  
\end{thm}
With this theorem, we define $\bfu^{\pm}_{E,\dd} = \bfE^{\pm}_{\dd}\bfu^{\pm}$ which are the keys in the analysis. Furthermore, for any subdomain $D\subset \Omega$, we shall employ the following notation for simplicity:
\begin{equation}
\label{Enorm}
\| \bfu \|_{E,\dd,s,D} = \| \bfu^+_{E,\dd} \|_{\bfH^s(\dd,D)} + \| \bfu^-_{E,\dd} \|_{\bfH^s(\dd,D)},
\end{equation}
and $\| \bfu \|_{E,s,D}$ simply denotes the Sobolev norm without the $\curl$ or $\div$ terms, i.e., only the $L^2$ norm.

\begin{figure}[h]
\centering
\begin{minipage}{.29\textwidth}
  \centering
  \includegraphics[width=1.6in]{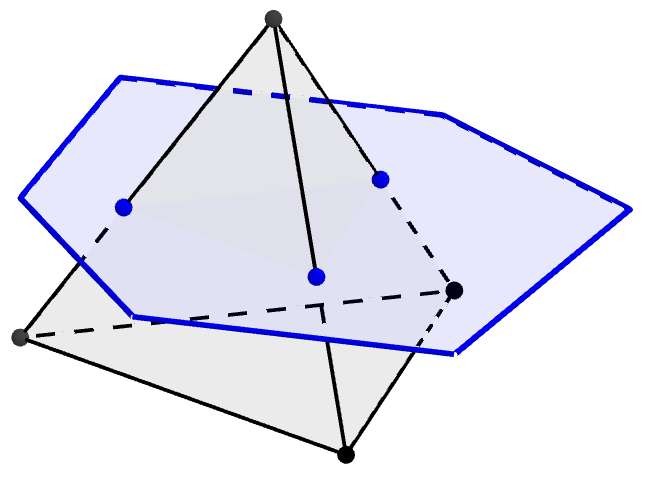}
  \caption{The approximate interface plane $\Gamma^{\omega_K}_h$.}
  \label{fig:patch}
\end{minipage}
\hspace{0.5cm}
\begin{minipage}{.3\textwidth}
  \centering
  \ \includegraphics[width=1.3in]{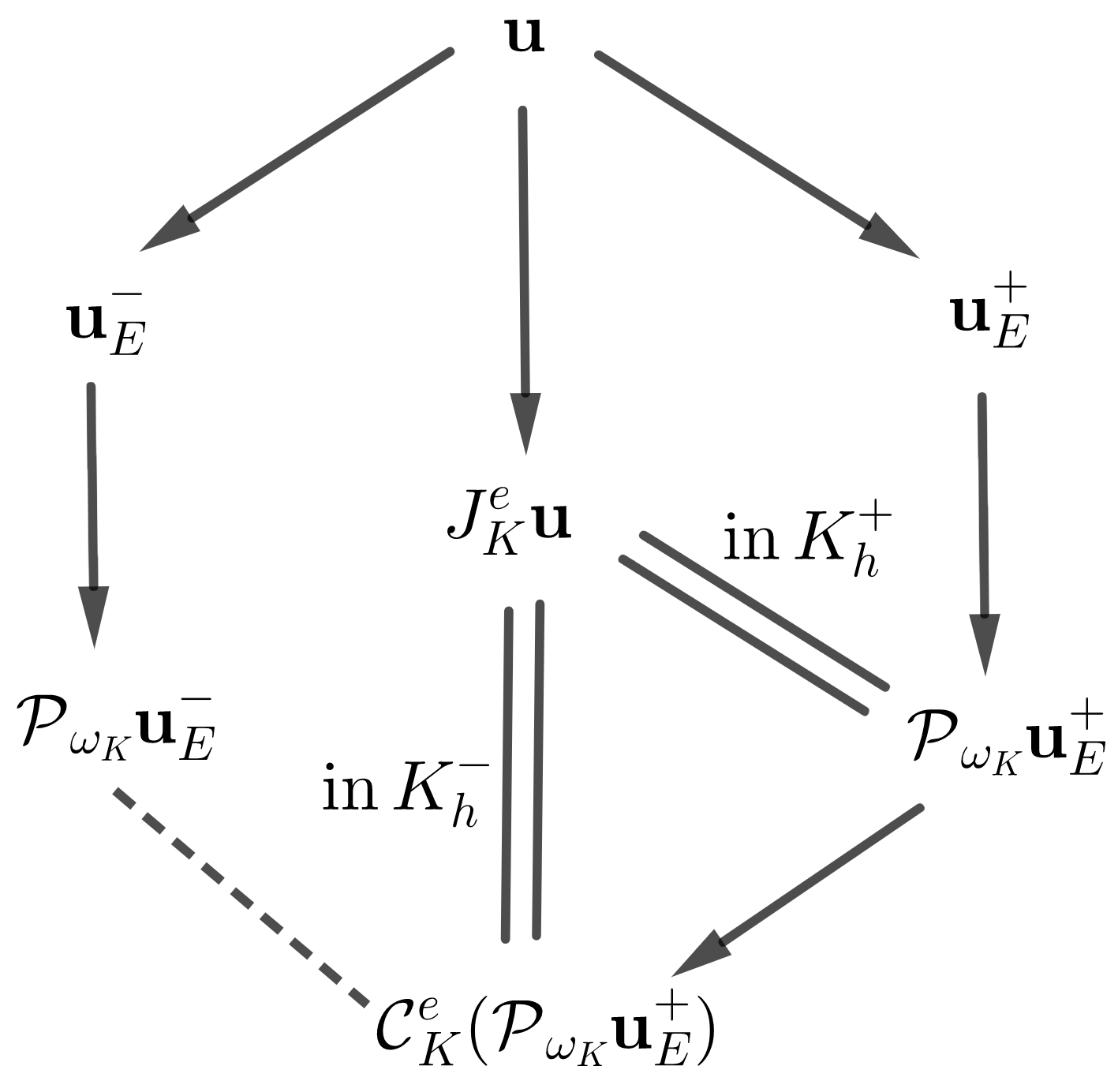}
  \caption{The diagram for interpolation.}
  \label{fig:Interpdiagram}
  \end{minipage}
  \hspace{0.5cm}
  \begin{minipage}{.3\textwidth}
  \centering
  \includegraphics[width=1.9in]{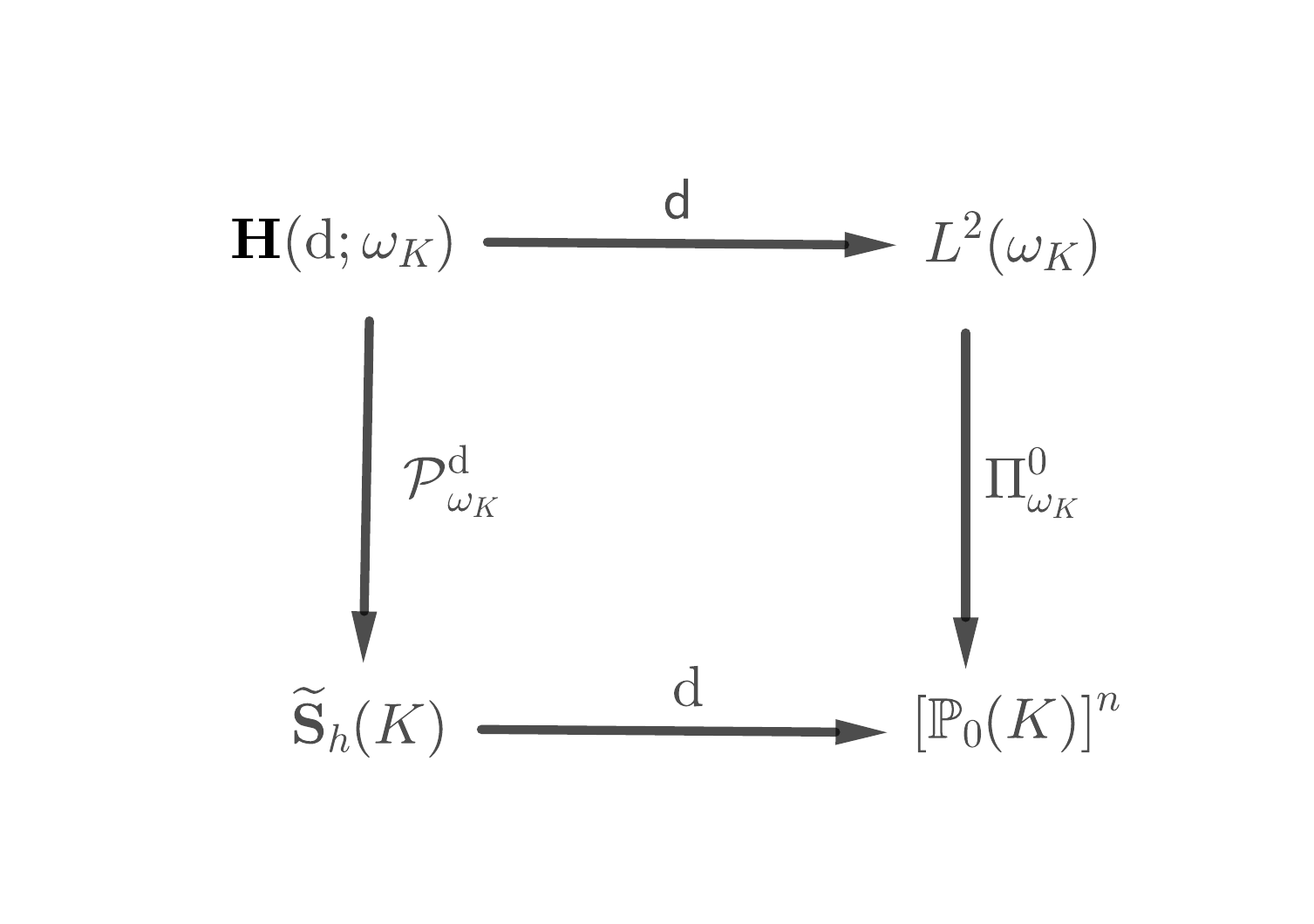}
  \caption{The commutative diagram for $\mathcal{P}^{\dd}_{\omega_K}$.}
  \label{fig:ppdiag}
  \end{minipage}
\end{figure}

Next, we introduce a patch $\omega_K$ of each interface element $K$: 
\begin{equation}
\label{patch_K}
\omega_K =  \cup\{ K': \partial K'\cap \partial K\neq\emptyset \}.
\end{equation}
Then, for each $\omega_K$, define $\Gamma^{\omega_K} := \Gamma\cap\omega_K$ and $\Gamma^{\omega_K}_h := \hat{\Gamma}^K_h\cap\omega_K$, where $\hat{\Gamma}^{K}_h$ is the entire plane containing $\Gamma^{K}_h$. Namely, $\Gamma^{\omega_K}_h$ is the extension of $\Gamma^{K}_h$ to the patch, see Figure~\ref{fig:patch} for an illustration.
Let $l^{\Gamma}_K$ be the maximal height within $\omega_K$ supporting $\Gamma^{\omega_K}_h$. Then, there exist $\gamma_1>0$ and $\gamma_2>0$ such that 
\begin{equation}
\label{geo_estimate_1}
| \Gamma^{\omega_K}_h | \ge \gamma_1 h^2_K ~~~~ \text{and}  ~~~~ l^{\Gamma}_K \ge \gamma_2 h_K;
\end{equation}
namely, $\Gamma^{\omega_K}_h$ as a polygon is always shape regular. The key of employing the patch is the following trace-type inequalities that have constants independent of interface location.
\begin{lemma}
\label{lem_patch_trace_inequa}
On each interface element $K$, for $v_h\in \mathbb{P}_1(K)$, there holds
\begin{equation}
\label{lem_patch_trace_inequa_eq}
|v_h(\bfz)| \lesssim h^{-1}_K \| v_h \|_{L^2(\Gamma^{\omega_K}_h)}, ~ \forall \bfz\in\Gamma^{\omega_K}_h, ~~~ \text{and} ~~~ \| v_h \|_{L^2(\Gamma^{\omega_K}_h)} \lesssim h^{-1/2}_K \| v_h \|_{L^2(\omega_K)}. 
\end{equation}
\end{lemma}
\begin{proof}
The first one is~\cite[Lemma 4.2]{2020GuoLin}. The second follows from~\cite[Lemma 4.1]{2020GuoLin} and an inverse inequality.
\end{proof}
Similar to Lemma~\ref{lemma_interface_flat}, the following geometric estimates regards the closeness between $\Gamma^{\omega_K}$ and $\Gamma^{\omega_K}_h$.
\begin{lemma}
\label{lem_geo_gamma}
Under the conditions of Lemma~\ref{lemma_interface_flat}, for each interface element $K$ and its patch $\omega_K$, there holds for each $\bfx\in\Gamma^{\omega_K}_h$ that
\begin{equation}
\label{lem_geo_gamma_eq0}
 \| \bfx - \bfx_{\bot} \| \le \tilde{c}_{\Gamma,1} h^2_K ~~~ \text{and} ~~~\| \bfn - \bar{\bfn} \| \le \tilde{c}_{\Gamma,2} h_K.
\end{equation}
\end{lemma}
\begin{proof}
The argument is similar to \eqref{lemma_interface_flat} by recognizing that \eqref{lemma_interface_flat_eq1} is actually true on $\omega_K$.
\end{proof}




We further need the following Sobolev-type inequalities relying on the closeness between $\Gamma^{\omega_K}$ and $\Gamma^{\omega_K}_h$.

\begin{lemma}
\label{lem_diff_u}
Given each interface element $K$, there holds
\begin{subequations}
\label{lem_diff_u_eq0}
\begin{align}
    &   \| a^- \bfu^-_{E,\dd}\times \bar{\bfn}_K - a^+ \bfu^+_{E,\dd}\times \bar{\bfn}_K \|_{L^2(\Gamma^{\omega_K}_h )} \lesssim h^{1/2}_K  \| \bfu \|_{E,1,\omega_K}  , \label{lem_function_gamma_eq01}  \\  
    &  \| b^- \bfu^-_{E,\dd}\cdot\bar{\bfn}_K - b^+ \bfu^+_{E,\dd}\cdot\bar{\bfn}_K \|_{L^2(\Gamma^{\omega_K}_h )} \lesssim h^{1/2}_K  \| \bfu \|_{E,1,\omega_K}  , \label{lem_function_gamma_eq02} 
\end{align}
where $\dd=\curl$ for $(a,b)=(1,\beta)$, $\bfu\in \bfH^1(\alpha,\beta,\curl;\Omega)$ or $\dd=\div$ for $(a,b)=(\alpha,1)$, $\bfu\in \bfH^1(\alpha,\div;\Omega)$. 
\end{subequations}
\end{lemma}
\begin{proof}
See Appedix~\ref{appen:sec:lem_diff_u}.
\end{proof}

\begin{remark}
\label{rem_diff_u}
As $\curl\bfH^1(\alpha,\beta,\Omega,\curl)\subset\bfH^1(\alpha,\div;\Omega)$, $\curl\bfu$ with $\bfu\in\bfH^1(\alpha,\beta,\curl,\Omega)$ actually satisfies \eqref{lem_diff_u_eq0} with $\dd = \div$. Furthermoer, with continuity of $\div\bfu$, it can be similarly shown that
\begin{equation}
\label{lem_function_gamma_eq003}
 \| \div \bfu^-_{E,\div} -  \div \bfu^+_{E,\div} \|_{L^2(\Gamma^{\omega_K}_h )} \lesssim h^{1/2}_K  \| \bfu \|_{E,\div,1,\omega_K} . 
\end{equation}
\end{remark}

\subsection{Special Quasi-Interpolation Operators}
\label{subsec:quasiInterp}
To develop quasi-interpolation for handling the discontinuity of IFE functions across the interface, 
by Table~\ref{tab_IFE_func} and \eqref{AB_space2}, let us rewrite the functions in $\bfS^e_h(K)$ in the following format
\begin{equation}
\label{IFEfun_ext}
\bfv_h=
\begin{cases}
      & \bfv^+_h = \bfa^+_h \times(\bfx-\bfx_K) + \bfb^+_h,~~~~~~~~~~~~~~~~~~~~~~~~~~~~~~~~~~ \text{in} ~ K^+_h, \\
      &  \bfv^-_h = \mathcal{C}^{\curl}_K(\bfv^+_h ) := (A_K\bfa^+_h) \times(\bfx-\bfx_K) + (B_K\bfb^+_h),~~~ \text{in} ~ K^-_h;
\end{cases}
\end{equation}
and similarly, the functions in $\bfS^f_h(K)$ can be represented as
\begin{equation}
\label{Hdiv_IFEfun_ext}
\bfv_h=
\begin{cases}
      & \bfv^+_h = c_0(\bfx-\bfx_K) + \bfa^+_h,~~~~~~~~~~~~~~~~~~~~~~~~~~~ \text{in} ~ K^+_h, \\
      &  \bfv^-_h = \mathcal{C}^{\div}_K(\bfv^+_h ) := c_0(\bfx-\bfx_K) + (A_K\bfa^+_h),~~~~ \text{in} ~ K^-_h.
\end{cases}
\end{equation}
The operators $\mathcal{C}^{\curl}_K$ and $\mathcal{C}^{\div}_K$ can be considered as discrete extensions mapping the polynomial component on the ``$+$" side to the ``$-$"  side through the jump conditions. Note that the ``$-$" piece of IFE functions can be completely determined by the ``$+$" piece. In addition, we  also need an interpolation: $\mathcal{P}^{\dd}_{\omega_K}:\bfH(\dd;\omega_K) \rightarrow \widetilde{\bfS}_h(K)$ where $\widetilde\bfS_h(K)=\mathcal{ND}_0(K)$ if $\dd = \curl$ or $\widetilde\bfS_h(K)=\mathcal{RT}_0(K)$ if $\dd=\div$ such that
\begin{equation}
\label{Interp_JK}
\mathcal{P}^{\curl}_{\omega_K} \bfv := (\Pi^0_{\omega_K} \curl\bfv) \times (\bfx - \bfx_K) / 2 + \Pi^0_{\omega_K} \bfv ~~ \text{and} ~~ \mathcal{P}^{\div}_{\omega_K} \bfv := (\Pi^0_{\omega_K} \div\bfv) (\bfx - \bfx_K)/3  + \Pi^0_{\omega_K} \bfv,
\end{equation}
where $\Pi^0_{\omega_K}$ is simply the projection onto $\left[\mathbb{P}_0(\omega_K)\right]^3$ or $\left[\mathbb{P}_0(\omega_K)\right]$. Note that $\mathcal{P}^{\dd}_{\omega_K}$ will be only applied to each Sobolev extension $\bfu^{\pm}_{E,\dd}$ which are smooth functions involving no discontinuity. In fact, such a quasi-interpolation also results in a commutative diagram in Figure~\ref{fig:ppdiag}.

Now, with $\mathcal{C}^{\dd}_K$ and $\mathcal{P}^{\dd}_{\omega_K}$, we define the quasi-interpolation $J^{\dd}_K :\bfH(\dd;K) \rightarrow \bfS_h(K)$ where $ \bfS_h(K) =  \bfS^e_h(K)$ if $\dd = \curl$ or $ \bfS_h(K) = \bfS^f_h(K)$ if $\dd=\div$, such that
\begin{equation}
\label{quasi_interp}
J^{\dd}_K \bfu=
\begin{cases}
      & J^{\dd,+}_K  \bfu = \mathcal{P}^{\dd}_{\omega_K}  \bfu^+_{E,\dd},~~~~~~~~~~~ \text{in} ~ K^+_h, \\
      & J^{\dd,-}_K \bfu =  \mathcal{C}^{\dd}_K( \mathcal{P}^{\dd}_{\omega_K} \bfu^+_{E,\dd} ) ,~~~~~ \text{in} ~ K^-_h.
\end{cases}
\end{equation}
\begin{remark}
\label{rem_motiv}
Let us use the diagram in Figure~\ref{fig:Interpdiagram} to briefly explain the motivation. Suppose $\mathcal{P}^{\dd}_{\omega_K}(\bfu^+_{E,\dd})$ has sufficient approximation to $\bfu^+_{E,\dd}$. Then, its discrete extension $\mathcal{C}^{\dd}_K( \mathcal{P}^{\dd}_{\omega_K} \bfu^+_{E,\dd} )$ may be also expected to be a good approximation to $\bfu^-_{E,\dd}$. So the quasi-interpolation $J^{\dd}_K $ behaves close to a Hermitian-type interpolation. 
\end{remark}
Remark~\ref{rem_motiv} motivates the estimation: the ``$+$" piece follows directly from the estimate of $\mathcal{P}^{\dd}_{\omega_K} \bfu^+_{E,\dd}$ and the ``$-$" piece is given by the closeness between $\mathcal{C}^{\dd}_K(\mathcal{P}^{\dd}_{\omega_K} \bfu^+_{E,\dd})$ and $\mathcal{P}_{\omega_K}\bfu^-_{E,\dd}$. It suggests the decomposition
\begin{equation}
\label{error_Decomp} 
\bfu^- - J^{\dd,-}_K \bfu = (\bfu^- - \mathcal{P}^{\dd}_{\omega_K} \bfu^-_{E,d}) + (\mathcal{P}^{\dd}_{\omega_K} \bfu^-_{E,d} - J^{\dd,-}_K \bfu)  := \bfxi_\bfu + \bfeta_\bfu.
\end{equation}
Now, we make the idea above more rigorous. Let us first estimate $\bfxi_\bfu$ in the following lemma.
\begin{lemma}
\label{lem_est_PK}
For $\bfw\in\bfH^1(\dd;\omega_K)$, there holds
\begin{equation}
\label{lem_est_PK_eq0}
\| \bfw - \mathcal{P}^{\dd}_{\omega_K} \bfw \|_{\bfH(\dd;\omega_K)} \lesssim  h_K \| \bfw \|_{\bfH^1(\dd;\omega_K)}.
\end{equation}
\end{lemma} 
\begin{proof}
The estimate for $L^2$ norm directly follows from the approximation result of $\Pi^0_{\omega_K}$ and $\|\bfx -\bfx_K\|\lesssim h_K$ for $\bfx\in\omega_K$. The estimate for $\dd$-semi norm additionally follows from the commutative diagram in Figure~\ref{fig:ppdiag}. 
\end{proof}

As for $\bfeta_\bfu$, 
here we present a detailed discussion for $\dd=\curl$, and the argument for the $\div$ case is similar. In order to deal with jumps, we define a special norm for $\bfv_h\in \mathcal{ND}_0(\omega_K)$:
\begin{equation}
\label{special_norm}
\vertiii{\bfv_h}^2_{\curl,\omega_K} := h^2_K | \bfv_h(\bfx_K)\cdot\bar{\bfn} |^2 + \|\bfv_h \times \bar{\bfn} \|^2_{L^2(\Gamma^{\omega_K}_h)} + \| \curl \bfv_h \times \bar{\bfn} \|^2_{L^2(\Gamma^{\omega_K}_h)} + \| \curl\bfv_h\cdot\bar{\bfn} \|^2_{L^2(\Gamma^{\omega_K}_h)}.
\end{equation}
We show it is equivalent to the standard $L^2$ norm.
\begin{lemma}
\label{lem_equiv}
The norm equivalence $h^{1/2}_K\vertiii{\cdot}_{\curl,\omega_K} \simeq \|\cdot\|_{\bfH(\curl;\omega_K)}$ holds in $\mathcal{ND}_0(\omega_K)$.
\end{lemma}
\begin{proof}
Since $\bfv_h\in \mathcal{ND}_0(\omega_K)$ is a polynomial, then Lemma~\ref{lem_patch_trace_inequa_eq} implies
\begin{equation}
\begin{split}
\label{special_norm_eq1}
&  | \bfv_h(\bfx_K)\cdot\bar{\bfn} | \lesssim  h^{-3/2}_K \| \bfv_h \|_{L^2(\omega_K)}, ~~~~~ \| \bfv_h \times \bar{\bfn} \|_{L^2(\Gamma^{\omega_K}_h)} \lesssim h^{-1/2}_K \| \bfv_h \|_{L^2(\omega_K)}, \\
& \| \curl \bfv_h \times \bar{\bfn} \|_{L^2(\Gamma^{\omega_K}_h)} \lesssim h^{-1/2}_K \| \curl \bfv_h \|_{L^2(\omega_K)},  ~~~~~ \| \curl \bfv_h \cdot \bar{\bfn} \|_{L^2(\Gamma^{\omega_K}_h)} \lesssim h^{-1/2}_K \| \curl \bfv_h \|_{L^2(\omega_K)}, 
\end{split}
\end{equation}
which readily show $h^{1/2}_K \vertiii{\bfv_h}_{\curl,\omega_K}\lesssim \| \bfv_h \|_{\bfH(\curl;\omega_K)}$. As for the reverse direction, we consider the expression
\begin{equation}
\label{special_norm_eq2}
\bfv_h = \curl\bfv_h\times(\bfx-\bfx_K)/2 + \bfv_h(\bfx_K).
\end{equation}
Note that both $\bfv_h\cdot \bar{\bfn}$ and $\bfv_h\times \bar{\bfn}$ (rotation of the tangential component) can be understood as polynomials defined on $\Gamma^{\omega_K}_h$. 
Therefore, using Lemma~\ref{lem_patch_trace_inequa_eq}, we obtain
\begin{equation}
\label{special_norm_eq3}
\| \bfv_h(\bfx_K) \| \lesssim \| \bfv_h(\bfx_K) \times \bar{\bfn} \| + | \bfv_h(\bfx_K) \cdot \bar{\bfn} | \lesssim h^{-1}_K \| \bfv_h \times \bar{\bfn} \|_{L^2(\Gamma^{\omega_K}_h)} +  | \bfv_h(\bfx_K) \cdot \bar{\bfn} |.
\end{equation}
Similarly, we also have
\begin{equation}
\label{special_norm_eq4}
\| \curl\bfv_h \| \lesssim \|  \curl\bfv_h\times \bar{\bfn} \| + |\curl \bfv_h\cdot\bar{\bfn} | \simeq h^{-1}_K \|\curl\bfv_h\times \bar{\bfn} \|_{L^2(\Gamma^{\omega_K}_h)} + h^{-1}_K \|\curl\bfv_h\cdot \bar{\bfn} \|_{L^2(\Gamma^{\omega_K}_h)},
\end{equation}
which leads to $\| \curl\bfv_h \|_{L^2(\omega_K)}\lesssim h^{1/2}_K \vertiii{\bfv_h}_{\curl,\omega_K}$. Next, 
by \eqref{special_norm_eq2}, we have
\begin{equation}
\begin{split}
\label{special_norm_eq5}
\| \bfv_h \|_{L^2(\omega_K)} \lesssim h^{5/2}_K \| \curl \bfv_h \| + h^{3/2}_K \| \bfv_h(\bfx_K) \|,
\end{split}
\end{equation}
and then we apply \eqref{special_norm_eq3} and \eqref{special_norm_eq4} to each term in \eqref{special_norm_eq5}, which yields the desired reuslt.
\end{proof}

We shall see that the norm $\vertiii{\cdot}_{\curl,\omega_K}$ is suitable for estimating $\bfeta_\bfu=\mathcal{P}^{\curl}_{\omega_K}  \bfu^-_{E,\curl} - \mathcal{C}^{\curl}_K( \mathcal{P}^{\curl}_{\omega_K} \bfu^+_{E,\curl} )$.
\begin{lemma}
\label{lem_PCKP}
For $\bfu\in\bfH^1(\alpha,\beta,\curl;\Omega)$, there holds
\begin{equation}
\label{lem_PCKP_eq0}
\| \mathcal{P}^{\curl}_{\omega_K}  \bfu^-_{E,\curl} - \mathcal{C}^{\curl}_K( \mathcal{P}^{\curl}_{\omega_K} \bfu^+_{E,\curl} ) \|_{\bfH(\curl;\omega_K)} \lesssim h_K \| \bfu \|_{E,\curl,1,\omega_K}.
\end{equation}
\end{lemma}
\begin{proof}
Note that $\bfeta_\bfu \in \mathcal{ND}_0(\omega_K)$.
Then, Lemma~\ref{lem_equiv} shows
\begin{equation}
\begin{split}
\label{lem_PCKP_eq1}
\| \bfeta_\bfu \|_{\bfH(\curl;\omega_K)} \lesssim h^{1/2}_K \big( & h_K | \bfeta_\bfu(\bfx_K)\cdot\bar{\bfn} | + \| \bfeta_\bfu \times \bar{\bfn} \|_{L^2(\Gamma^{\omega_K}_h)} \\
+ &\| \curl \bfeta_\bfu \times \bar{\bfn} \|_{L^2(\Gamma^{\omega_K}_h)} + \| \curl \bfeta_\bfu \cdot\bar{\bfn} \|_{L^2(\Gamma^{\omega_K}_h)} \big).
\end{split}
\end{equation}
We need to estimate each term in the right hand side of \eqref{lem_PCKP_eq1}. As their estimates are similar, we only present the first term. By the normal jump condition for $\bfS^e_h$ in Table~\ref{tab_IFE_func} and the first inequality in \eqref{geo_estimate_1}, we have
\begin{equation*}
\begin{split}
\label{lem_PCKP_eq2}
 | \bfeta_\bfu (\bfx_K)\cdot\bar{\bfn} | & = \verti{ \mathcal{P}^{\curl}_{\omega_K}  \bfu^-_{E,\curl}(\bfx_K)\cdot\bar{\bfn} - \frac{\beta^+}{\beta^-}  \mathcal{P}^{\curl}_{\omega_K} \bfu^+_{E,\curl} (\bfx_K)\cdot\bar{\bfn} } \\
 & \lesssim h^{-1}_K \| \beta^-\mathcal{P}^{\curl}_{\omega_K}  \bfu^-_{E,\curl} \cdot\bar{\bfn} - \beta^+\mathcal{P}^{\curl}_{\omega_K}  \bfu^+_{E,\curl} \cdot\bar{\bfn} \|_{L^2(\Gamma^{\omega_K}_h)}  \\
 & \lesssim h^{-1}_K \sum_{s=\pm}  \underbrace{ \| \beta^s\mathcal{P}^{\curl}_{\omega_K}  \bfu^s_{E,\curl} \cdot\bar{\bfn} - \beta^s \bfu^s_{E,\curl} \cdot\bar{\bfn} \|_{L^2(\Gamma^{\omega_K}_h)}}_{({\rm I})}  + \underbrace{ \| \beta^+ \bfu^+_{E,\curl} \cdot\bar{\bfn} - \beta^- \bfu^-_{E,\curl} \cdot\bar{\bfn} \|_{L^2(\Gamma^{\omega_K}_h)}}_{({\rm II})}.
 \end{split}
\end{equation*}
For $({\rm I})$, using the trace inequality with the second inequality in \eqref{geo_estimate_1} and Lemma~\ref{lem_est_PK}, we obtain
\begin{equation}
\label{lem_PCKP_eq3}
({\rm I}) \lesssim h^{-1/2}_K \| \mathcal{P}^{\curl}_{\omega_K}  \bfu^s_E -  \bfu^s_E \|_{L^2(\omega_K)} + h^{1/2}_K | \mathcal{P}^{\curl}_{\omega_K}  \bfu^s_E -  \bfu^s_E |_{H^1(\omega_K)} \lesssim h^{1/2}_K \| \bfu^s_E \|_{\bfH^1(\curl;\omega_K)}
\end{equation} 
where we have used $\| \nabla \mathcal{P}^{\curl}_{\omega_K}  \bfu^s_E \|_{L^2(\omega_K)} \simeq \| \curl \mathcal{P}^{\curl}_{\omega_K}  \bfu^s_E \|_{L^2(\omega_K)} \lesssim \| \bfu^s_E \|_{\bfH^1(\curl;\omega_K)}$. The estimate of $({\rm II})$ follows from \eqref{lem_function_gamma_eq02} in Lemma~\ref{lem_diff_u} with $\dd=\curl$ and $b=\beta$. The other terms in \eqref{lem_PCKP_eq1} follow from similar arguments with the corresponding estimates in Lemma~\ref{lem_diff_u} where for the $\curl$ terms we need Remark~\ref{rem_diff_u}.
\end{proof}

The result similar to Lemma~\ref{lem_equiv} holds for $\dd=\div$. 
\begin{lemma}
\label{lem_PCKP_div}
For $\bfu\in\bfH^1(\alpha,\div;\Omega)$, there holds
\begin{equation}
\label{lem_PCKP_div_eq0}
\| \mathcal{P}^{\div}_{\omega_K}  \bfu^-_{E,\div} - \mathcal{C}^{\div}_K( \mathcal{P}^{\div}_{\omega_K} \bfu^+_{E,\div} ) \|_{\bfH(\div;\omega_K)} \lesssim h_K \| \bfu \|_{E,\div,1,\omega_K}.
\end{equation}
\end{lemma}

Now, we are ready to estimate $J^{\dd}_K$.

\begin{lemma}
\label{lem_est_JK}
For $\bfu\in\bfH^1(\alpha,\beta,\curl;\Omega)$ or $\bfu\in\bfH^1(\alpha,\div;\Omega)$, there holds
\begin{equation}
\label{lem_est_JK_eq0}
\| \bfu^{\pm}_{E,\dd} - J^{\dd,\pm}_K \bfu \|_{\bfH(\curl;\omega_K)} \lesssim  h_K \| \bfu \|_{E,\dd,1,\omega_K}.
\end{equation}
\end{lemma}
\begin{proof}
It can be seen from the Diagram~\ref{fig:Interpdiagram} and definition \eqref{quasi_interp} that the estimate of $\bfu^+_{E,\dd} - J^{\dd,+}_K \bfu$ directly follows from Lemma~\ref{lem_est_PK}. The estimate on ``$-$" piece follows from applying the triangular inequality to \eqref{error_Decomp} with Lemmas~\ref{lem_est_PK} and~\ref{lem_PCKP}.
\end{proof}

\subsection{Face and Edge Interpolation} 
\label{subsec:faceedgeinterp}
Lemma~\ref{lem_est_JK} only gives the local approximation capabilities. We then need to estimate the errors of the global spaces through the edge and face interpolations in \eqref{interpolation}. Note that the nodal interpolation of the $H^1$ IFE space has been analyzed in~\cite{2005KafafyLinLinWang}. 

Let us first present the following estimates for the bounds of the face and edge shape functions.
\begin{lemma}
\label{lem_bound_HdivIFE}
Under the conditions of Lemma~\ref{lem_Hdiv_unisol}, there holds that for $i=1,2,...,4$
\begin{equation}
\label{lem_bound_HdivIFE_eq0}
\| \bfphi^f_i \|_{L^{\infty}(K)} \lesssim  h^{-2}_K, ~~~ \text{and} ~~~  \| \div \bfphi^f_i \|_{L^{\infty}(K)} \lesssim h^{-3}_K.
\end{equation}
\end{lemma}
\begin{proof}
The second estimate in \eqref{lem_bound_HdivIFE_eq0} is trivial since $\div \bfphi^f_i $ is a constant and
$
\int_K \div \bfphi^f_i  \dd \bfx = \int_{\partial K} \bfphi^f_i\cdot\bfn \dd s.
$
For the first one, using the estimate in the proof of Lemma~\ref{lem_Hdiv_unisol} for the eigenvalues of $M_K$ for trirectangular tetrahedron and regular tetrahedron, we obtain $\verti{ \text{det}(M_K) } \ge C(\alpha)h^6_K $, where $C(\alpha)$ denotes a constant only depending on $\alpha$. Then, there holds $\| M^{-1}_K \|_{\infty} \le C(\tilde{\alpha})h^{-2}_K$, which, by \eqref{ConstructHdivIFE_eq2}, shows $\|\bfa_h\|\lesssim h^{-2}_K$. Applying it with the second one in \eqref{lem_bound_HdivIFE_eq0} to the general formula $\bfv_h=c(\bfx-\bfx_K) + \bfa_h \in \bfS^f_h(K)$ yields the desired result.
\end{proof}

\begin{lemma}
\label{lem_bound_HcurlIFE}
Under the conditions of Lemma~\ref{Hcurl_lem_unisol}, there holds that for $i=1,2,...,6$
\begin{equation}
\label{lem_bound_HcurlIFE_eq0}
\| \bfphi^e_i \|_{L^{\infty}(K)} \lesssim  h^{-1}_K, ~~~ \text{and} ~~~ \| \curl \bfphi^e_i \|_{L^{\infty}(K)} \lesssim h^{-2}_K.
\end{equation}
\end{lemma}
\begin{proof}
Lemma~\ref{lem_bound_HdivIFE} with \eqref{Hcurl_lem_unisol_eq_new4} and \eqref{Hcurl_lem_unisol_eq_new2} already yields the second inequality in \eqref{lem_bound_HcurlIFE_eq0}.
Next, note that $\phi^n$ is constructed from \eqref{Hcurl_lem_unisol_eq_new3}. By Lemma~\ref{lem_bound_HdivIFE} and $\|\bfx-\bfx_K \|\lesssim h_K$, we know $|\phi^n(\bfz)|\lesssim 1$, $\forall \bfz\in \mathcal{N}_K$.
Using the boundedness for $H^1$ IFE shape functions~\cite{2005KafafyLinLinWang}, we have $\| \nabla\phi^n \|_{\infty} \lesssim h^{-1}_K$. Then, applying it to \eqref{Hcurl_lem_unisol_eq_new4} implies the first estimate in \eqref{lem_bound_HcurlIFE_eq0}.
\end{proof}

Let us briefly explain the idea of interpolation estimationg. The key is to employ $J^{\dd}_K$ as the bridge. Let $I_h$ be $I^e_h$ if $\dd=\curl$ or $I^f_h$ if $\dd = \div$. Then, the following triangular inequality always holds
\begin{equation}
\label{u_Iu_Ju}
\| \bfu - I_h \bfu \|_{L^2(K)} \le  \| \bfu - J^{\dd}_K \bfu \|_{L^2(K)}  + \| J^{\dd}_K \bfu - I_h \bfu \|_{L^2(K)}.
\end{equation}
For the first term in the right-hand side of \eqref{u_Iu_Ju}, Lemma~\ref{lem_est_JK} almost gives the desired result, but special attention must be paid to the different partition of $\bfu$ and $J_K\bfu$ which are defined with $\Gamma$ and $\Gamma^K_h$, respectively.

\begin{lemma}
\label{lem_u_Ju_est}
For $\bfu\in\bfH^1(\alpha,\beta,\curl;\Omega)$ or $\bfu\in\bfH^1(\alpha,\div;\Omega)$, there holds
\begin{equation}
\label{lem_u_Ju_est_eq0}
\sum_{K\in\mathcal{T}^i_h}\| \bfu - J^{\dd}_K \bfu \|_{L^2(K)} \lesssim h( \| \bfu \|_{\bfH^1(\dd;\Omega^-)} +  \| \bfu \|_{\bfH^1(\dd;\Omega^+)} ).
\end{equation}
\end{lemma}
\begin{proof}
On each interface element $K$, we have
\begin{equation}
\label{thm_L2_eq1_2}
\| \bfu - J^{\dd}_K \bfu \|_{L^2(K)} \le \sum_{s=\pm} \| \bfu^s_{E,\dd} - J^{\dd,s}_K \bfu \|_{L^2(K^s_h)} + \| \bfu^+_{E,\dd} - \bfu^-_{E,\dd} \|_{L^2(K_{\textint})},
\end{equation}
where the estimate of the first term is readily given by Lemma~\ref{lem_est_JK}. 
Summing \eqref{thm_L2_eq1_2} over all the interface elements and applying \eqref{delta_strip_1}, Lemma~\ref{lem_delta} as well as Theorem~\ref{thm_ext}, we have the desired result. 
\end{proof}

Next, we turn into the term $\| J^{\dd}_K \bfu - I_h \bfu \|_{L^2(K)}$ in \eqref{u_Iu_Ju}. We need the following estimate
\begin{lemma}
\label{lem_Ju_est}
For $\bfu\in\bfH^1(\alpha,\beta,\curl;\Omega)$ or $\bfu\in\bfH^1(\alpha,\div;\Omega)$, given a $\delta$-strip with $\delta \lesssim h^2$, there holds
\begin{equation}
\label{lem_Ju_est_eq0}
\sum_{K\in\mathcal{T}^i_h} \| J^{\dd,\pm}_K \bfu \|_{L^2(S_{\delta}\cap \omega_K)} \lesssim h \| \bfu \|_{\bfH^1(\dd;\Omega^+)} .
\end{equation}
\end{lemma}
\begin{proof}
For the ``$+$" component, with the estimates
$
|S_{\delta}\cap\omega_K | \lesssim h^4_K ~ \text{and} ~ \| \bfx - \bfx_K \|_{L^2(S_{\delta}\cap\omega_K)}\lesssim h^{3}_K,
$
we have
\begin{equation}
\label{lem_Sdelta1_eq1}
\begin{split}
\| J^{\dd,+}_K \bfu \|_{L^2(S_{\delta}\cap\omega_K)} & \lesssim h^3_K \| \Pi^0_{\omega_K} \dd \bfu^+_{E,\dd} \| + h^2_K \| \Pi^0_{\omega_K}  \bfu^+_{E,\dd} \| \\
& \lesssim h^{3/2}_K \| \Pi^0_{\omega_K} \dd \bfu^+_{E,\dd} \|_{L^2(\omega_K)} + h^{1/2}_K \| \Pi^0_{\omega_K}  \bfu^+_{E,\dd} \|_{L^2(\omega_K)} \\
&  \lesssim h^{3/2}_K \|  \dd \bfu^+_{E,\dd} \|_{L^2(\omega_K)} + h^{1/2}_K \|  \bfu^+_{E,\dd} \|_{L^2(\omega_K)}
\end{split}
\end{equation}
where we have used the boundedness property for the projection $\Pi^0_{\omega_K}$ in the last inequality. Then, the desired estimate follows from summing \eqref{lem_Sdelta1_eq1} over all the interface elements, using the finite overlapping of $\omega_K$, and applying the $\delta$-strip argument from Lemma~\ref{lem_delta} with $\delta\simeq h_K$ as well as Theorem~\ref{thm_ext}. For the ``$-$" side, by \eqref{lem_ABeigen_eq1}, we have the boundedness $\| A_K \|\lesssim 1$ and $\| B_K \|\lesssim 1$. Therefore, there also holds that
\begin{equation}
\label{lem_Sdelta1_eq2}
\| J^{\dd,-}_K \bfu \|_{L^2(S_{\delta}\cap\omega_K)}  \lesssim h^3_K \| \Pi^0_{\omega_K} \dd \bfu^+_{E,\dd} \| + h^2_K \| \Pi^0_{\omega_K}  \bfu^+_{E,\dd} \|
\end{equation}
of which the estimete is the same as \eqref{lem_Sdelta1_eq1}.
\end{proof}

Now, with the estimate in Lemma~\ref{lem_Ju_est}, we proceed to show the second term in \eqref{u_Iu_Ju} and split the discussion for the edge and face elements.

\subsubsection{Face Interpolation}
\begin{figure}[h]
\centering
\begin{subfigure}{.35\textwidth}
    \includegraphics[width=1.4in]{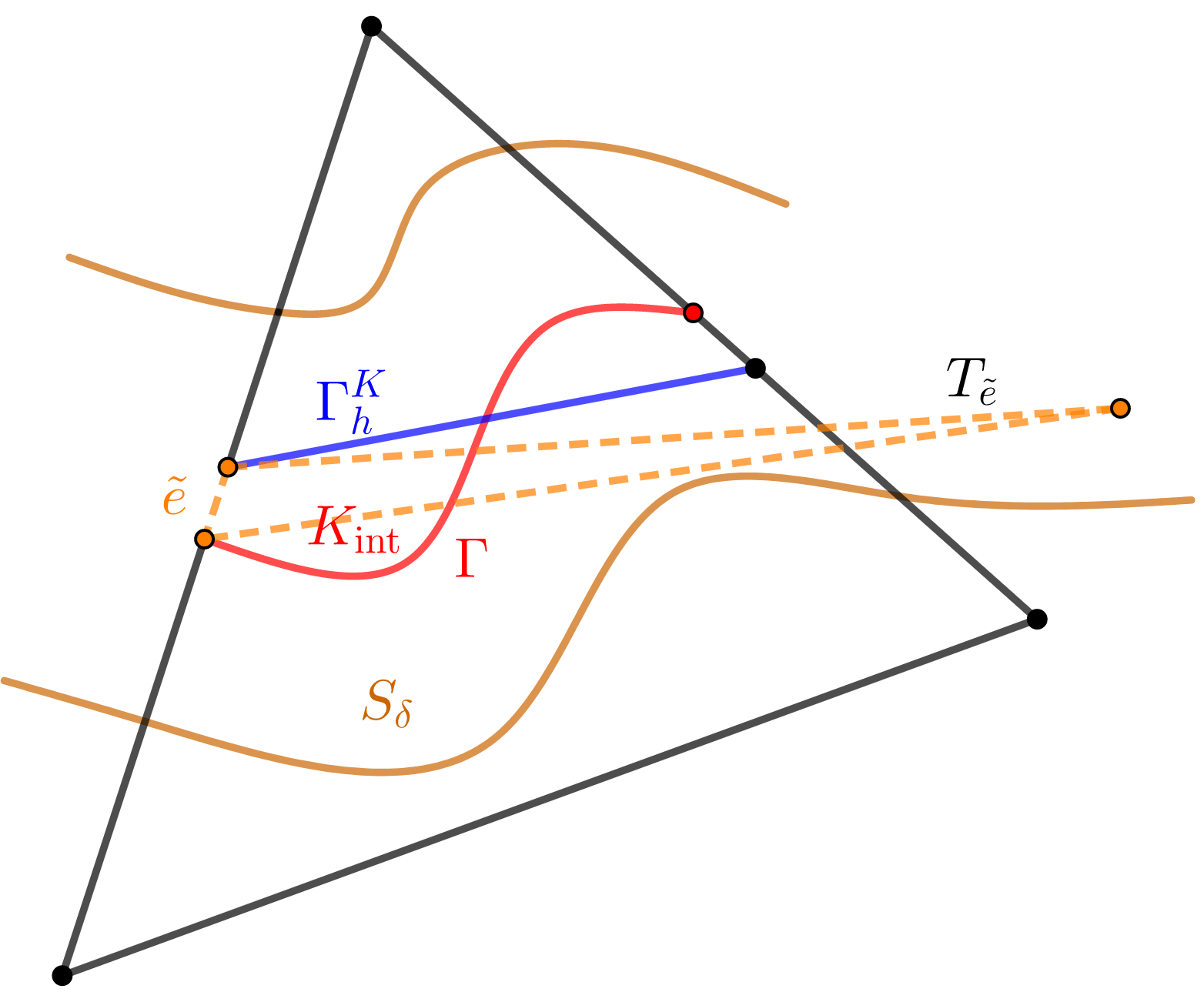}
     \label{fig:assump21} 
\end{subfigure}
~~~~~
\begin{subfigure}{.4\textwidth}
     \includegraphics[width=2in]{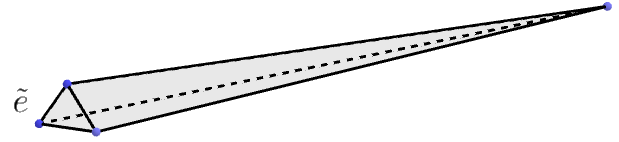}
     \label{fig:assump22} 
\end{subfigure}
     \caption{Illustration of the Assumptions \hyperref[asp:A2]{A2} and \hyperref[asp:A3]{A3}. Left: a 2D illustration. Right: a tetrahedron $T_{\tilde{e}}$ has a long height supporting a small face, which has the volume $\mathcal{O}(|\tilde{e}^2h_K|)$.}
  \label{fig:assump2} 
\end{figure}

Let us first consider the face interpolation, and in the following discussion we fix $\dd = \div$. We need the following assumption to handle the mismatching portion on faces:
\begin{itemize}
  \item[(\textbf{A2})] \label{asp:A2} For each interface element $K\in\mathcal{T}^i_h$ and any of its face $F$, define $\tilde{F}=F \cap K_{\textint}$. Assume that there exists a $\delta$-strip $S_{\delta_1}$ with $\delta_1\simeq O(h^2)$ such that for every $K\in\mathcal{T}^i_h$, there is a pyramid $P_{\tilde{F}}\subset \omega_K \cap S_{\delta_1}$ with $\tilde{F}$ as its base satisfying that the height $l_{\tilde{F}}$ of $P_{\tilde{F}}$ supporting $\tilde{F}$ is $\mathcal{O}(h_K)$.
\end{itemize}
We refer readers to the 2D illustration of this assumption in the left plot of Figure~\ref{fig:assump2}. 

\begin{theorem}
\label{thm_Hdiv_L2}
For $\bfu\in \bfH^1(\alpha,\div;\Omega)$, under the Assumption \hyperref[asp:A2]{A2}, there holds
\begin{equation}
\label{thm_Hdiv_L2_eq0}
\| \bfu - I^f_h \bfu \|_{L^2(\Omega)} \lesssim h( \| \bfu \|_{\bfH^1(\div;\Omega^-)} +  \| \bfu \|_{\bfH^1(\div;\Omega^+)} ).
\end{equation}
\end{theorem}
\begin{proof}
We only consider interface elements. Given $K\in\mathcal{T}^i_h$, by \eqref{u_Iu_Ju} and Lemma~\ref{lem_u_Ju_est}, it remains to estimate $\| J^{\dd}_K \bfu - I^f_h \bfu \|_{L^2(K)}$. We need an auxiliary function $\bfw $ such that $\bfw^{\pm} = J^{\dd,\pm}_K \bfu - \bfu^{\pm}_{E,\dd}$ in $K^{\pm}$, i.e., it is partitioned with $\Gamma$. Note that $\bfw$ is slightly different from $J^{\dd}_K\bfu - \bfu$ as $J^{\dd}_K\bfu$ is defined with $\Gamma^K_h$, and thus this difference is only on $K_{\textint}$. To control the mismatch, by Assumption \hyperref[asp:A3]{A3} and the trace inequality for polynomials~\cite{2003WarburtonHesthaven}, we have
\begin{equation}
\label{thm_Hdiv_L2_eq2}
\| J^{\dd,+}_K \bfu - J^{\dd,-}_K \bfu \|_{L^2(\tilde{F})} \lesssim h^{-1/2}_K ( \| J^{\dd,+}_K \bfu \|_{L^2(S_{\delta_1\cap\omega_K})} +  \| J^{\dd,-}_K \bfu \|_{L^2(S_{\delta_1\cap\omega_K})} ).
\end{equation}
Note that each component of $\bfw^{\pm}$ can be naturally extended to the whole element and whole patch. Noting that $\bfw^{\pm} \in \bfH^{1/2}(F)$ for each face of $F$, then by the trace inequality and Lemma~\ref{lem_est_JK} we have
\begin{equation}
\label{thm_Hdiv_L2_eq3}
\| \bfw^{\pm} \|_{L^2(F^{\pm})} \lesssim \| \bfw^{\pm} \|_{L^2(F)} \lesssim h^{-1/2}_K \| \bfw \|_{L^2(K)} + h^{1/2}_K | \bfw |_{H^1(K)} \lesssim  h^{1/2}_K \| \bfu \|_{E,\div,1,\omega_K}.
\end{equation}
Now, by the definition of the face interpolation, we obtain
\begin{equation}
\begin{split}
\label{thm_Hdiv_L2_eq4}
& \| J^{\dd}_K \bfu - I^f_h \bfu \|_{L^2(K)} = \| I^f_h (J^{\dd}_K \bfu - \bfu) \|_{L^2(K)} \\
 \lesssim & \left( \sum_{i=1}^4 |F_i|^{1/2} (\| \bfw^+\cdot\bfn \|_{L^2(F_i)} +  \| \bfw^-\cdot\bfn \|_{L^2(F_i)})  + |\tilde{F}_i|^{1/2} \| J^{\dd,+}_K \bfu - J^{\dd,-}_K \bfu \|_{L^2(\tilde{F}_i)} \right) \| \bfphi^f_i \|_{L^2(K)}  \\
 \lesssim &  h_K \| \bfu \|_{E,\div,1,\omega_K} +  (\| J^{\dd,+}_K \bfu \|_{L^2(S_{\delta_1\cap\omega_K})} +  \| J^{\dd,-}_K \bfu \|_{L^2(S_{\delta_1\cap\omega_K})} ),
\end{split}
\end{equation}
where we have used \eqref{thm_Hdiv_L2_eq2}, \eqref{thm_Hdiv_L2_eq3} and Lemma~\ref{lem_bound_HdivIFE}. Summing \eqref{thm_Hdiv_L2_eq4} over all the elements, using the finite overlapping property of $\omega_K$, and applying Lemma~\ref{lem_Ju_est} and Theorem~\ref{thm_ext} yields the desired result.
\end{proof}

In addition, the estimate of $\div(\bfu - I^f_h\bfu)$ is trivial due to the commutative property by Lemma~\ref{lem_IFE_commut}.
\begin{theorem}
\label{thm_Hdiv_div}
For $\bfu\in \bfH^1(\alpha,\div;\Omega)$, there holds
\begin{equation}
\label{thm_Hdiv_div_eq0}
\| \div(\bfu - I^f_h \bfu) \|_{L^2(\Omega)} \lesssim h( \| \bfu \|_{\bfH^1(\div;\Omega^-)} +  \| \bfu \|_{\bfH^1(\div;\Omega^+)} ).
\end{equation}
\end{theorem}


\subsubsection{Edge Interpolation}

Next, we proceed to estimate the approximation capabilities of the global space $\bfS^e_h$ through the interpolation $I^e_h$. We need a similar assumption to \hyperref[asp:A2]{A2}.
\begin{itemize}
  \item[(\textbf{A3})] \label{asp:A3} For each interface element $K\in\mathcal{T}^i_h$ and any of its edge $e$, define $\tilde{e}=e \cap K_{\textint}$, i.e., the portion of the edge contained in the mismatching region $K_{\textint}$. Assume that there exists a $\delta$-strip $S_{\delta_2}$ with $\delta_2\simeq O(h^2)$ such that for every $K\in\mathcal{T}^i_h$, there is a tetrahedron $T_{\tilde{e}} \subset \omega_K \cap S_{\delta_2}$ with $\tilde{e}$ as one of its edges satisfying $|T_{\tilde{e}}| \gtrsim |\tilde{e}|^2h_K$. 
\end{itemize}
We illustrate this assumption in Figure~\ref{fig:assump2}. In the following disucssion, we fix $\dd=\curl$.
\begin{theorem}
\label{thm_L2}
For $\bfu\in \bfH^1(\alpha,\beta,\curl;\Omega)$, under Assumption \hyperref[asp:A3]{A3}, there holds
\begin{equation}
\label{thm_L2_eq0}
\| \bfu - I^e_h \bfu \|_{L^2(\Omega)} \lesssim h( \| \bfu \|_{\bfH^1(\curl;\Omega^-)} +  \| \bfu \|_{\bfH^1(\curl;\Omega^+)} ).
\end{equation}
\end{theorem}
\begin{proof}
We still only consider interface elements.
Similar to the argument of Theorem~\ref{thm_Hdiv_L2}, we only need to estimate $\| J^{\dd}_K \bfu - I^e_h \bfu \|_{L^2(K)}$. Define the same auxiliary function $\bfw $ partitioned with $\Gamma$ and $\bfw^{\pm} = J^{\dd,\pm}_K \bfu - \bfu^{\pm}_{E,\dd}$ in $K^{\pm}$. Let us estimate $J^{\dd,+}_K \bfu - J^{\dd,-}_K \bfu$ on the small subedge $\bar{e}$. By Assumption \hyperref[asp:A3]{A3} and the trace inequality for polynomials~\cite{2003WarburtonHesthaven}, we obtain
\begin{equation}
\label{thm_L2_eq1_1}
\| J^{\dd,+}_K \bfu - J^{\dd,-}_K \bfu \|_{L^2(\bar{e})} \lesssim h^{-1/2}_K|\tilde{e}|^{-1/2} ( \| J^{\dd,+}_K \bfu \|_{L^2(S_{\delta_2\cap\omega_K})} +  \| J^{\dd,-}_K \bfu \|_{L^2(S_{\delta_2\cap\omega_K})} ).
\end{equation}
In addition, as $\bfu^{\pm}_{E,\dd}\in \bfH^1(\curl;\Omega)$, the first moment $\int_e \bfu^{\pm}_{E,\dd}\cdot\bft \dd s$ is well defined on each edge $e$. Furthermore, we have $\bfu^{\pm}_{E,\dd}\in \bfH^{1/2}(F)\cap \bfH(\rot;F)$ on each face $F$. Thus, by Lemma 4.4 in~\cite{2020BeiroMascotto} and the trace inequality, we obtain for each face $F$ of $K$ and any one of its edge $e$ that
\begin{equation}
\begin{split}
\label{thm_L2_eq3}
\| \bfw^{\pm} \|_{L^1(e^{\pm})} & \lesssim \| \bfw^{\pm} \|_{L^1(e)} \lesssim \| \bfw^{\pm} \|_{L^2(F)} + h^{1/2}_F \| \bfw^{\pm} \|_{H^{1/2}(F)} + h_F \| \rot_F\bfw^{\pm} \|_{L^2(F)} \\
& \lesssim h^{-1/2}_K \| \bfw^{\pm} \|_{L^2(K)} + h^{1/2}_K \| \bfw^{\pm} \|_{H^1(K)} \\
&+ h^{1/2}_K \| \curl \bfw^{\pm} \|_{L^2(K)} + h^{3/2}_K \| \curl \bfw^{\pm} \|_{H^1(K)}  \lesssim h^{1/2}_K \| \bfu \|_{E,\curl,1,\omega_K},
\end{split}
\end{equation}
where we have used Lemma~\ref{lem_est_JK} in the forth inequality with the following estimates
\begin{equation*}
\begin{split}
& \| \bfw^{\pm} \|_{H^1(K)} \lesssim \| \bfu^{\pm}_{E,\dd} \|_{H^1(K)} + \| \curl J^{\dd,\pm}_K\bfu \|_{L^2(K)} \lesssim \| \bfu \|_{E,\curl,1,\omega_K}, \\
&  \| \curl \bfw^{\pm} \|_{H^1(K)} \lesssim \| \curl \bfu^{\pm}_{E,\dd}  \|_{H^1(K)}  \lesssim \| \bfu \|_{E,\curl,1,\omega_K}.
\end{split}
\end{equation*}
Then, it follows from the definition of $J^{\dd}_K$ together with \eqref{thm_L2_eq1_1}, \eqref{thm_L2_eq3} and Lemma~\ref{lem_bound_HcurlIFE} that
\begin{equation}
\begin{split}
\label{thm_L2_e4}
 & \| J^{\dd}_K\bfu - I^e_h  \bfu \|_{L^2(K)} = \|  I^e_h ( J^{\dd}_K\bfu - \bfu ) \|_{L^2(K)}  \\
  \le & \left( \sum_{i=1}^6 \int_{e_i} \verti{ \bfw\cdot\bft }\dd s + |\tilde{e}_i|^{1/2} \| J^{\dd,+}_K \bfu - J^{\dd,-}_K \bfu \|_{L^2(\tilde{e}_i)}   \right) \| \bfphi^e_i \|_{L^2(K)} \\
 \lesssim & h^{1/2}_K \sum_{i=1}^6  ( \| \bfw^+ \|_{L^1(e_i)} + \| \bfw^- \|_{L^1(e_i)} ) + ( \| J^{\dd,+}_K\bfu \|_{L^2(S_{\delta_2\cap\omega_K})} +  \| J^{\dd,-}_K\bfu \|_{L^2(S_{\delta_2\cap\omega_K})} ) \\
  \lesssim & h_K \| \bfu \|_{E,\curl,1,\omega_K} + ( \| J^{\dd,+}_K\bfu \|_{L^2(S_{\delta_2\cap\omega_K})} +  \| J^{\dd,-}_K \bfu\|_{L^2(S_{\delta_2\cap\omega_K})} ).
 \end{split}
\end{equation}
Summing \eqref{thm_L2_e4} over all the interface elements, using the finite overlapping property of $\omega_K$, and applying Lemma~\ref{lem_Ju_est} and Theorem~\ref{thm_ext} leads to the desired estimate.
\end{proof}

\begin{remark}
One key in the proof is to use Assumption \hyperref[asp:A3]{A3} with the trace inequality for polynomials to control $J^{\dd,\pm}_K\bfu$ on edges. If this is done for $\bfu^{\pm}_{E,\dd}$, one must be careful about bounding $\|\bfu_{E,\dd}\|_{L^2(\tilde{e})}$ since it may not be well-defined for $\bfu\in\bfH^1(\omega_K)$. 
\end{remark}

Finally, Theorem~\ref{thm_Hdiv_L2} immediately yields the estimate for $\curl (\bfu -  I^e_h\bfu)$, $\bfu\in \bfH^1(\alpha,\beta,\curl;\Omega)$.
\begin{theorem}
\label{thm_Hcurl_L2}
Under Assumption \hyperref[asp:A3]{A3}, for $\bfu\in \bfH^1(\alpha,\beta,\curl;\Omega)$, there holds
\begin{equation}
\label{thm_Hcurl_L2_eq0}
\| \curl( \bfu - I^e_h \bfu ) \|_{L^2(\Omega)} \lesssim h( \| \bfu \|_{\bfH^1(\curl;\Omega^-)} +  \| \bfu \|_{\bfH^1(\curl;\Omega^+)} ).
\end{equation}
\end{theorem}
\begin{proof}
Note that $\div\curl\bfu=0$, and thus $\curl \bfu\in \bfH^1(\alpha,\div;\Omega)$. Then the estimate follows from the commutative property by Lemma~\ref{lem_IFE_commut} and Theorem~\ref{thm_Hdiv_L2}.
\end{proof}


\section{Numerical Experiments}
\label{sec:numer}

In this section, we present a group of numerical experiments to demonstrate the performance of the proposed method. 
We fix the cubic domain to be $\Omega = (-1,1)^3$ but with varying interface shapes. The background Cartesian mesh is generated by $N^3$ uniform cubes of which each is partitioned into tetrahedra. In~\cite{2020GuoLinZou} for the 2D case, we have compared the PG formulation with other penalty-type schemes and observed its significant advantages in producing optimal convergent solutions. 
A similar behavior can be also observed in the 3D case, and here we omit their comparison.

\begin{figure}[h]
\centering
\begin{subfigure}{.31\textwidth}
    \includegraphics[width=1.3in]{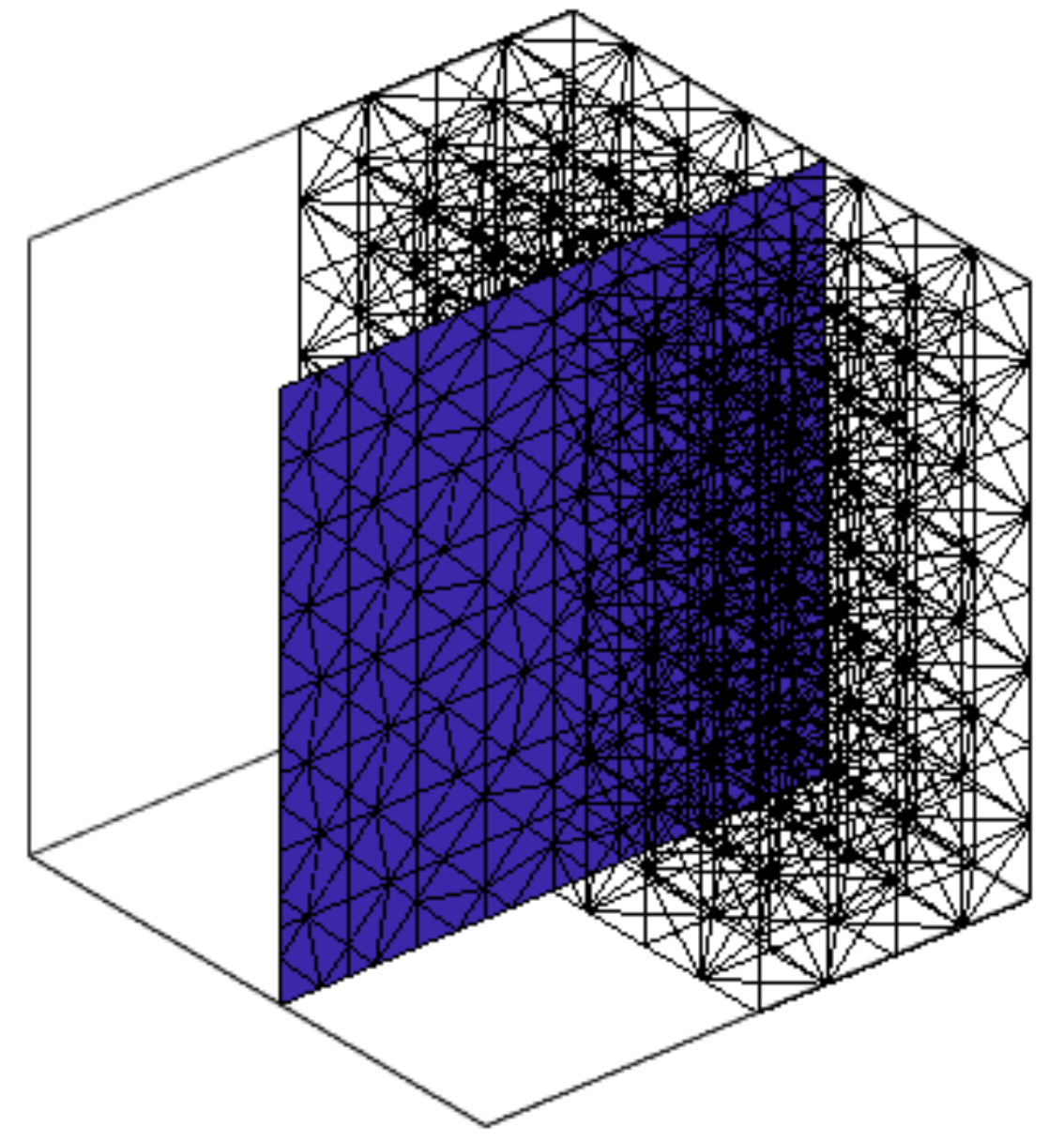}
     \label{fig:interfexample1} 
\end{subfigure}
~
\begin{subfigure}{.31\textwidth}
     \includegraphics[width=1.3in]{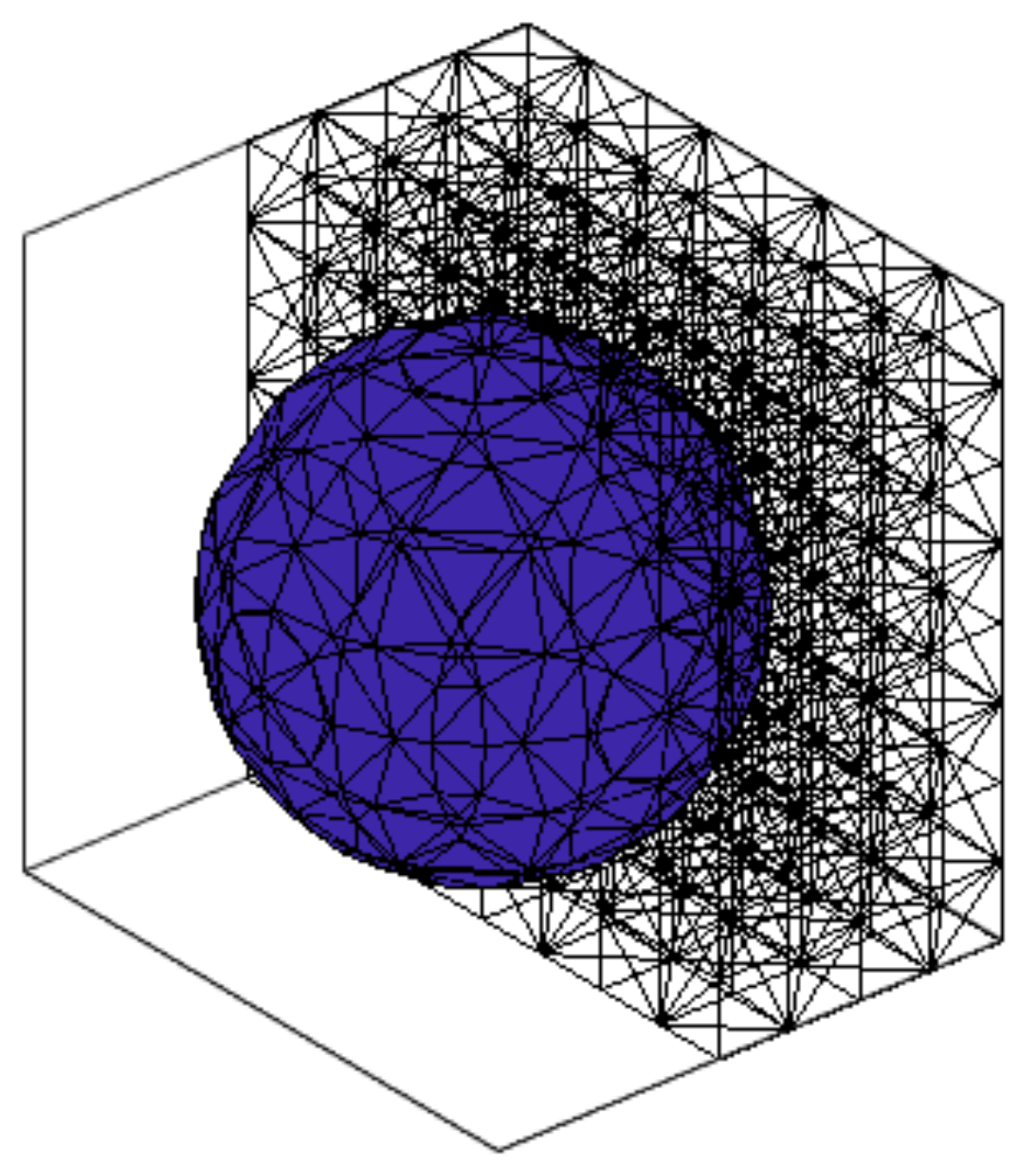}
     \label{fig:interfexample2} 
\end{subfigure}
~
 \begin{subfigure}{.31\textwidth}
    \includegraphics[width=1.3in]{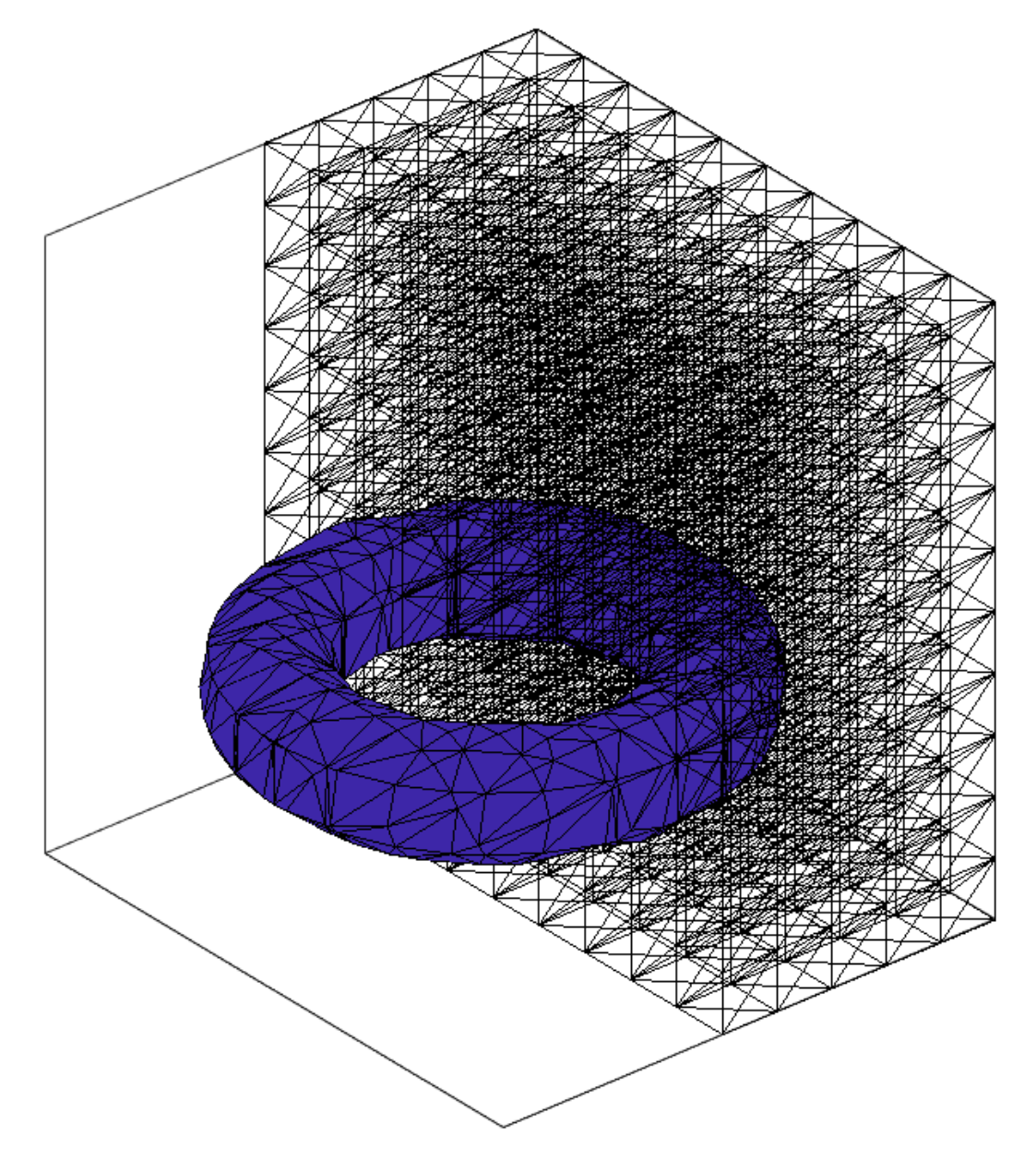}
     \label{fig:interfexample3} 
\end{subfigure}
     \caption{Interface shapes: flat interface (left), a spheric interface (middle) and a torus interface (right).}
  \label{fig:interfexample} 
\end{figure}

\subsection{Test 1(Stability)}
\begin{figure}[h]
\centering
\begin{subfigure}{.22\textwidth}
    \includegraphics[width=1.4in]{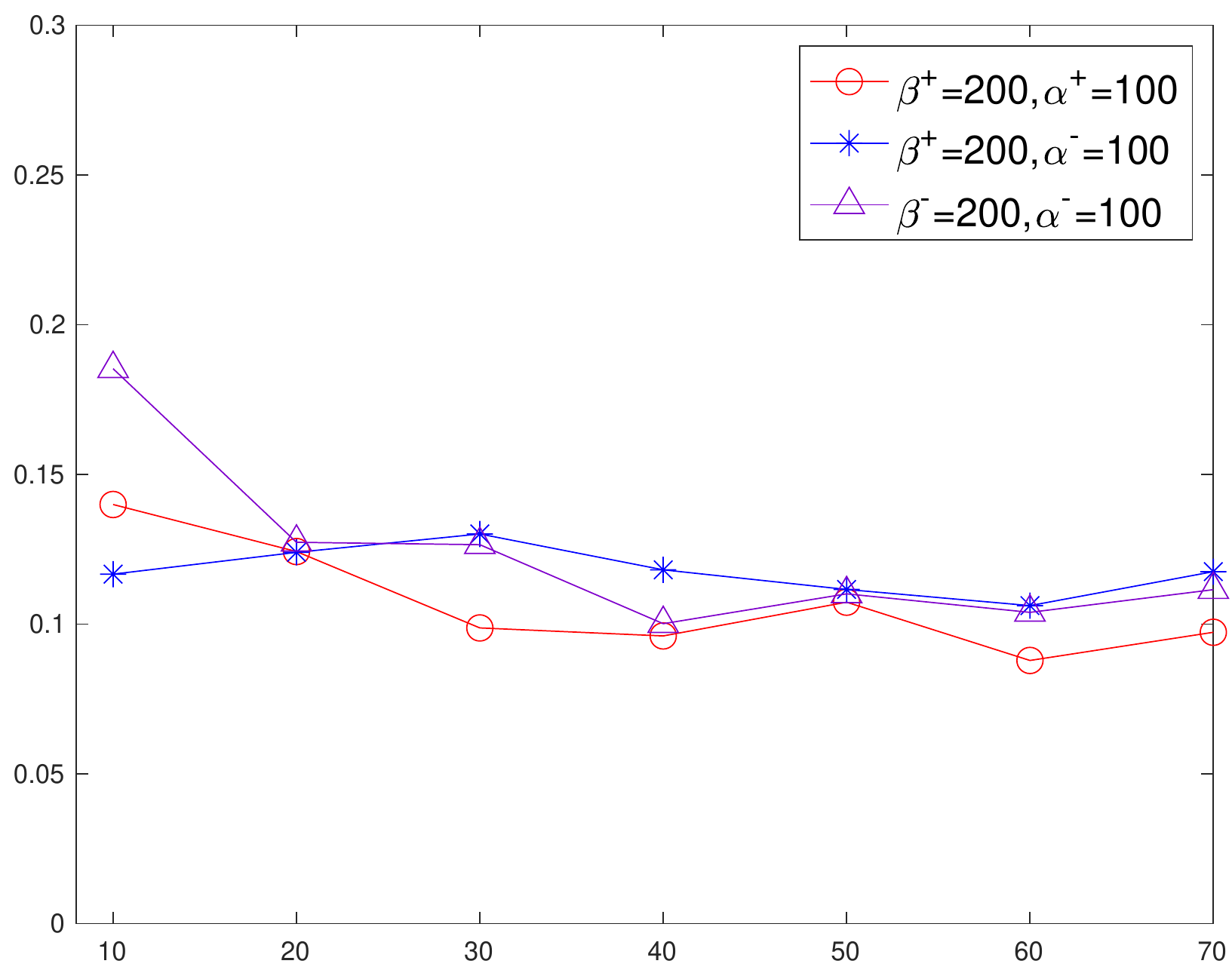}
     \label{fig:infsupconstant} 
\end{subfigure}
~
\begin{subfigure}{.22\textwidth}
    \includegraphics[width=1.4in]{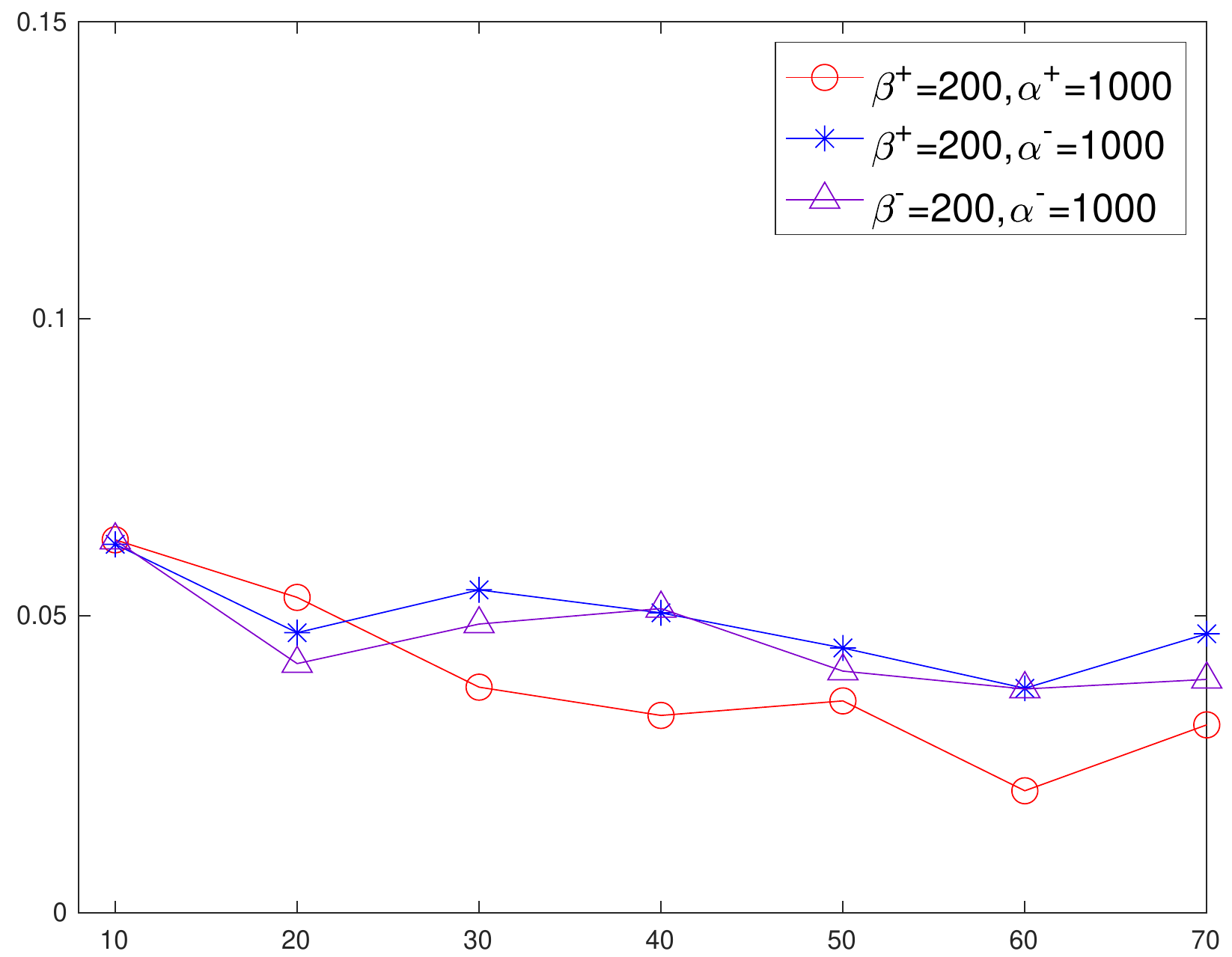}
     \label{fig:infsupconstant2} 
\end{subfigure}
~
\begin{subfigure}{.22\textwidth}
    \includegraphics[width=1.4in]{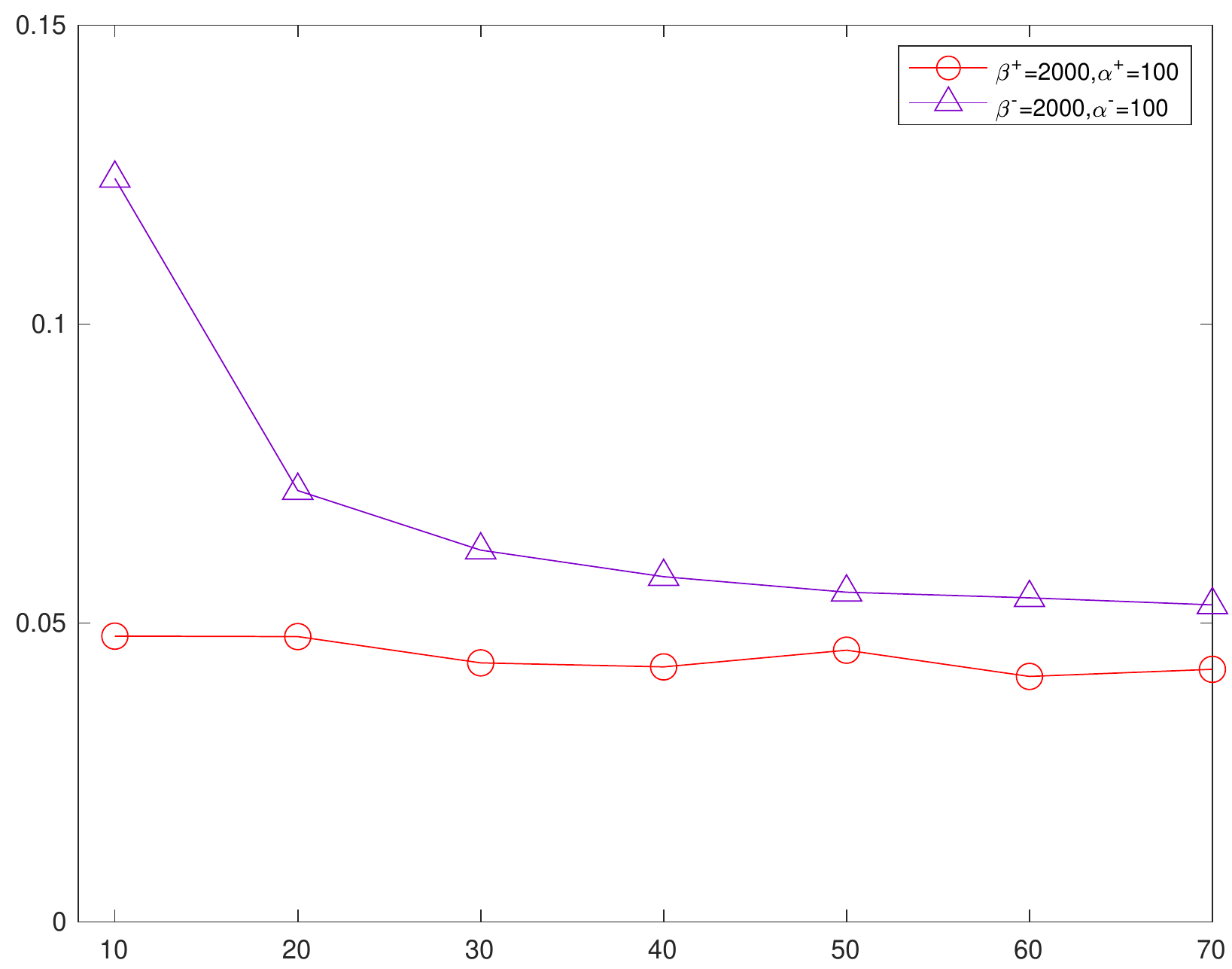}
     \label{fig:infsupconstant31} 
\end{subfigure}
~
\begin{subfigure}{.22\textwidth}
    \includegraphics[width=1.4in]{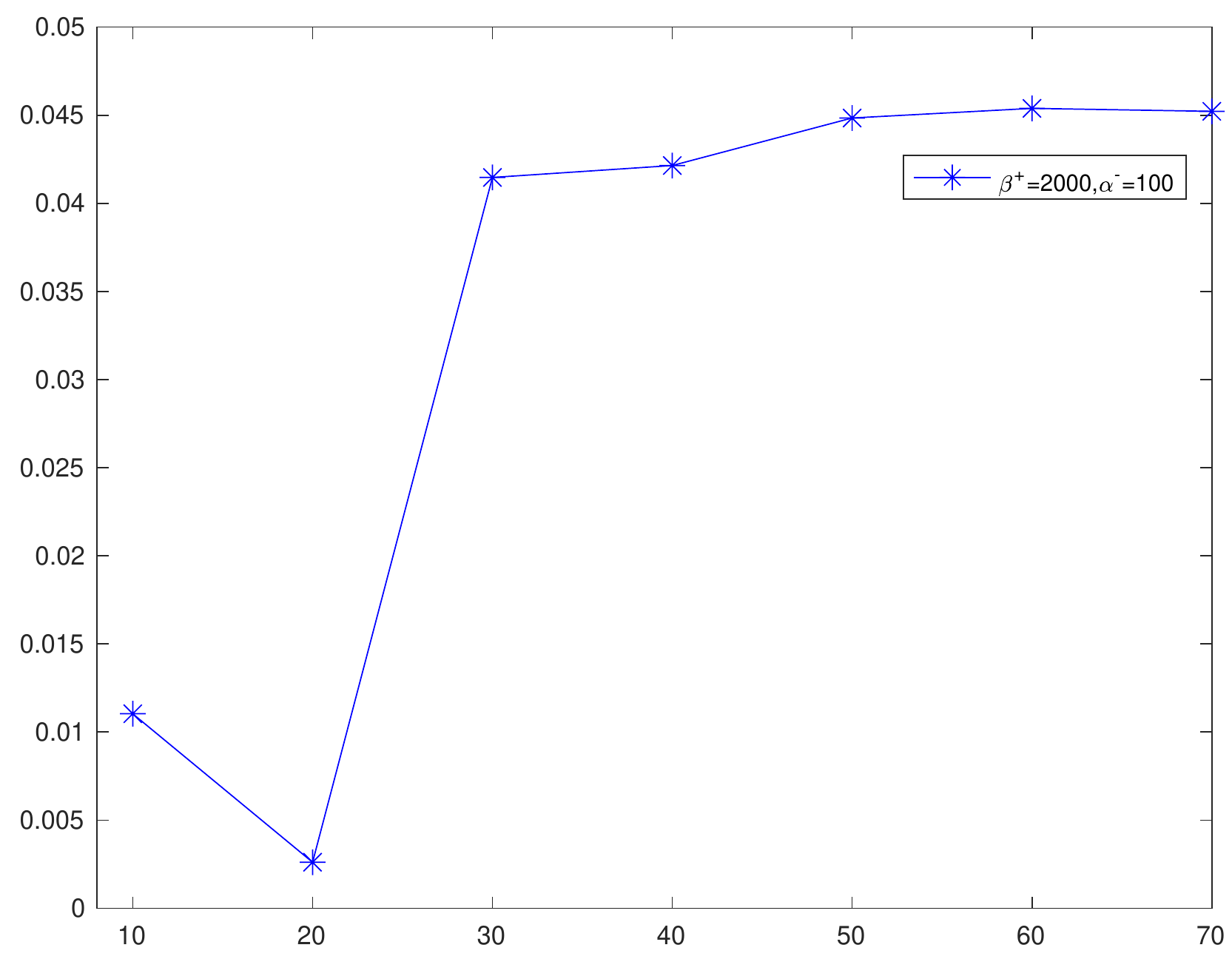}
     \label{fig:infsupconstant32} 
\end{subfigure}
     \caption{Estimates of $\eta_s$ for various $\alpha$, $\beta$ and mesh size.}
  \label{fig:Example1infsup} 
\end{figure}

In the first test, we study the numerical stability of the method with respect to mesh size and interface location. We first provide a group of numerical estimates for the quantity $\eta_s$ in \eqref{infsup_cond} indicating a lower bound of the \textit{inf-sup} constant. 
We consider the spheric interface shape shown in the middle plot of Figure~\ref{fig:interfexample} for the three groups of parameters. 
\begin{equation}
\begin{split}
\label{alphabeta}
\left( \frac{\max(\beta^+,\beta^-)}{\min(\beta^+,\beta^-)} , ~~ \frac{\max(\alpha^+,\alpha^-)}{\min(\alpha^+,\alpha^-)}   \right) = (200,100), (200,1000), (2000,100)
\end{split}
\end{equation}
where the minimal values are always set as $1$, see the detailed parameter illustration in the right corner of plots in Figure~\ref{fig:Example1infsup}. Then, we consider a sequence of meshes with $N=10,20,30,40,50,60,70$. We plot the values of $\eta_s$ in Figure~\ref{fig:Example1infsup}. Here, the computational results for the first two groups of parameters are presented in the first two plots, but for the third group, we separate the parameter $(\beta^+,\alpha^-)=(2000,100)$ in a single plot to show more details. Overall, we can see that $\eta_s$ is bounded below. More specifically, we can make the following numerical observations:
\begin{itemize}
\item $\eta_s$ slightly decreases at the beginning, and then the curve proceeds to be flat as the mesh becomes finer. 
\item $\eta_s$ highly depends on the contrast of the parameters, and the larger contrast leads to smaller $\eta_s$. Particularly, the values decrease from the first to the third group of parameters.
\item Furthermore, for the third group that the parameters have the ratio $(2000,100)$, when $(\alpha^-,\beta^+)=(100,2000)$, $\eta_s\approx 2.6\times 10^{-3}$ for the mesh $N=20$, which is indeed close to zero. But it turns back to around $0.05$ for finer meshes.  
\end{itemize}
These observations are consistent with the analysis (cf.\,\cite{2020GuoLinZou}) that the \textit{inf-sup} condition only holds for sufficiently small mesh size.

Next, we test the condition number of $\bfB^{PG}$ with respect to various interface locations. Consider a flat interface $(\bfx-\bfx_0)\cdot\bfn = 0$ passing through $\bfx_0=[\epsilon,0,0]$ with $\epsilon=10^{-j}$, $j=1,2,...,6$ and the normal vector $\bfn=[1,0,0]^{\top}$ such that the interface gradually moves to one face/edge of interface elements. Thus the small subelements exist for every interface element. The condition numbers are shown in Table~\ref{table:example1} for various $\epsilon$. In fact, as the interface element shrinks to a non-interface element, the IFE functions will converge to standard FE functions which is ``good" for computation. 
\begin{table}[H]
\begin{center}
\begin{tabular}{|c |c | c|c | c|c | c|}
\hline
$\epsilon$   & $10^{-1}$      & $10^{-2}$             & $10^{-3}$          & $10^{-4}$            & $10^{-5}$     & $10^{-6}$   \\ \hline
\text{cond}($\bfB^{PG}$)                         &    88267          &      83052       &    78834         &     84570          &    78253   &       81324           \\ \hline
\end{tabular}
\end{center}
\caption{Condition numbers for the flat interface with various $\epsilon$.}
\label{table:example1}
\end{table}


\subsection{Example 2 (Optimal Convergence)}
\begin{figure}[h]
\centering
\begin{subfigure}{.31\textwidth}
    \includegraphics[width=1.6in]{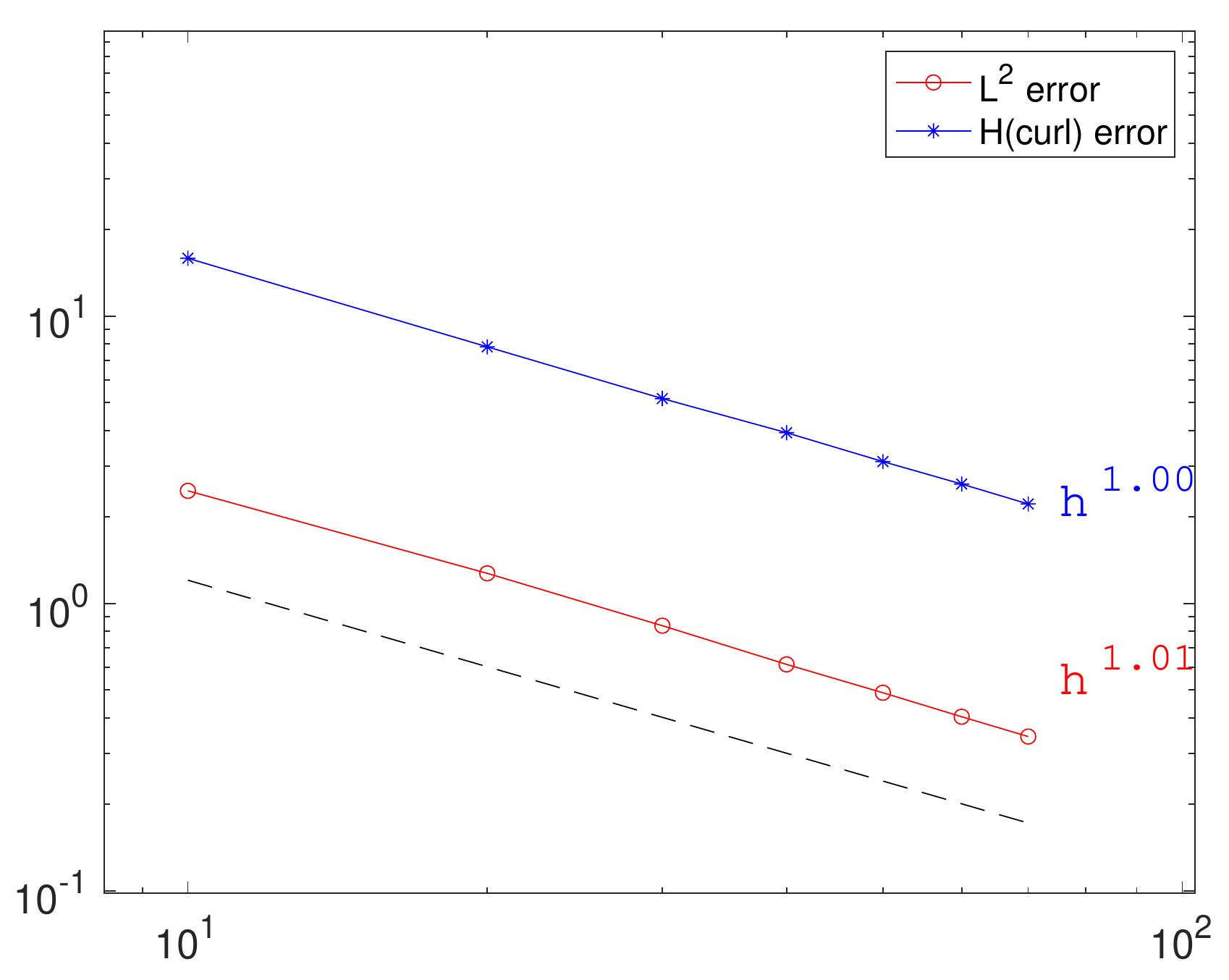}
     \caption{$\alpha^-,\alpha^+)=(\beta^-,\beta^+)=(1,100)$}
     \label{fig:Example2m1p100} 
\end{subfigure}
~
\begin{subfigure}{.31\textwidth}
     \includegraphics[width=1.6in]{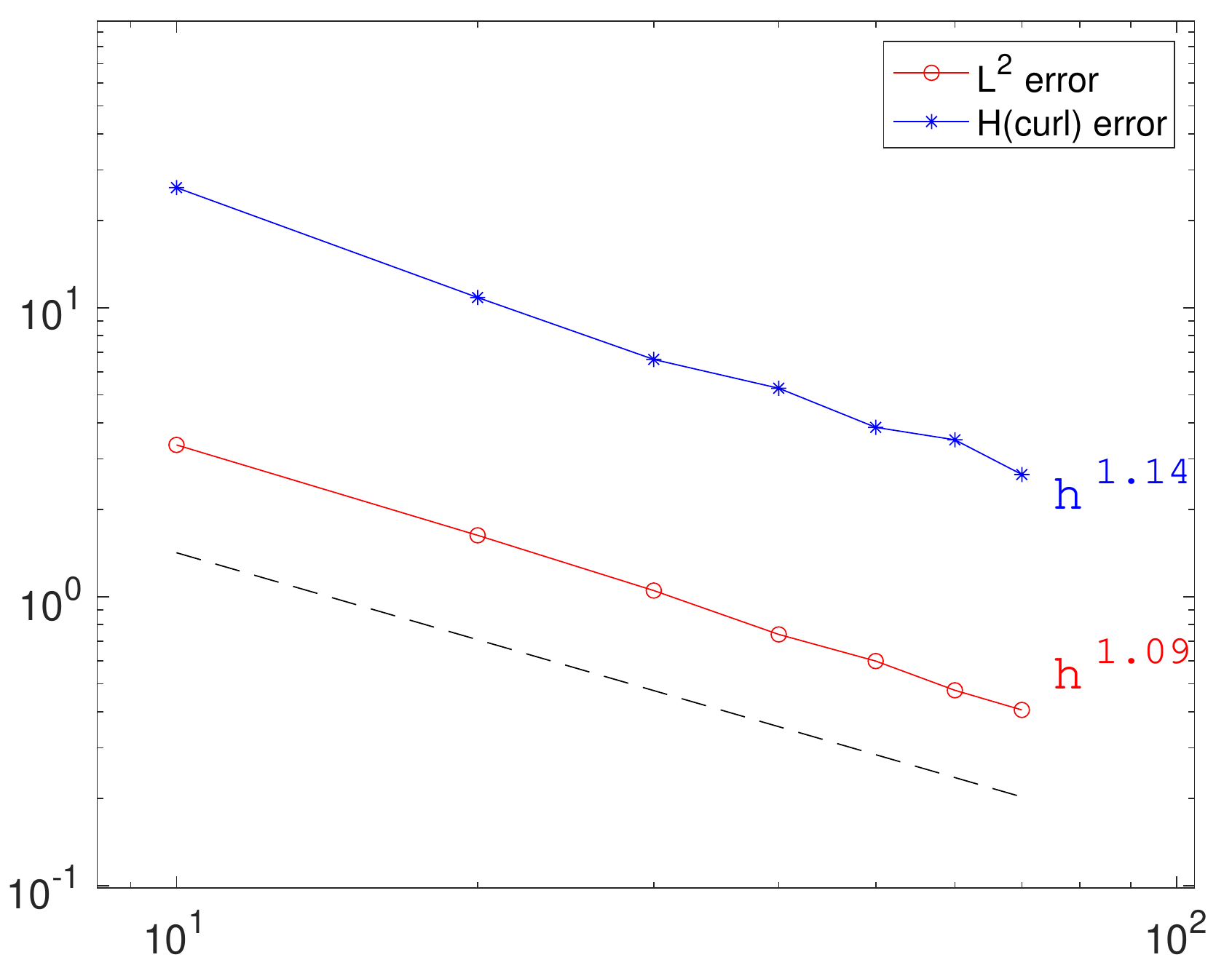}
     \caption{$(\alpha^-,\alpha^+)=(\beta^-,\beta^+)=(1,1000)$}
     \label{fig:Example2m1p1000} 
\end{subfigure}
     \caption{Errors for Example 2 where the dashed line is the reference line indicating $\mathcal{O}(h)$.}
  \label{fig:Example2Error} 
\end{figure}

In the second test, we investigate the optimal convergence rates with respect to mesh size and speed of the fast solver. Consider a spheric interface with a radius $r_1$ shown in the left plot of Figure~\ref{fig:interfexample}. We still employ \eqref{Max_model3} as the model equation, and use the benchmark example from~\cite{2020GuoLinZou} of which the analytical solution is given by
\begin{equation}
\label{exactu_2}
\bfu = 
\begin{cases}
      & \frac{1}{\beta^- } \bfx +  \frac{1}{\alpha^- } n_1R_1(\bfx)[(x_2-x_3),(x_3-x_1),(x_1-x_2)]^{\top}  ~~~~~~~~~~~~ \text{in}~ \Omega^-,  \\
      & \frac{1}{\beta^+ } \bfx +  \frac{1}{\alpha^+} n_2R_1(\bfx)R_2(\bfx) [(x_2-x_3),(x_3-x_1),(x_1-x_2)]^{\top}   ~~~~ \text{in}~ \Omega^+,
\end{cases}
\end{equation}
where $\bfx = [x_1,x_2,x_3]^{\top}$, $R_1(\bfx) = r^2_1 - \| \bfx \|^2$, $R_2(\bfx) = r^2_2 - \| \bfx \|^2$, and $n_1 = n_2(r^2_2 - r^2_1)$. We let $(\alpha^-,\alpha^+)=(\beta^-,\beta^+)=(1,100)$ and $(1,1000)$, and denote $\rho =\alpha^+/\alpha^- = \beta^+/\beta^- $. We present the numerical results for $N=10,...,70$ in Figure~\ref{fig:Example2Error} which clearly show the optimal convergence rates.

To solve the linear system, we employ the preconditioner described in Section~\ref{subsec:solver} with both GMRES and CG iteration strategies. Of course, it involves a direct solver on a smaller matrix of which the size depends on the expanding width $l$, {i.e., the size of $\bfB^{(2,2)}_l$ in \eqref{PG_IFE_system_mat}}. But the block diagonal smoother in \eqref{smoother} can help in greatly reducing the iteration steps for convergence, and consequently significantly reduce the computational time. {In fact, as $l$ is small, in practice we store the $LU$ factorization of $\bfB^{(2,2)}_l$ to avoid repeatly compute the inverse.} In Table~\ref{table:Example2NumItegmres}, we present both the number of iterations and computational time for GMRES. {It is highlighted that the preconditioning strategy even works for a CG method which in general does not converge for non-SPD matrices. The resuls are presented in Table~\ref{table:Example2NumIteCG}, and we can observe that the CG method must need $l=1$ at least for convergence. But once the block diagonal smoother in \eqref{smoother} is used, the CG method has better performance than GMRES. We may refer to our detailed discussions in Section~\ref{subsec:solver}, and point out again the important fact that the majority portion of $\bfB^{PG}$ is symmetric, and the problematic portion is only around the interface which can be well-handled by the direct inverse $(\bfB^{(2,2)}_l)^{-1}$.} Overall, for both methods we observe that the expanding width $l$ depends on the contrast $\rho$. For a larger contrast, there needs a relatively larger $l$ to ensure fast convergence. With an appropriately chosen expanding width, the algorithm converges in a small number of iterations which indeed takes significantly less computational time.

\begin{table}[H]
\begin{center}
\begin{tabular}{|c |c | c|c | c|c | c| c|}
\hline
\# DoFs   &         $7930$      & $59660$             & $197190$          & $462520$   & $897650$     & $1544580$   & $2445310$   \\ \hline
$\rho = 100$, $l=0$      &    63(2.63s)   &   66(14s)   &   66(29s)   &     67(51s)     &    68(108s)   &     68(191s)   &   68(283s)           \\ \hline
$\rho = 100$, $l=1$      &    35(2.11s)   &   40(7.89s)   &   39(14s)   &     40(28s)     &    40(56s)   &     38(100s)   &   41(147s)           \\ \hline
$\rho = 1000$, $l=0$    &    212(8.87s)    &  244(60s)      &   263(117s)      &     223(185s)      &     246(404s)    &   272(817s)      &    234(1017s)           \\ \hline
$\rho = 1000$, $l=1$    &    55(8.44s)   &   63(13.67s)   &   76(33s)  &      48(35s)    &     46(70s)   &     52(143s)    &  54(243s)          \\ \hline
$\rho = 1000$, $l=2$    &    42(1.9s)   &   47(11.23s)   &   46(21s)   &     52(39s)     &    46(83s)   &     51(144s)   &   51(209s)           \\ \hline
$\rho = 1000$, $l=3$    &    41(32s)   &   39(11s)   &   42(20s)   &     39(40s)     &    42(75s)   &     42(135s)   &   42(191s)           \\ \hline
\end{tabular}
\end{center}
\caption{Iteration numbers and computational time for the GMRES solver.}
\label{table:Example2NumItegmres}
\end{table}

\begin{table}[H]
\begin{center}
\begin{tabular}{|c |c | c|c | c|c | c| c|}
\hline
\# DoFs   &         $7930$      & $59660$             & $197190$          & $462520$   & $897650$     & $1544580$   & $2445310$   \\ \hline
$\rho = 100$, $l=0$      &    --   &   84(18s)   &   78(29s)   &     --     &    --   &     --   &   --           \\ \hline
$\rho = 100$, $l=1$      &    37(1.6s)   &   39(9.8s)   &   40(18s)   &     41(30s)     &    42(56s)   &     42(99s)   &   43(140s)           \\ \hline
$\rho = 100$, $l=2$      &    36(1.3s)   &   39(10s)   &   41(17s)   &     42(32s)     &    42(66s)   &     42(110s)   &   42(170s)           \\ \hline
$\rho = 1000$, $l=0$    &    --   &   --   &   --   &     --     &    --   &     --    &  --          \\ \hline
$\rho = 1000$, $l=1$    &    40(1.4s)   &   42(17s)   &   42(30s)   &     42(30s)     &    43(59s)   &     45(100s)    &  43(140s)           \\ \hline
$\rho = 1000$, $l=2$    &    37(1.4s)   &   39(8.4s)   &   40(18s)   &     41(29s)     &    42(61s)  &   44(110s)   &   42(150s)           \\ \hline
\end{tabular}
\end{center}
\caption{Iteration numbers and computational time for the CG solver.}
\label{table:Example2NumIteCG}
\end{table}


\subsection{Example 3 (Application to a time-domain Maxwell equation)}

In this test, we consider the scattering problem modeled by \eqref{Max_mode1} in a time domain $[0,T]$. We partition it into $M$ sub-intervals: $t_0=0<t_1<t_2<\cdots<t_M=T$ with equal length $\tau = T/M$. We then apply the implicit time-discretization introduced in~\cite{1999CiarletZou} which reads
\begin{equation}
\label{timeScheme_1}
(\epsilon \partial^2_{\tau}\bfu^n_h, \bfv_h  )_{\Omega} + (\sigma \partial_{\tau}\bfu^n_h, \bfv_h)_{\Omega} + (\mu^{-1}\curl \bfu^n_h, \curl\bfv_h)_{\Omega} = (\bfJ_t,\bfv_h)_{\Omega}, ~~~~ \forall \bfv_h\in \bfS^e_h,
\end{equation}
where, for $n\ge 2$, $\partial_{\tau}\bfu^n_h = (\bfu^{n}_h - \bfu^{n-1}_h)/\tau$ and  $ \partial^2_{\tau}\bfu^n_h = (\partial_{\tau}\bfu^{n}_h - \partial_{\tau } \bfu^{n-1}_h)/\tau$.
Then, at each time step it can be discretized by the PG-IFE scheme with the parameters $\alpha = \mu^{-1}$ and $\beta = \epsilon \tau^{-2} + \sigma \tau^{-1}$ and solved with the proposed fast solver. The optimal convergence rates can be also observed if the analytical solution is available, which is omitted for avoiding redundancy. Here, we consider a simulation problem of electromagnetic waves propagating through a torus where the analytical solution is difficult to derive if not impossible, but instead we compare the results of the IFE method and the standard FE method on fitted meshes. 

The torus interface shown in the right plot of Figure~\ref{fig:interfexample} is described by a level-set signed-distance function:
$$
\gamma(\bfx) = (\sqrt{(x_1- z_1)^2+(x_2 - z_2)^2}-r_2)^2 + (x_3 - z_3)^2 - r_1^2,
$$
where $z_1  = z_2= 0$, $z_3 = -0.3$ and $r_1=0.2$, $r_2=\pi/5$. We further set the medium parameters as $\epsilon^+=0.05$, $\epsilon^- = 2\epsilon^+$, $\sigma^+=0.1$, $\sigma^-=1$ and $\mu^+=4\pi$, $\mu^-=3\mu^+$. Initially, we impose a Gaussian pulse
\begin{equation}
\label{GSPulse}
\bfu = \exp(-b(a(x_1-z_0)-\omega t)^2 )[0,1,0]^{\top}
\end{equation}
where $a=\omega\sqrt{\epsilon^+\mu^+}$ modeling the speed of the electromagnetic wave. In addition, the boundary condition is set to match \eqref{GSPulse} except on the faces at $\bfx=-1,1$ where it is set as $\mathbf{ 0}$. Set the source as $\bfJ=\mathbf{ 0}$, let the time domain be $[0,1.5]$ and let the step size be $\tau = 1.5/128$. The IFE method is computed on a $64^3$ unfitted mesh which contains 1572864 elements and 1872064 DoFs, while the simulation of the FE method is performed on a fitted mesh containing 1679542 elements and 2002947 DoFs. 

\begin{figure}[htbp]
\begin{tabular}{ >{\centering\arraybackslash}m{1.3in}>{\centering\arraybackslash}m{1.3in} >{\centering\arraybackslash}m{1.3in}  >{\centering\arraybackslash}m{1.3in}  }
	\centering
	$t=12\tau$ &
	$t=24\tau$ &
	$t=36\tau$ &
	$t=48\tau$ \\
	\includegraphics[width=1.3in]{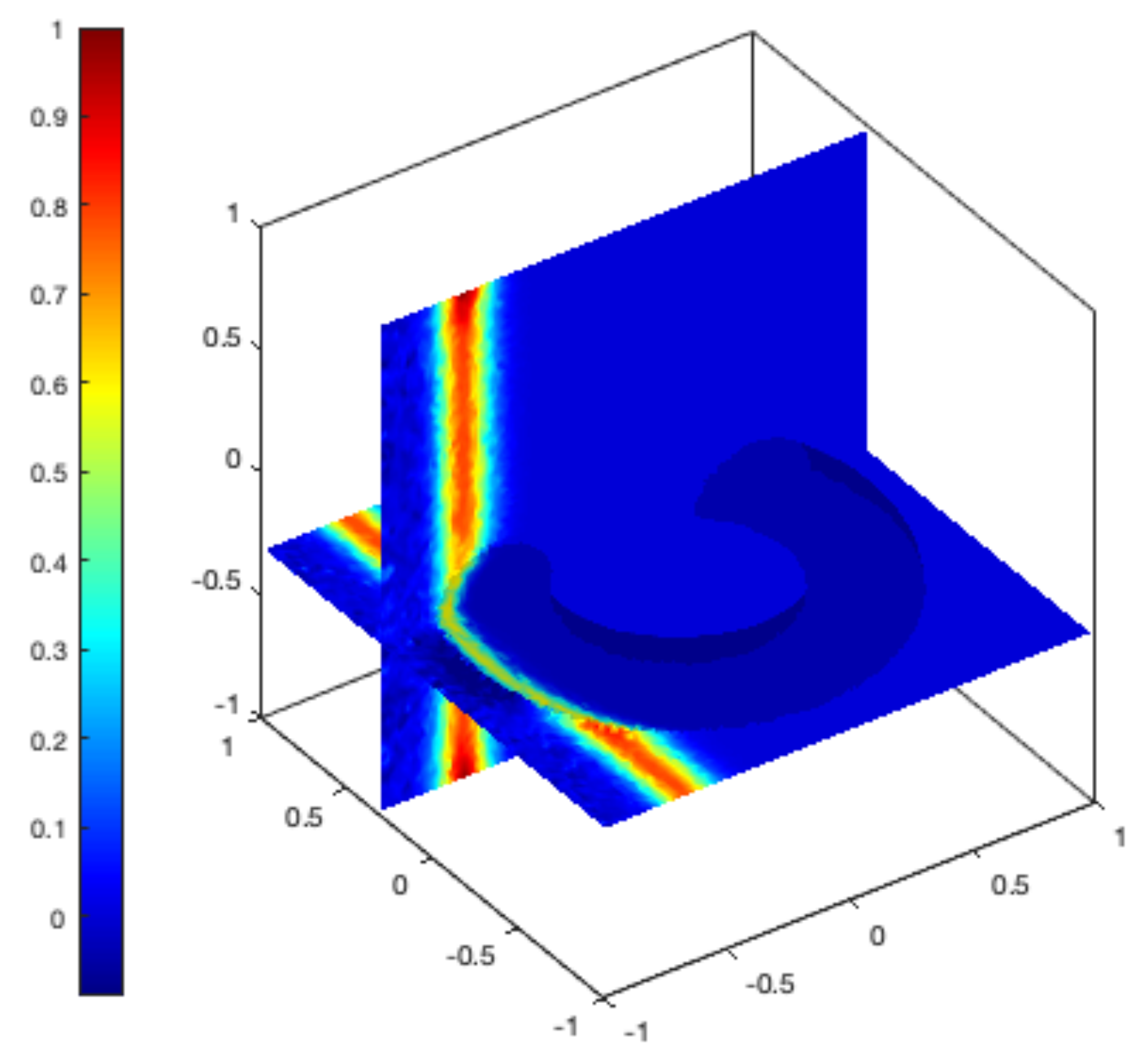}&
	\includegraphics[width=1.2in]{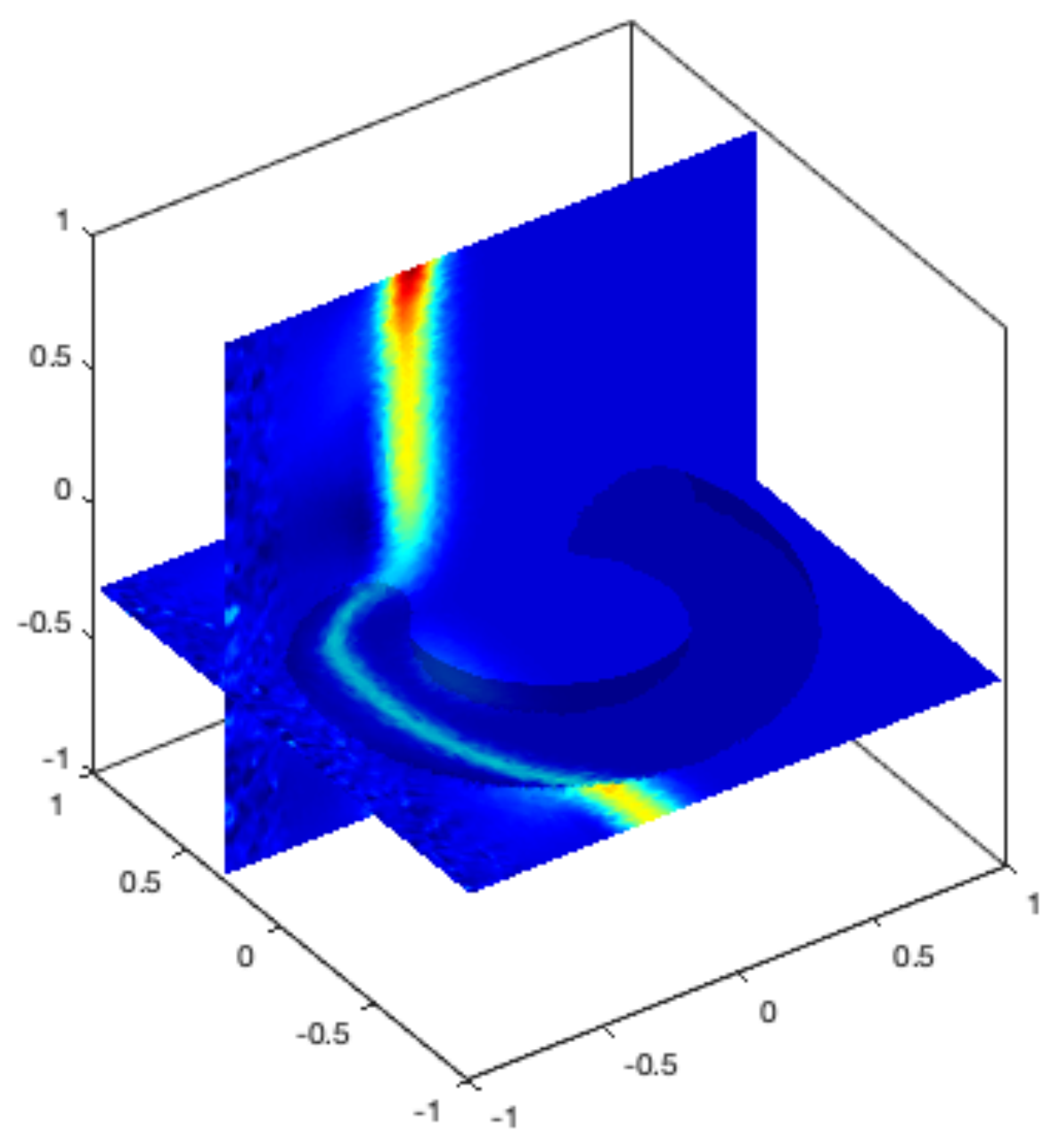}&
	\includegraphics[width=1.2in]{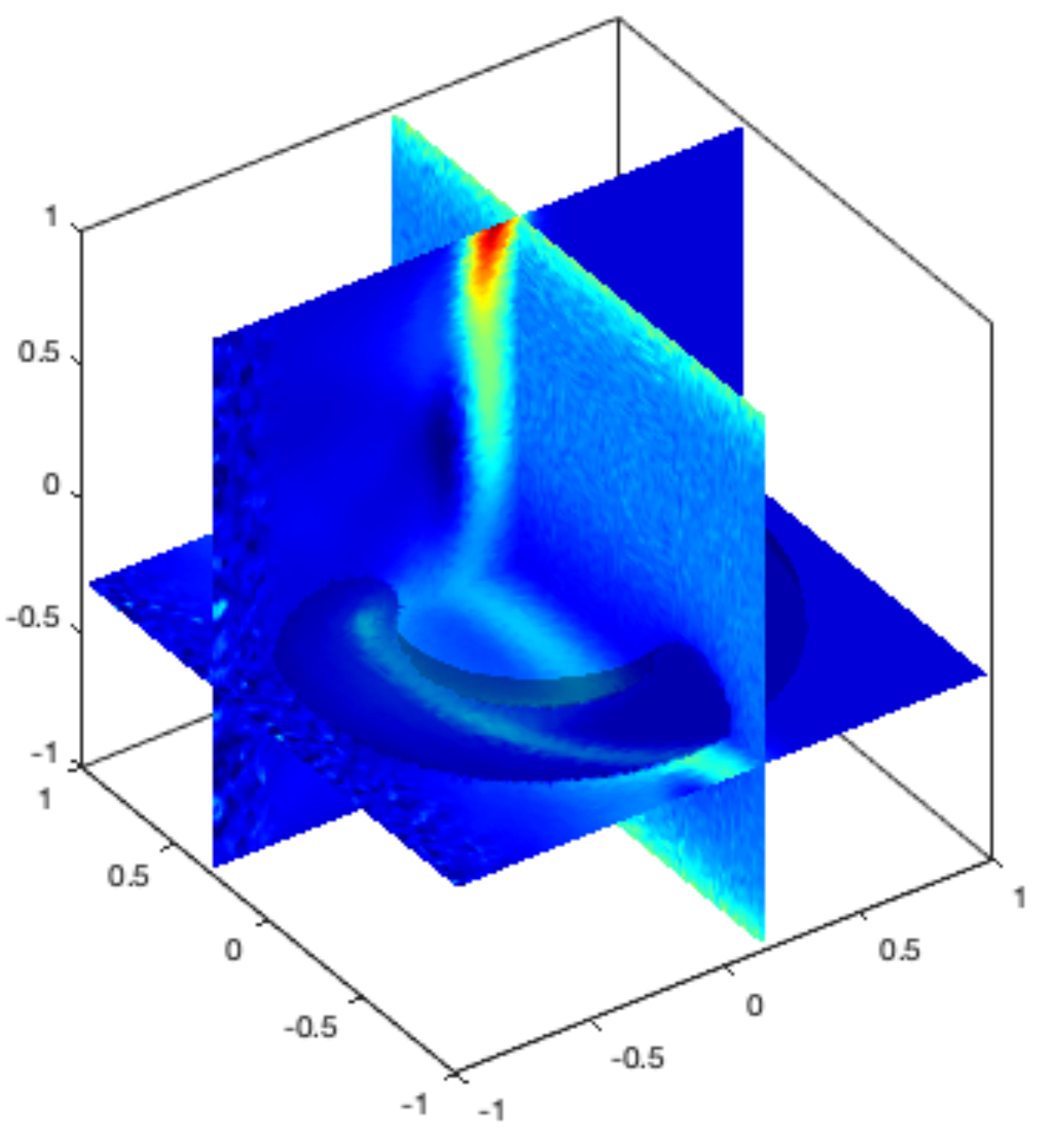}&
	\includegraphics[width=1.2in]{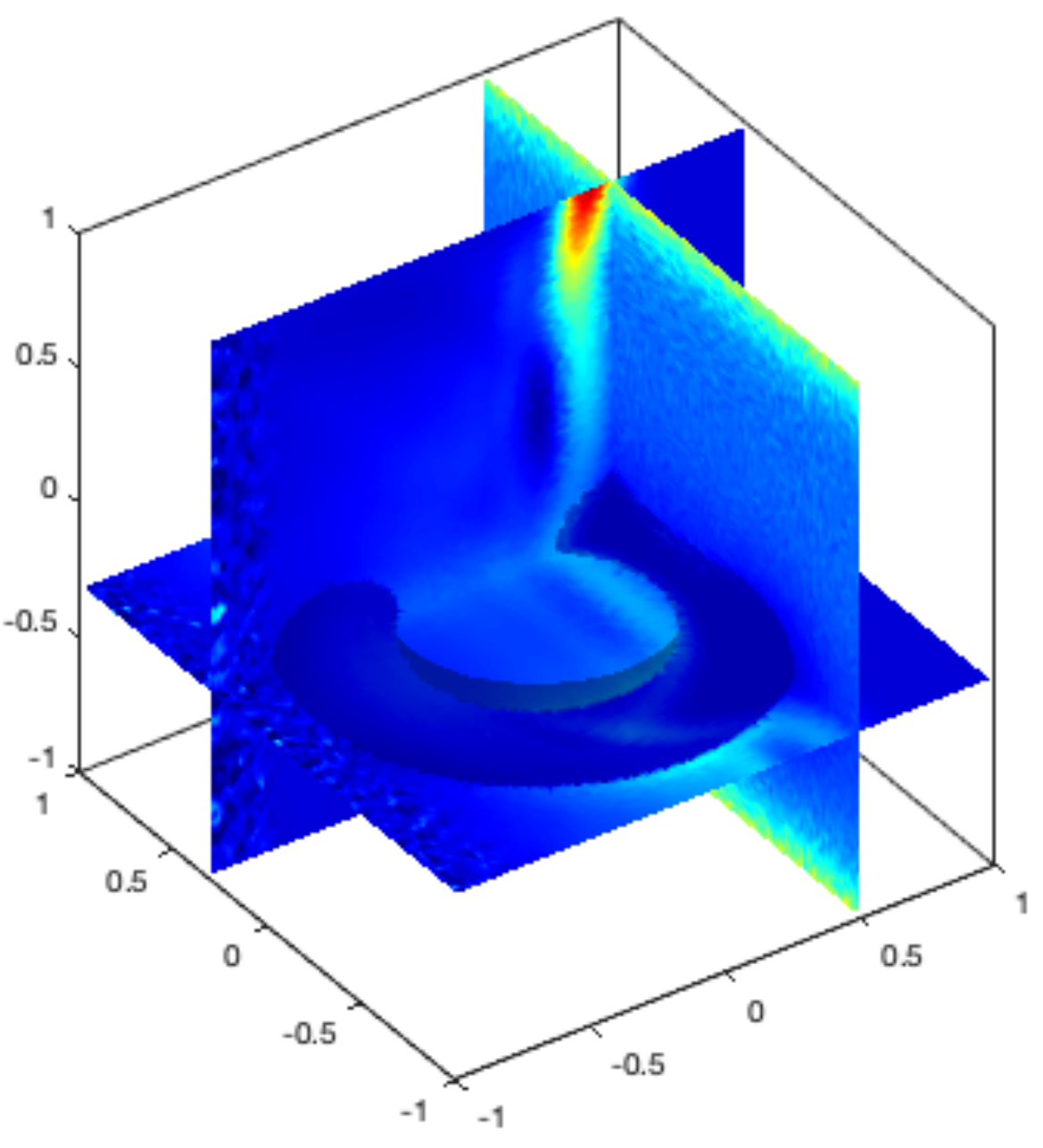}\\
	\includegraphics[width=1.3in]{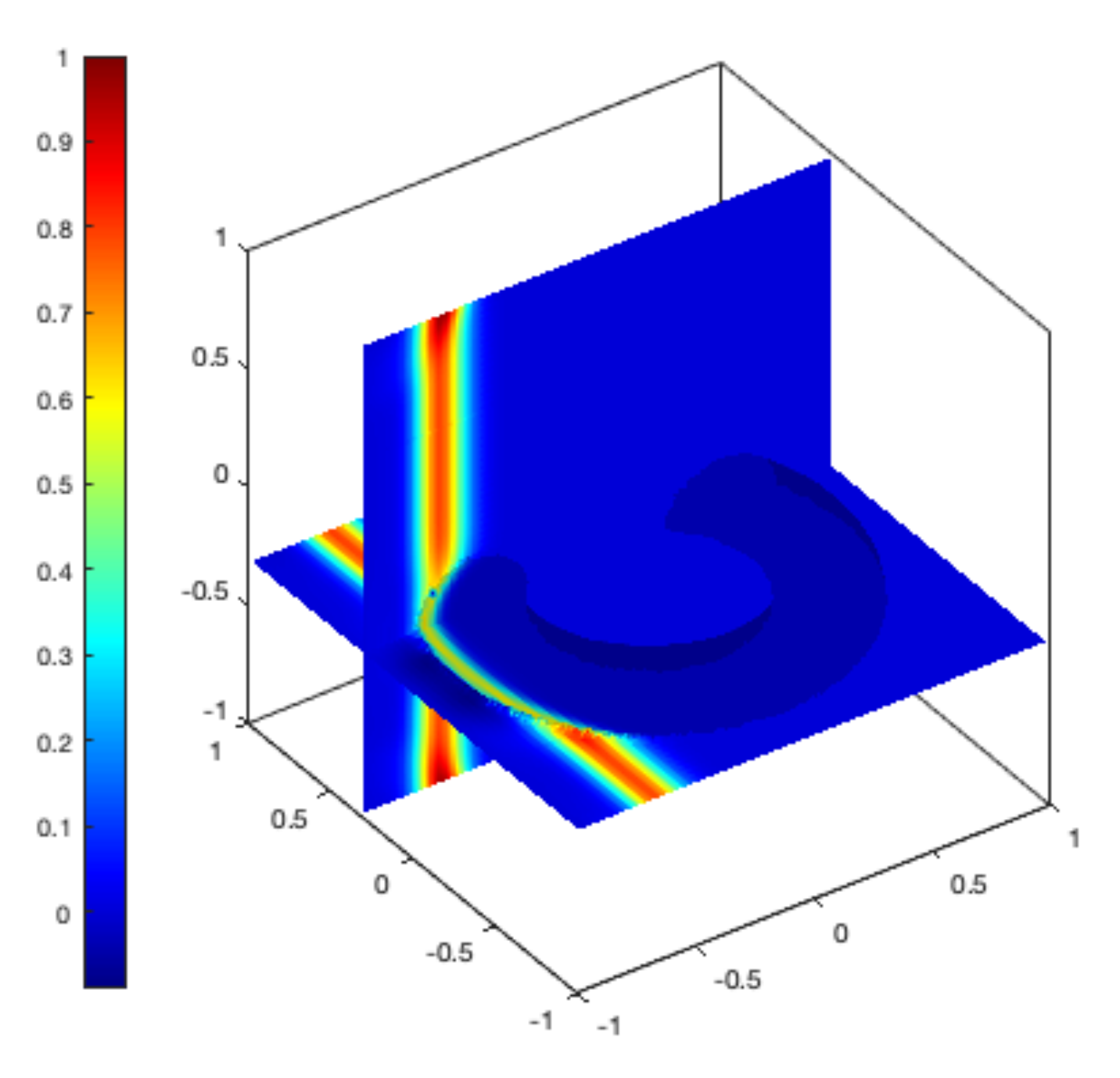}&
	\includegraphics[width=1.2in]{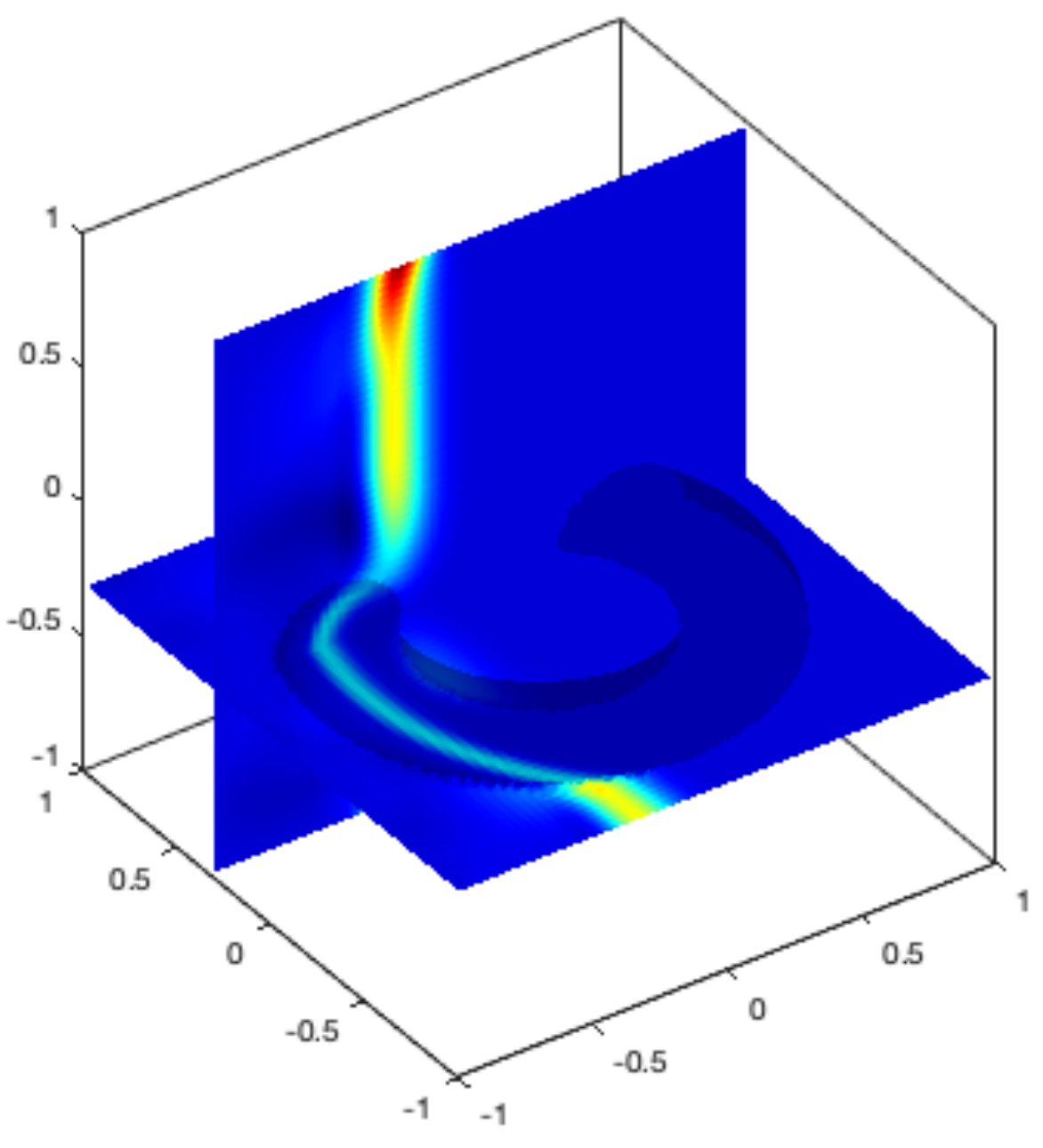}&
	\includegraphics[width=1.2in]{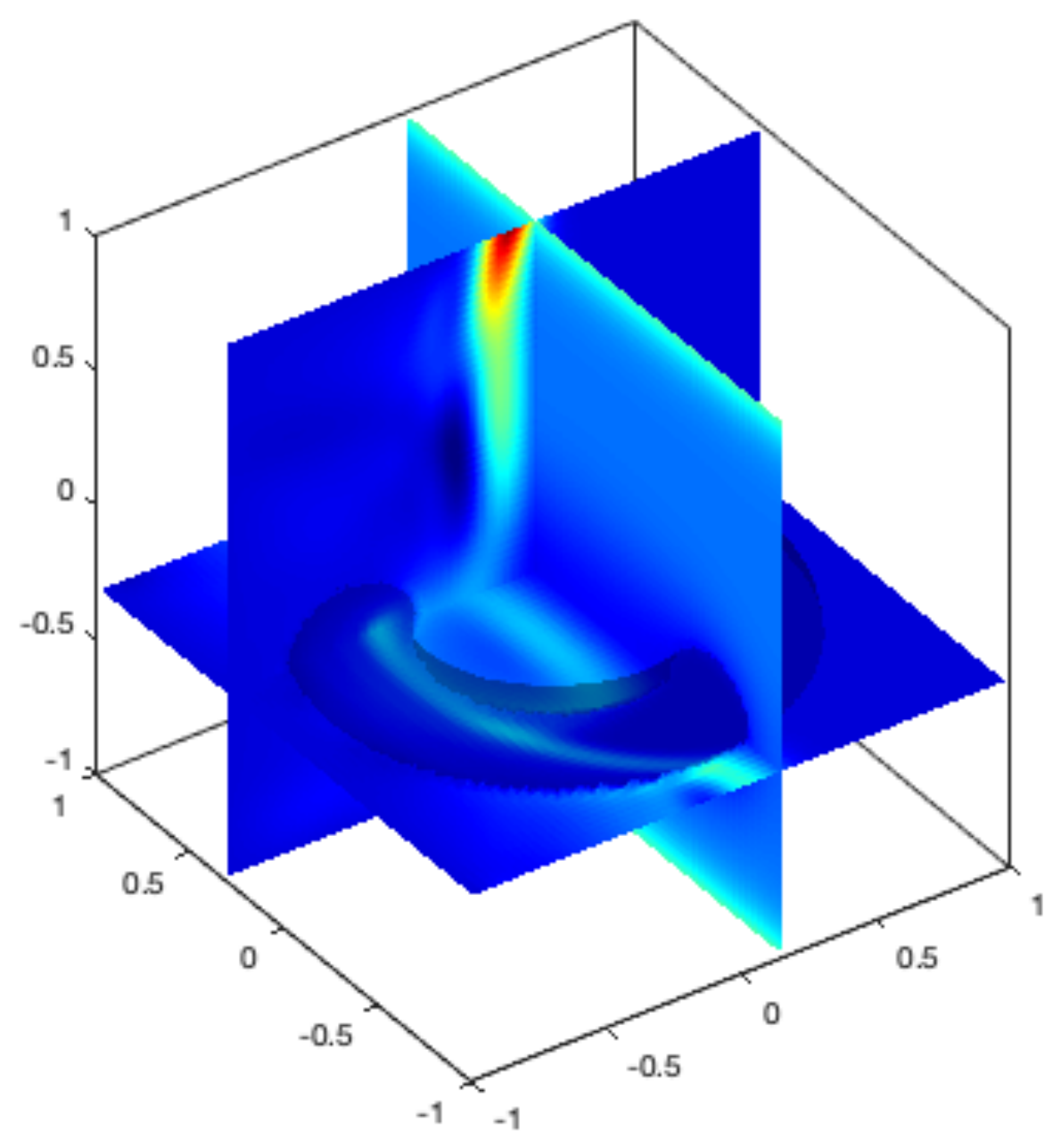}&
	\includegraphics[width=1.2in]{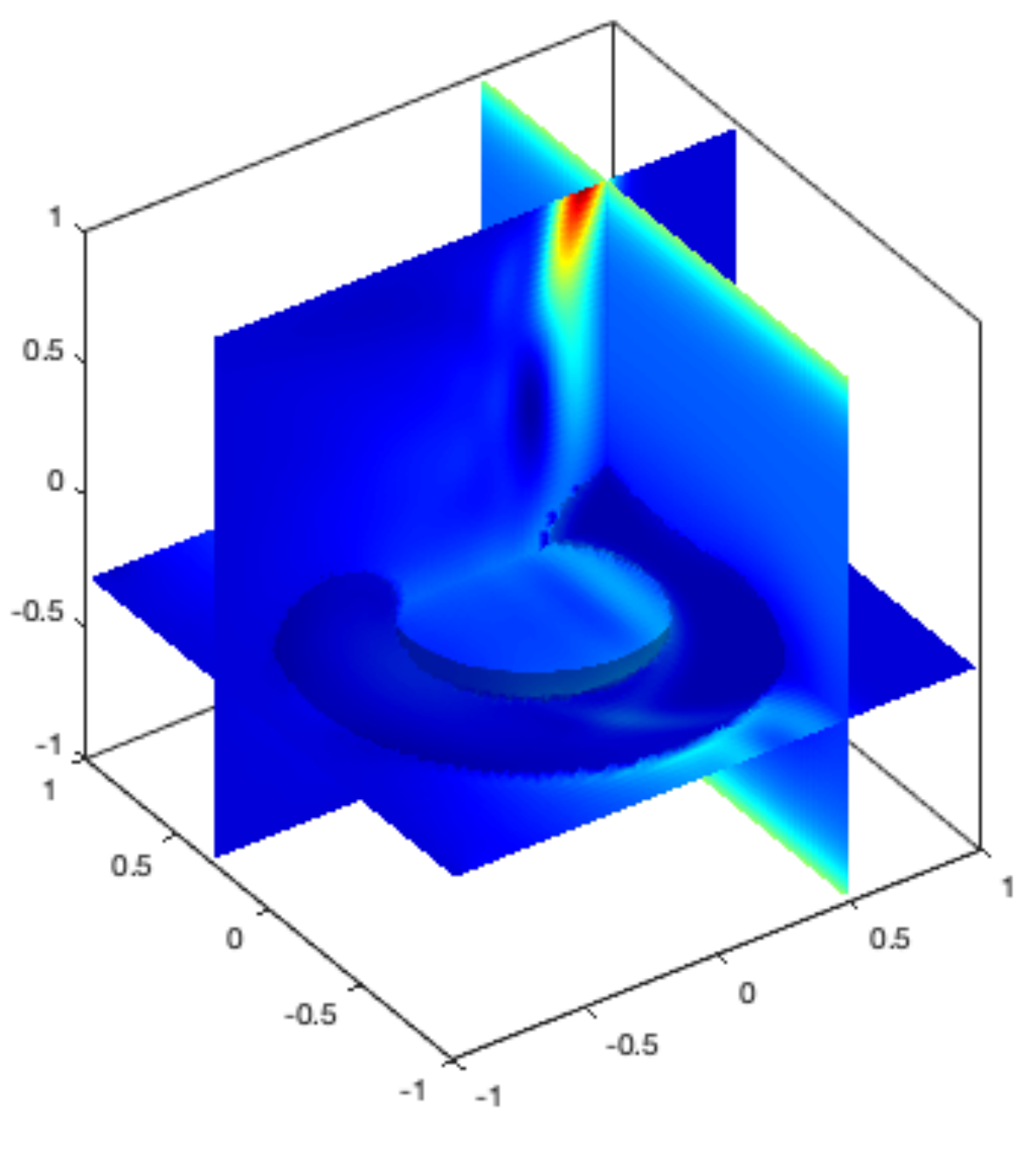}\\
\end{tabular}
\caption{Comparison of FE (top) and IFE (bottom) solutions at slices $x_1=5/32$, $x_2=5/16$, $x_3=-5/16$.} 
\label{example3slice}
\end{figure}

\begin{figure}[htbp]
\begin{tabular}{ >{\centering\arraybackslash}m{1.3in}>{\centering\arraybackslash}m{1.3in} >{\centering\arraybackslash}m{1.3in}  >{\centering\arraybackslash}m{1.3in}  }
	\centering
	$t=12\tau$ &
	$t=24\tau$ &
	$t=36\tau$ &
	$t=48\tau$ \\
	\includegraphics[width=1.2in]{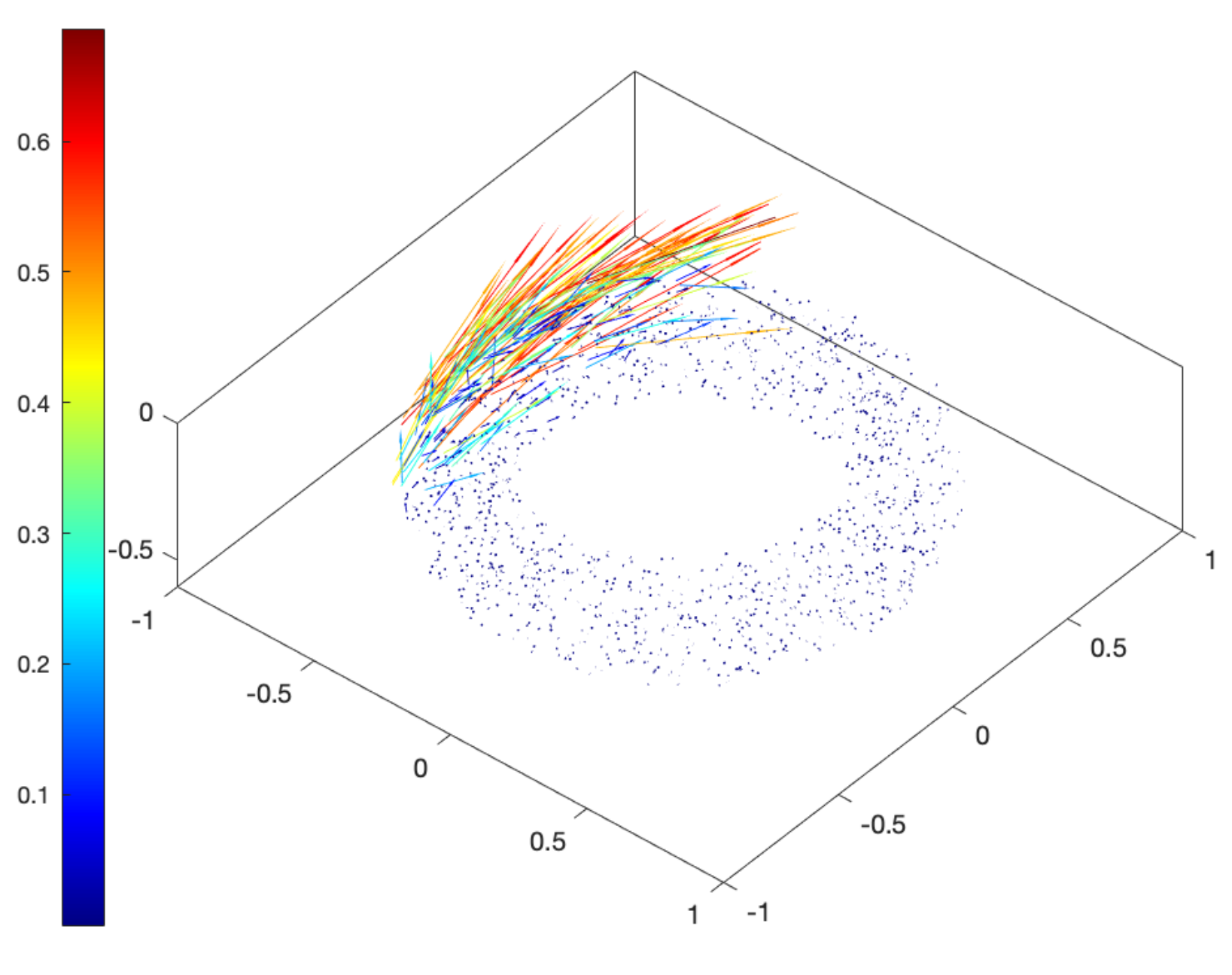}&
	\includegraphics[width=1.2in]{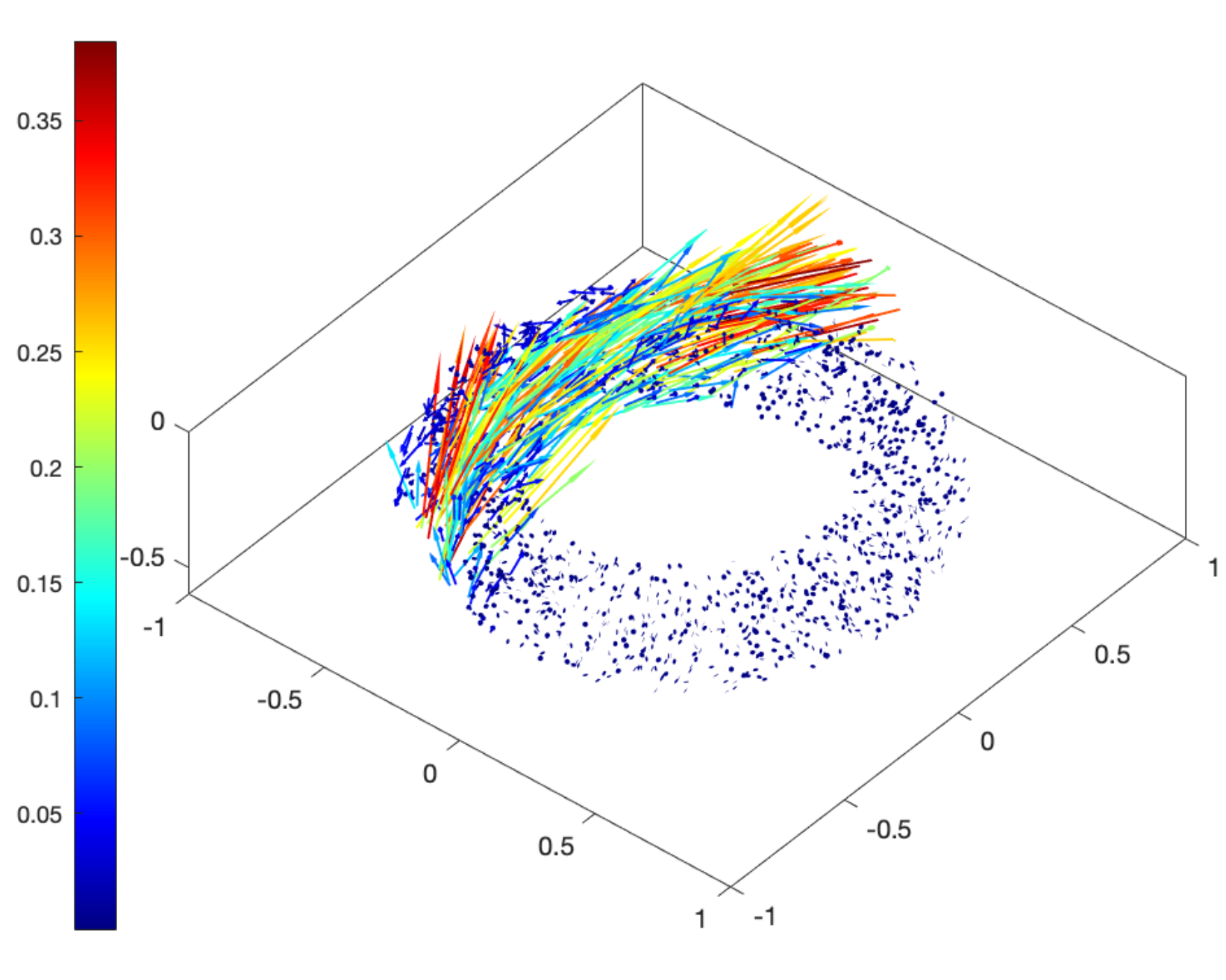}&
	\includegraphics[width=1.2in]{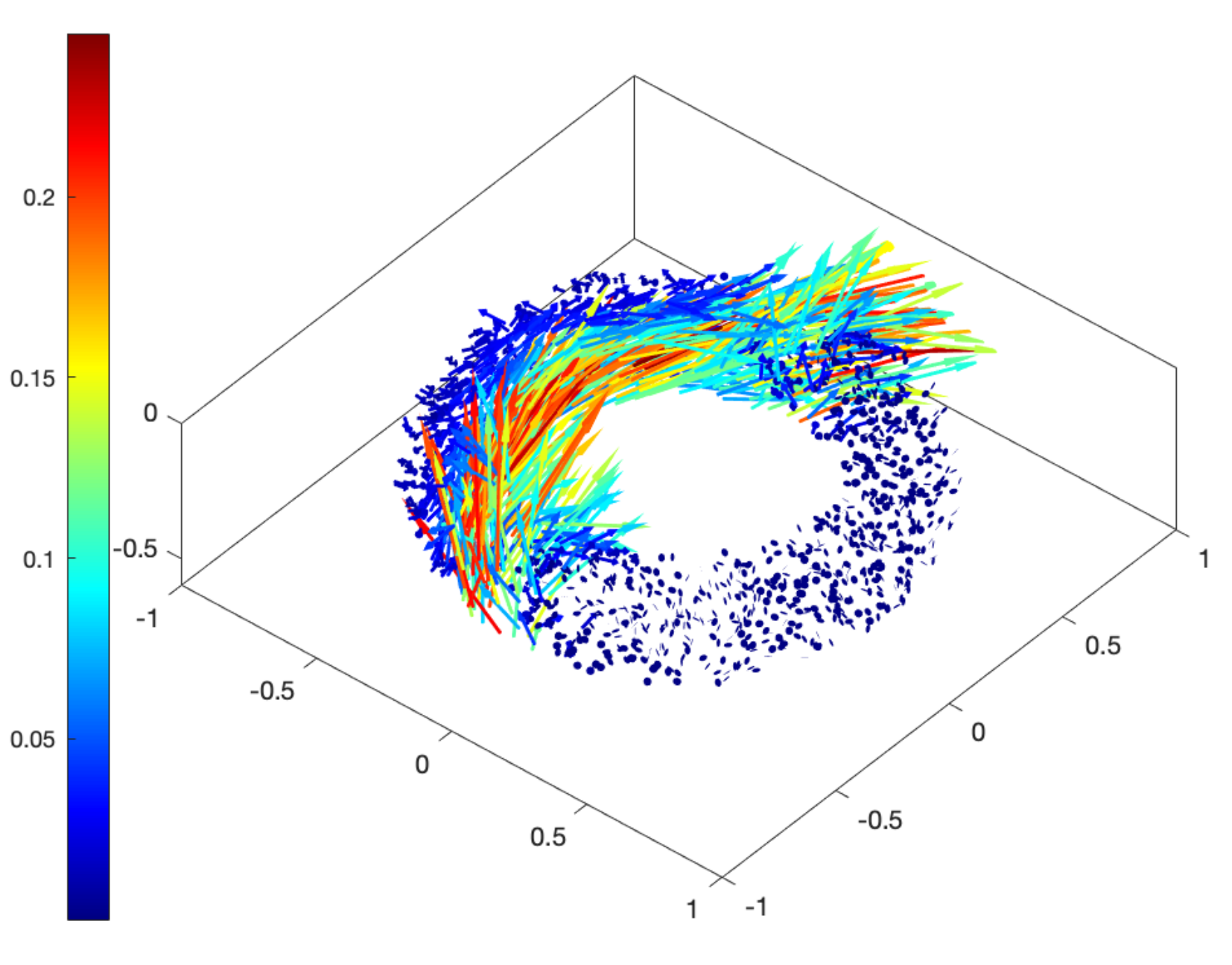}&
	\includegraphics[width=1.2in]{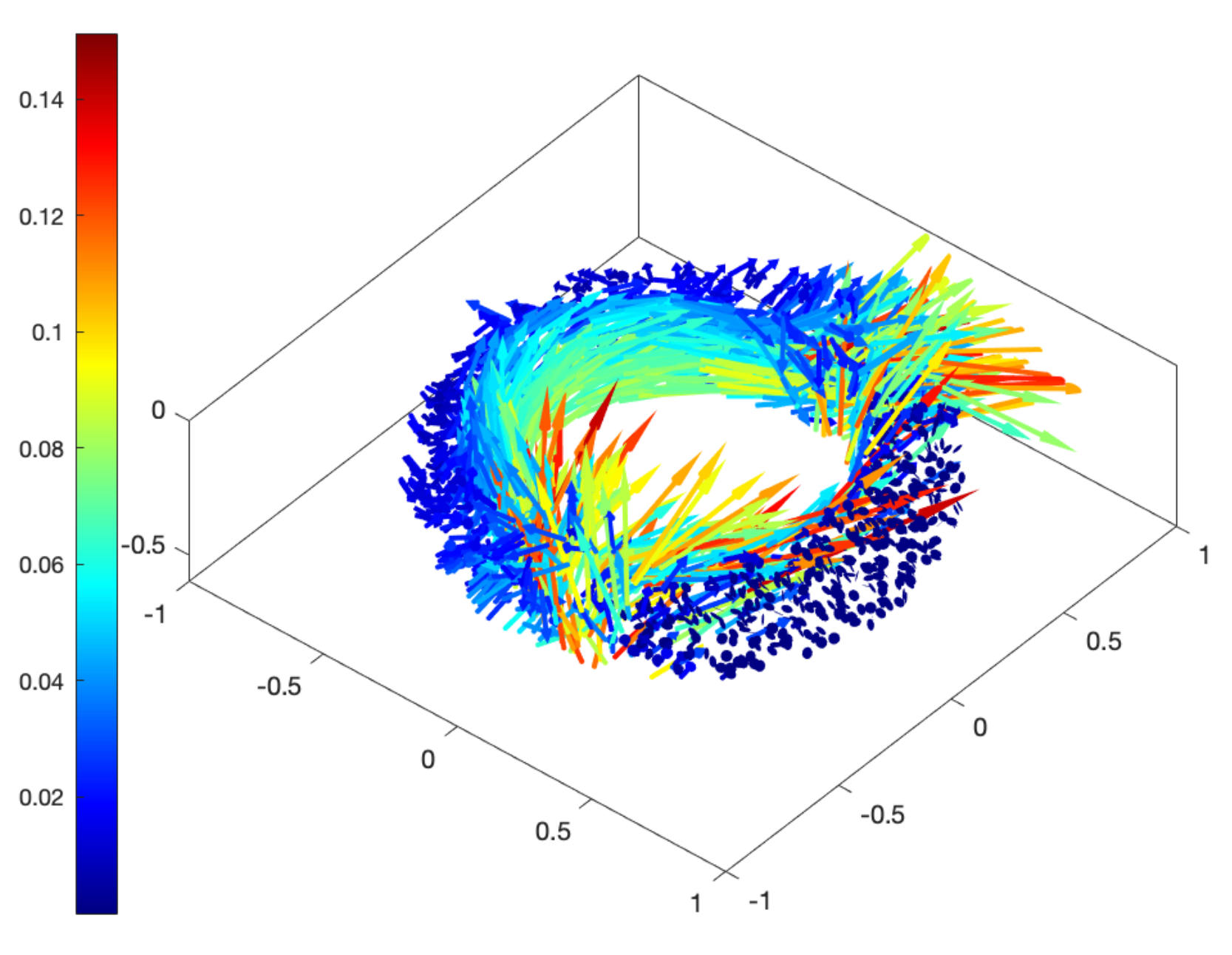}\\
	\includegraphics[width=1.2in]{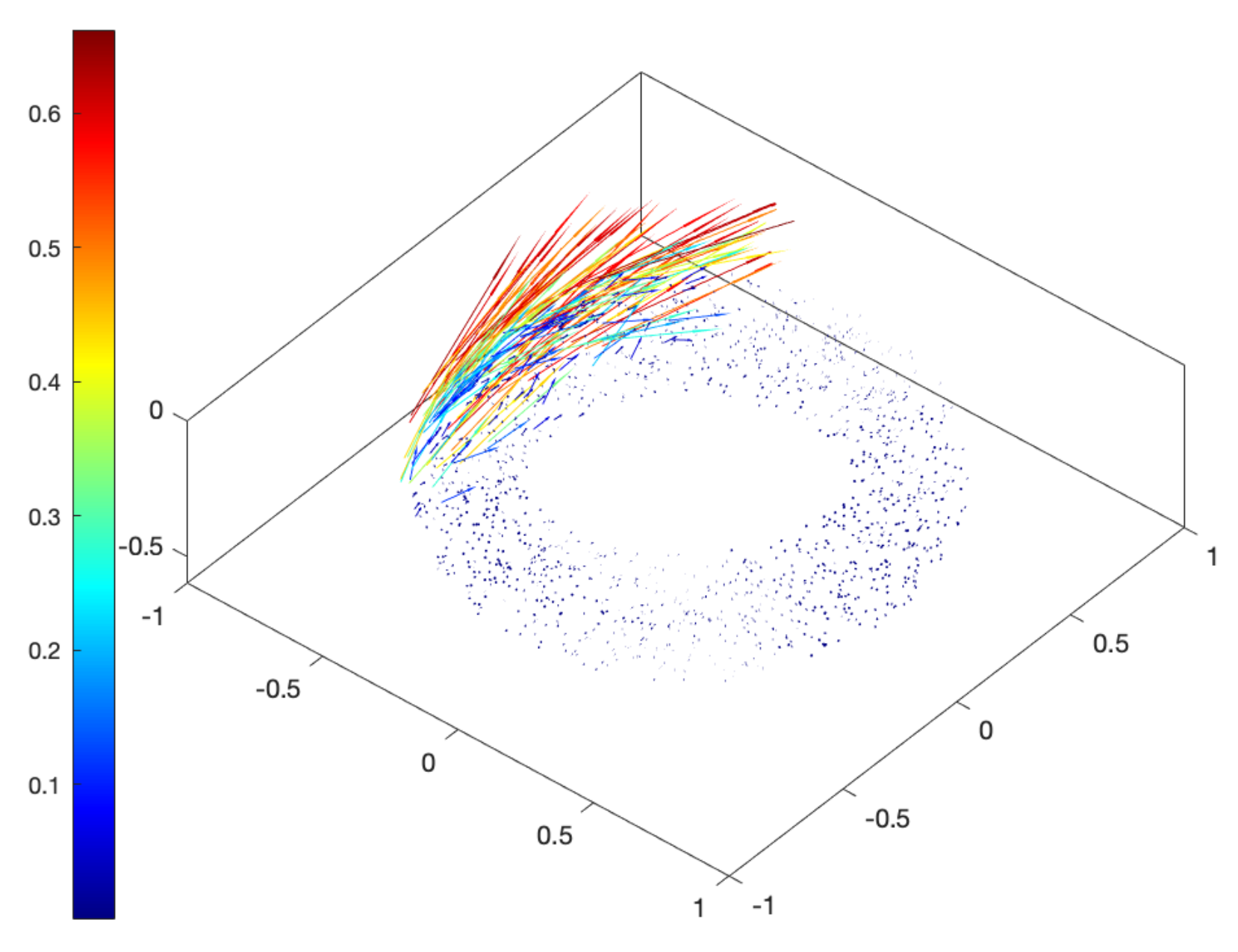}&
	\includegraphics[width=1.2in]{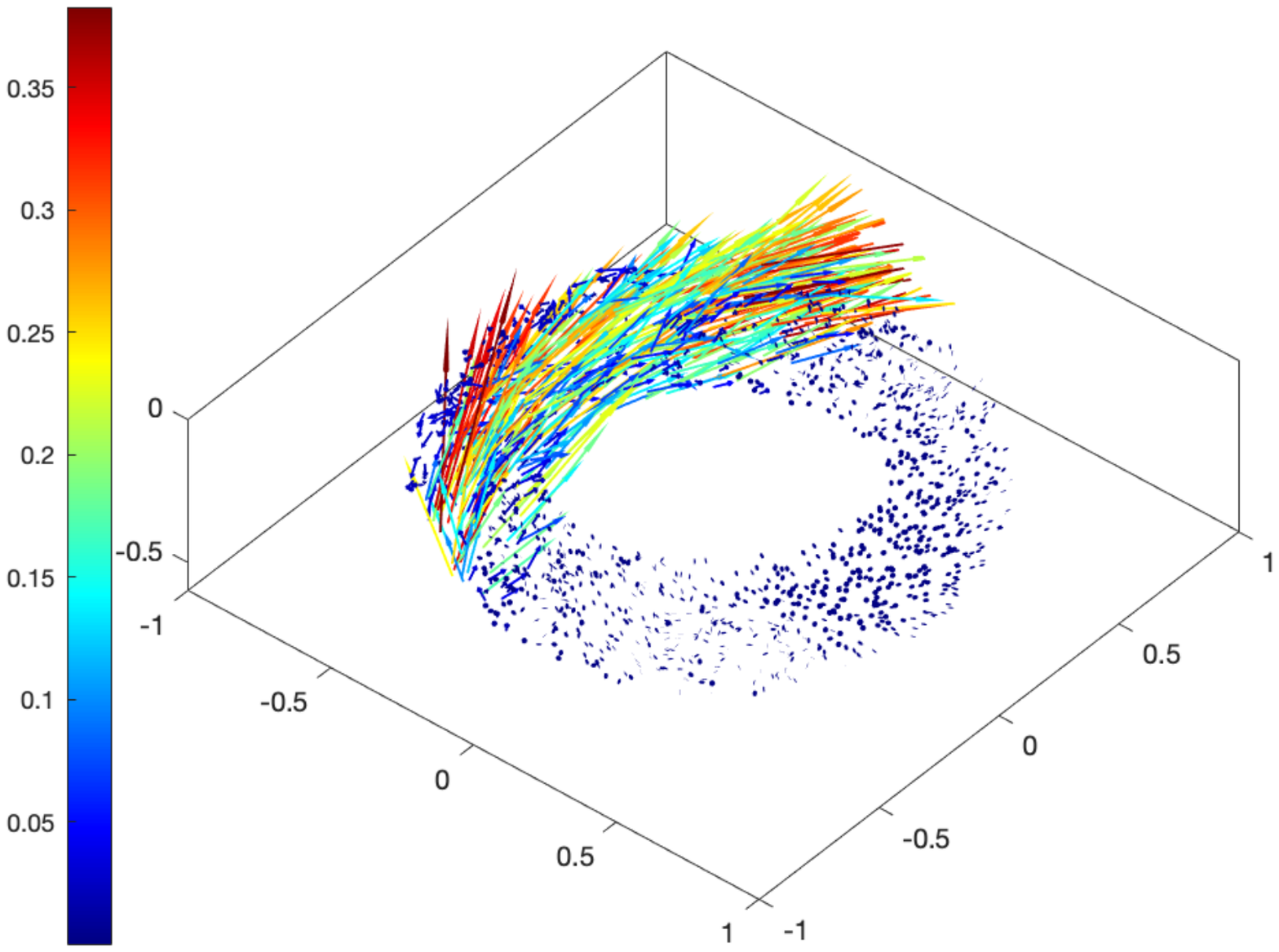}&
	\includegraphics[width=1.2in]{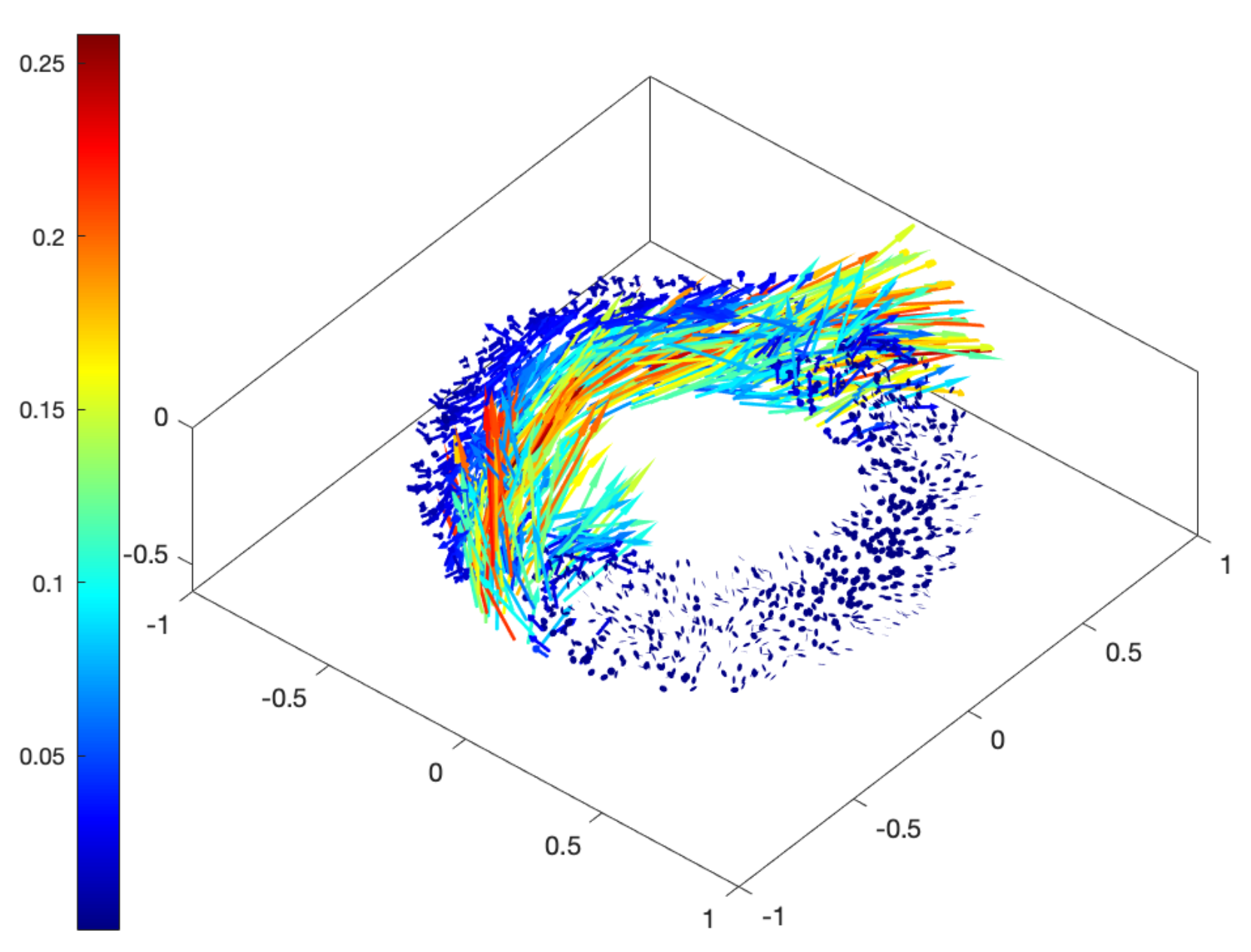}&
	\includegraphics[width=1.2in]{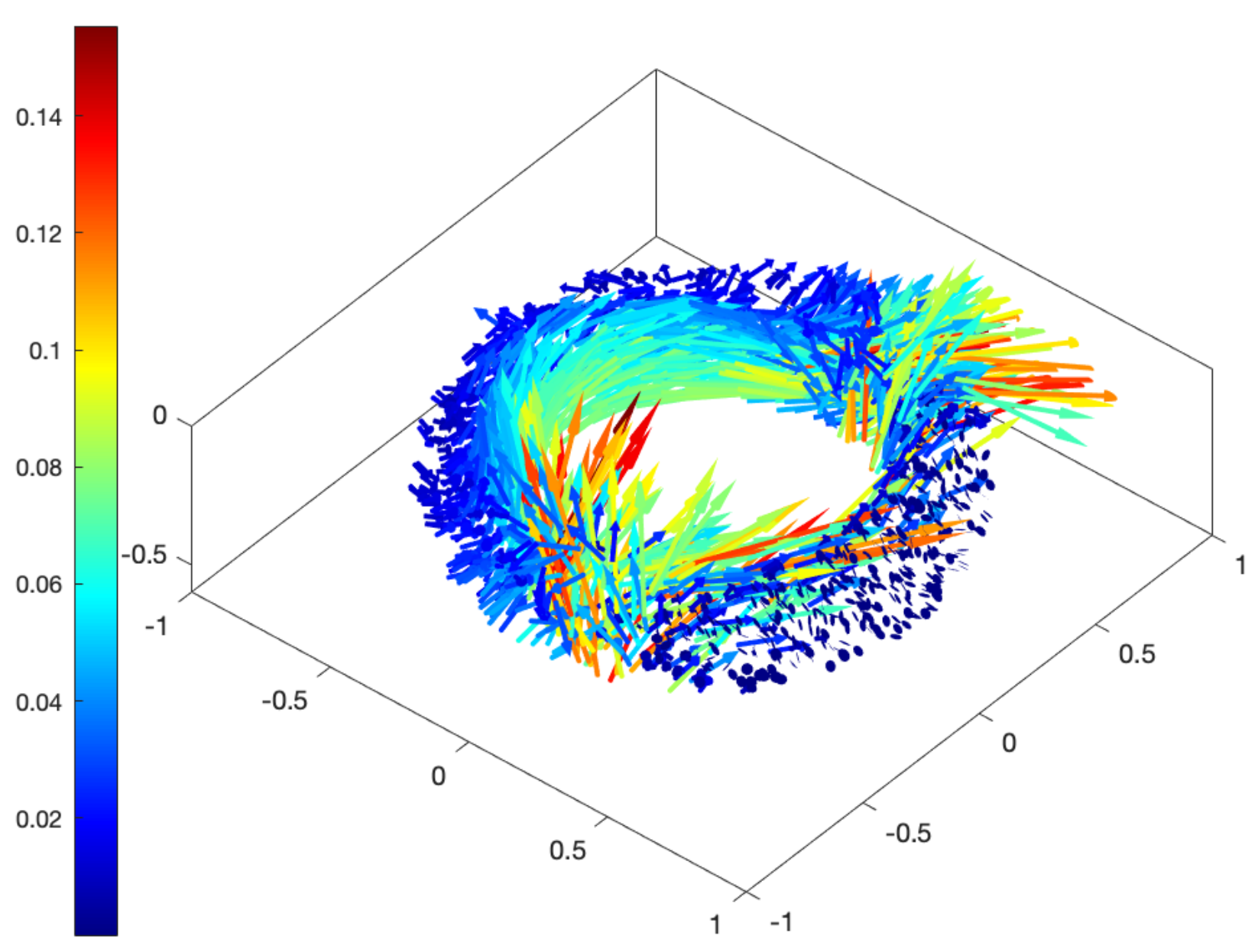}\\
\end{tabular}
\caption{Comparison of FE (top) and IFE (bottom) solutions of vector fields.} 
\label{example3torus}
\end{figure}

Denote the numerical solutions by $\bfu^{IFE}_h$ and $\bfu^{FE}_h$, respectively. To compare the simulation results, we first plot the evolution of the $x_2$ component at some slices in Figure~\ref{example3slice} at $t=12\tau$, $t=24\tau$, $t=36\tau$ and $t=48\tau$. The results of the IFE and FE methods are almost identical demonstrating the effectiveness of the proposed method. In addition, we can clearly observe the delayed electromagnetic wave inside the torus agreeing with the physics. Next, we plot the electromagnetic fields inside the torus in Figure~\ref{example3torus} such that we can also examine the direction. The plots still clearly show that the two methods have very similar simulation results.

\begin{appendices}
\section{Proof of Lemma~\ref{lemma_interface_flat}}
Here, we only show \eqref{lemma_interface_flat_eq02}, since it can yield \eqref{lemma_interface_flat_eq01} by the argument in~\cite{2020GuoLin}. Since $\varphi^{\Gamma}$ is sufficiently smooth, we obtain by the interpolation approximation that
\begin{equation}
\label{lemma_interface_flat_eq1}
\| \nabla \varphi^{\Gamma} - \nabla \varphi^{\Gamma}_h \|_{L^{\infty}(K)} \le C h_K
\end{equation}
where $C$ depends on the second derivative of $\varphi^{\Gamma} $. As a signed-distance function, $\varphi^{\Gamma}$ satisfies $\|\nabla \varphi^{\Gamma} \|=1$, and thus \eqref{lemma_interface_flat_eq1} implies
\begin{equation}
\label{lemma_interface_flat_eq2}
\verti{ 1 - \| \nabla \varphi^{\Gamma}_h(\bfx) \| } \le Ch_K, ~~~ \forall \bfx\in K.
\end{equation}
Note that the unit normal vectors of $\Gamma$ and $\Gamma^K_h$ are $\nabla \varphi^{\Gamma}$ and $ \nabla \varphi^{\Gamma}_h/\|  \nabla \varphi^{\Gamma}_h \|$, respectively. 
Then, for sufficiently small $h$, \eqref{lemma_interface_flat_eq1} and \eqref{lemma_interface_flat_eq2} together yield
$$
\| \nabla \varphi^{\Gamma} - \nabla \varphi^{\Gamma}_h/ \| \nabla \varphi^{\Gamma}_h \| \|_{L^{\infty}(K)} \le \frac{\| \nabla \varphi^{\Gamma} - \nabla \varphi^{\Gamma}_h + \nabla \varphi^{\Gamma} (\| \nabla \varphi^{\Gamma}_h \| -1)  \|_{L^{\infty}(K)} }{\min_{\bfx\in K}\| \nabla \varphi^{\Gamma}_h(\bfx) \|} \lesssim h_K
$$
which has finished the proof.


  \section{Proof of Lemma~\ref{lem_diff_u}}
  \label{appen:sec:lem_diff_u}
Here, we only prove \eqref{lem_function_gamma_eq01} with $\dd=\curl$, i.e., $a=1$, and the others are similar. Let $\bfv = \bfu^-_E\times \bar{\bfn}_K - \bfu^+_E\times \bar{\bfn}_K$. Without loss of generality, we define a local system $(\hat{x}_1,\hat{x}_2,\hat{x}_3)$ with $\hat{x}_3$ being the direction normal to $\Gamma^{\omega_K}_h$, i.e., $\bar{\bfn}_K$, and $\hat{x}_1,\hat{x}_2$ spanning the plane containing $\Gamma^{\omega_K}_h$ such that the interface surface can be represented by $\hat{x}_3=f(\hat{x}_1,\hat{x}_2)$ with $f$ being a smooth function. Then,
\begin{equation}
\label{lem_diff_u_eq1}
\bfv(\hat{x}_1,\hat{x}_2,f(\hat{x}_1,\hat{x}_2)) = \bfv(\hat{x}_1,\hat{x}_2,0) + \int_{0}^{f(\hat{x}_1,\hat{x}_2)} \partial_{\hat{x}_3} \bfv(\hat{x}_1,\hat{x}_2,\hat{x}_3) \dd \hat{x}_3
\end{equation}
which, by H\"older's inequality, yields
\begin{equation}
\begin{split}
\label{lem_diff_u_eq2}
\| \bfv \|^2_{L^2(\Gamma^{\omega_K}_h)} & =  \int_{\Gamma^{\omega_K}_h} \|\bfv(\hat{x}_1,\hat{x}_2,0)\|^2 \dd \hat{x}_1\dd \hat{x}_2 \lesssim \int_{\Gamma^{\omega_K}_h} \|\bfv(\hat{x}_1,\hat{x}_2,f(\hat{x}_1,\hat{x}_2))\|^2 \dd \hat{x}_1\dd \hat{x}_2\\
& + \int_{\Gamma^{\omega_K}_h} \int_0^{f(\hat{x}_1,\hat{x}_2)} \|\partial_{\hat{x}_3} \bfv(\hat{x}_1,\hat{x}_2,f(\hat{x}_1,\hat{x}_2)) \|^2 \dd \hat{x}_3 |f(\hat{x}_1,\hat{x}_2)| \dd \hat{x}_1\dd \hat{x}_2
\end{split}
\end{equation}
of which the last two integrals are denoted by $({\rm I})$ and $({\rm II})$, respectively. With \eqref{lem_geo_gamma_eq0}, it is clear that
\begin{equation}
\label{lem_diff_u_eq3}
({\rm II}) \le \max_{(\hat{x}_1,\hat{x}_2)\in \Gamma^{\omega_K}_h}|f(\hat{x}_1,\hat{x}_2)| \| \partial_{\hat{x}_3} \bfv \|^2_{L^2(\omega_K)} \lesssim h^2_K\| \nabla \bfv \|^2_{L^2(\omega_K)} .
\end{equation}
As for $({\rm I})$, we note that $\sqrt{1+ \| \nabla f \|^2} \ge 1$, and thus obtain
\begin{equation}
\begin{split}
\label{lem_diff_u_eq4}
({\rm I}) & \le \int_{\Gamma^{\omega_K}_h} \|\bfv(\hat{x}_1,\hat{x}_2,f(\hat{x}_1,\hat{x}_2))\|^2 \sqrt{1+ \| \nabla f \|^2} \dd \hat{x}_1\dd \hat{x}_2 = \int_{\Gamma^{\omega_K}} \| \bfv \|^2 \dd s = \| \bfv \|^2_{L^2(\Gamma^{\omega_K})} \\
 & = \| (\bfu^-_E - \bfu^+_E)\times (\bfn-\bar{\bfn}_K) \|^2_{L^2(\Gamma^{\omega_K})} \lesssim \max_{\Gamma^{\omega_K}} \| \bfn - \bar{\bfn}_K \|^2 \| \bfu^-_E - \bfu^+_E \|^2_{L^2(\Gamma^{\omega_K})} \\
 & \lesssim h^2_K \| \bfu \|_{E,1,\omega_K}
 \end{split}
\end{equation}
where we have used $(\bfu^-_E - \bfu^+_E)\times \bfn =0$ on $\Gamma^{\omega_K}$, \eqref{lem_geo_gamma_eq0} and the trace inequality from Lemma 3.1 in~\cite{2016WangXiaoXu}.


\begin{figure}[h]
\centering
\begin{subfigure}{.31\textwidth}
    \includegraphics[width=1.3in]{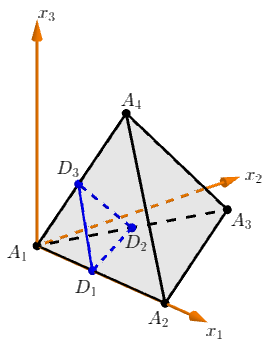}
     \label{fig:regulartet2} 
\end{subfigure}
~
\begin{subfigure}{.31\textwidth}
     \includegraphics[width=1.3in]{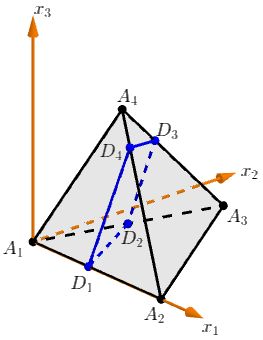}
     \label{fig:regulartet3} 
\end{subfigure}
     \caption{Two types of interface-cutting configuration for a regular tetrahedron. Left (Case 1): 3 intersection points and right (Case 2): 4 intersection points.}
  \label{fig:regulartetInter} 
\end{figure}

\section{Proof of Lemma~\ref{lem_Hdiv_unisol} for a regular tetrahedron}
\label{append:unisol}

We treat the cases of 3 or 4 intersection points shown in Figure~\ref{fig:regulartetInter}, separatly, where $D_i$ are the intersection points of the approximate interface $\Gamma^K_h$ with the edges. Without loss of generality, we assume the tetrahedron has the vertices $A_1(0,0,0)$, $A_2(1,0,0)$, $A_3(1/2,\sqrt{3}/2,0)$, $A_4(1/2,\sqrt{3}/6,\sqrt{6}/3)$. The three faces are chosen as $A_1A_2A_3$, $A_1A_3A_4$ and $A_1A_2A_4$. In addition, for Case 1 we let the intersection points be $D_1=d_1A_2$, $D_2=d_2A_3$, $D_3=d_3A_4$, and for Case 2 we let them be $D_1=d_1A_2$, $D_2=d_2A_3$, $D_3=(1-d_3)A_4+d_3A_3$ and  $D_4=(1-d_4)A_4+d_4A_2$, where $d_i\in[0,1]$, $i=1,2,3,4$, are parameters of their location. For Case 2, since $D_i$ are coplanar, $d_4$ can be expressed as
\begin{equation}
\label{Hcurl_lem_unisol_eq3}
d_4 = \frac{d_1d_3(d_2-1)}{d_1d_2 + d_2d_3 - d_2 - d_1d_3}.
\end{equation}
We use Mathematica to express the determinant as
\begin{equation}
\label{Hcurl_lem_unisol_eq4}
\text{det}(M_K) = \frac{f_1(d_1,d_2,d_3) + f_2(d_1,d_2,d_3) \tilde{\alpha} + f_3(d_1,d_2,d_3) \tilde{\alpha}^2}{f_0(d_1,d_2,d_3)}
\end{equation}
where $f_i(d_1,d_2,d_3)$, $i=0,1,2,3$, are some polynomials of $d_1,d_2,d_3$. We omit the explicit formulas for the functions $f_i$ here since they are complicated. When $\tilde{\alpha}=1$, i.e., $A_K=I$, it is straightforward to calculate $\text{det}(M_K)=-1/16$, and thus
\begin{equation}
\label{Hcurl_lem_unisol_eq5}
f_1(d_1,d_2,d_3) + f_2(d_1,d_2,d_3)  + f_3(d_1,d_2,d_3) = -\frac{1}{16} f_0(d_1,d_2,d_3).
\end{equation}
For Case 1, we can verify that $f_i(d_1,d_2,d_3) \le 0, ~ i =1,2,3$ and $f_0(d_1,d_2,d_3) \ge 0$.
Thus, by \eqref{Hcurl_lem_unisol_eq5} we have
\begin{equation}
\begin{split}
\label{Hcurl_lem_unisol_eq7}
&f_1(d_1,d_2,d_3) + f_2(d_1,d_2,d_3) \tilde{\alpha} + f_3(d_1,d_2,d_3) \tilde{\alpha}^2 \\
 \le & \min\{\tilde{\alpha}^2,1 \} \sum_{i=1,2,3} f_i(d_1,d_2,d_3) =- \frac{\min\{\tilde{\alpha}^2,1 \}}{16} f_0(d_1,d_2,d_3).
\end{split}
\end{equation}
For Case 2, we can veriify that $f_i(d_1,d_2,d_3) \ge 0, ~ i =1,2,3$ and $f_0(d_1,d_2,d_3) \le 0$.
Similar to \eqref{Hcurl_lem_unisol_eq7}, we obtain
\begin{equation}
\label{Hcurl_lem_unisol_eq9}
f_1(d_1,d_2,d_3) + f_2(d_1,d_2,d_3) \tilde{\alpha} + f_3(d_1,d_2,d_3) \tilde{\alpha}^2 \ge - \frac{\min\{\tilde{\alpha}^2,1 \}}{16} f_0(d_1,d_2,d_3).
\end{equation}
Combining \eqref{Hcurl_lem_unisol_eq7} and \eqref{Hcurl_lem_unisol_eq7}, for both cases we arrive at
\begin{equation}
\label{Hcurl_lem_unisol_eq10}
\text{det}(M_K) \le  - \frac{\min\{\tilde{\alpha}^2,1 \}}{16}  < 0.
\end{equation}

\end{appendices}

%

\begin{thebibliography}{10}

\bibitem{2006AlbaneseMonk}
{\sc R.~Albanese and P.~B. Monk}, {\em The inverse source problem for
  {M}axwell's equations}, Inverse Problems, 22 (2006), pp.~1023--1035.

\bibitem{2002AmmariBaoFleming}
{\sc H.~Ammari, G.~Bao, and J.~L. Fleming}, {\em An inverse source problem for
  {M}axwell's equations in magnetoencephalography}, SIAM J. Appl. Math., 62
  (2002).

\bibitem{2015AmmariChenChenVolkov}
{\sc H.~Ammari, J.~Chen, Z.~Chen, D.~Volkov, and H.~Wang}, {\em Detection and
  classification from electromagnetic induction data}, J. Comput. Phys., 301
  (2015), pp.~201 -- 217.

\bibitem{2020VeigaDassiMascotto}
{\sc L.~B. ao~da Veiga, F.~Dassi, G.~Manzini, and L.~Mascotto}, {\em Virtual
  elements for maxwell's equations}, Comput. Math. Appl.,  (2021).

\bibitem{1997ARNOLDRICHARDWINTHER}
{\sc D.~N. Arnold, R.~S. Falk, and R.~Winther}, {\em Preconditioning in
  {H}(div) and applications}, Math. Comp., 66 (1997).

\bibitem{1994BabuskaCalozOsborn}
{\sc I.~Babu{\v{s}}ka, G.~Caloz, and J.~E. Osborn}, {\em Special finite element
  methods for a class of second order elliptic problems with rough
  coefficients}, SIAM J. Numer. Anal., 31 (1994), pp.~945--981.

\bibitem{2001Bathe}
{\sc K.~Bathe}, {\em The inf--sup condition and its evaluation for mixed finite
  element methods}, Computers \& Structures, 79 (2001), pp.~243--252.

\bibitem{2009BeckerBurmanHansbo}
{\sc R.~Becker, E.~Burman, and P.~Hansbo}, {\em A {N}itsche extended finite
  element method for incompressible elasticity with discontinuous modulus of
  elasticity}, Comput. Methods Appl. Mech. Engrg., 198 (2009), pp.~3352--3360.

\bibitem{2001BenBuffaMaday}
{\sc F.~Ben~Belgacem, A.~Buffa, and Y.~Maday}, {\em The mortar finite element
  method for 3d {M}axwell equations: First results}, SIAM J. Numer. Anal., 39
  (2001), pp.~880--901.

\bibitem{1996BrauerRuehl}
{\sc J.~R. Brauer, J.~J. Ruehl, M.~A. Juds, M.~J.~V. Heiden, and A.~A.
  Arkadan}, {\em Dynamic stress in magnetic actuator computed by coupled
  structural and electromagnetic finite elements}, IEEE Trans. Magn., 32
  (1996).

\bibitem{2008BrennerCuiLiSung}
{\sc S.~C. Brenner, J.~Cui, F.~Li, and L.~Y. Sung}, {\em A nonconforming finite
  element method for a two-dimensional curl--curl and grad-div problem}, Numer.
  Math., 109 (2008), pp.~509--533.

\bibitem{2005BuffaCostabelDauge}
{\sc A.~Buffa, M.~Costabel, and M.~Dauge}, {\em Algebraic convergence for
  anisotropic edge elements in polyhedral domains}, Numer. Math., 101 (2005),
  pp.~29--65.

\bibitem{2015BurmanClaus}
{\sc E.~Burman, S.~Claus, P.~Hansbo, M.~G. Larson, and A.~Massing}, {\em
  Cut{FEM}: Discretizing geometry and partial differential equations},
  Internat. J. Numer. Methods Engrg., 104 (2015), pp.~472--501.

\bibitem{2016BurmanGuzmanSanchezSarkisHcontrast}
{\sc E.~Burman, J.~Guzm{\'a}n, M.~A. S{\'a}nchez, and M.~Sarkis}, {\em Robust
  flux error estimation of an unfitted nitsche method for high-contrast
  interface problems}, IMA J. Numer. Anal., 38 (2018), pp.~646--668.

\bibitem{2015CaiCao}
{\sc Z.~Cai and S.~Cao}, {\em A recovery-based a posteriori error estimator for
  {H}(curl) interface problems}, Comput. Methods Appl. Mech. Engrg., 296
  (2015), pp.~169 -- 195.

\bibitem{2018CaoChen}
{\sc S.~Cao and L.~Chen}, {\em Anisotropic error estimates of the linear
  virtual element method on polygonal meshes}, SIAM J. Numer. Anal., 56 (2018),
  pp.~2913--2939.

\bibitem{2021CaoChenGuo}
{\sc S.~Cao, L.~Chen, and R.~Guo}, {\em A virtual finite element method for two
  dimensional {M}axwell interface problems with a background unfitted mesh},
  Math. Models Methods Appl. Sci., 31 (2021), pp.~2907--2936.

\bibitem{2022CaoChenGuo}
\leavevmode\vrule height 2pt depth -1.6pt width 23pt, {\em Immersed virtual
  element methods for {M}axwell interface problems in three dimensions}, ArXiv
  preprint,  (2022).

\bibitem{2016CasagrandeHiptmairOstrowski}
{\sc R.~Casagrande, R.~Hiptmair, and J.~Ostrowski}, {\em An a priori error
  estimate for interior penalty discretizations of the curl-curl operator on
  non-conforming meshes}, J. Math. Ind., 6 (2016), p.~4.

\bibitem{2016CasagrandeWinkelmannHiptmairOstrowski}
{\sc R.~Casagrande, C.~Winkelmann, R.~Hiptmair, and J.~Ostrowski}, {\em {DG}
  treatment of non-conforming interfaces in 3{D} {C}url-{C}url problems}, in
  Scientific Computing in Electrical Engineering, Cham, 2016, Springer
  International Publishing, pp.~53--61.

\bibitem{1993ChapelleBathe}
{\sc D.~Chapelle and K.~Bathe}, {\em The inf-sup test}, Computers \&
  Structures, 47 (1993), pp.~537--545.

\bibitem{2020CHENLIANGZOU}
{\sc J.~Chen, Y.~Liang, and J.~Zou}, {\em Mathematical and numerical study of a
  three-dimensional inverse eddy current problem}, SIAM J. Appl. Math., 80
  (2020), pp.~1467--1492.

\bibitem{Chen:2008ifem}
{\sc L.~Chen}, {\em $i$fem: an integrated finite element methods package in
  {MATLAB}}, tech. rep., 2009.

\bibitem{2017ChenWeiWen}
{\sc L.~Chen, H.~Wei, and M.~Wen}, {\em An interface-fitted mesh generator and
  virtual element methods for elliptic interface problems}, J. Comput. Phys.,
  334 (2017), pp.~327--348.

\bibitem{2000ChenDuZou}
{\sc Z.~Chen, Q.~Du, and J.~Zou}, {\em Finite element methods with matching and
  nonmatching meshes for {M}axwell equations with discontinuous coefficients},
  SIAM J. Numer. Anal., 37 (2000), pp.~1542--1570.

\bibitem{2007ChenWangZheng}
{\sc Z.~Chen, L.~Wang, and W.~Zheng}, {\em An adaptive multilevel method for
  time‐harmonic {M}axwell equations with singularities}, SIAM J. Sci.
  Comput., 29 (2007), pp.~118--138.

\bibitem{2009ChenXiaoZhang}
{\sc Z.~Chen, Y.~Xiao, and L.~Zhang}, {\em The adaptive immersed interface
  finite element method for elliptic and {M}axwell interface problems}, J.
  Comput. Phys., 228 (2009), pp.~5000 -- 5019.

\bibitem{2010ChuGrahamHou}
{\sc C.-C. Chu, I.~G. Graham, and T.-Y. Hou}, {\em A new multiscale finite
  element method for high-contrast elliptic interface problems}, Math. Comp.,
  79 (2010), pp.~1915--1955.

\bibitem{1999CiarletZou}
{\sc P.~Ciarlet, Jr and J.~Zou}, {\em Fully discrete finite element approaches
  for time-dependent {M}axwell's equations}, Numer. Math., 82 (1999),
  pp.~193--219.

\bibitem{2000CostabelDauge}
{\sc M.~Costabel and M.~Dauge}, {\em Singularities of electromagnetic fields in
  polyhedral domains}, Arch. Ration. Mech. Anal., 151 (2000), pp.~221--276.

\bibitem{1999MartinMoniqueSerge}
{\sc M.~Costabel, M.~Dauge, and S.~Nicaise}, {\em Singularities of {M}axwell
  interface problems.}, ESAIM: M2AN, 33 (1999), pp.~627--649.

\bibitem{2004CostabelDaugeNicaise}
\leavevmode\vrule height 2pt depth -1.6pt width 23pt, {\em Corner Singularities
  of Maxwell Interface and Eddy Current Problems}, Birkh{\"a}user Basel, Basel,
  2004, pp.~241--256.

\bibitem{2020BeiroMascotto}
{\sc L.~B. da~Veiga and L.~Mascotto}, {\em Interpolation and stability
  properties of low order face and edge virtual element spaces}, IMA J. Numer.
  Anal.,  (2022).

\bibitem{2014ErcanJaewook}
{\sc E.~M. Dede, J.~Lee, and T.~Nomura}, {\em Multiphysics Simulation:
  Electromechanical System Applications and Optimization}, Springer, 2014.

\bibitem{2000Edelsbrunner}
{\sc H.~Edelsbrunner}, {\em Triangulations and meshes in computational
  geometry}, Acta Numerica, 9 (2000), pp.~133--213.

\bibitem{2007GroReusken}
{\sc S.~Gro{\ss} and A.~Reusken}, {\em An extended pressure finite element
  space for two-phase incompressible flows with surface tension}, J. Comput.
  Phys., 224 (2007), pp.~40 -- 58.
\newblock Special Issue Dedicated to Professor Piet Wesseling on the occasion
  of his retirement from Delft University of Technology.

\bibitem{2020Guo}
{\sc R.~Guo}, {\em Solving parabolic moving interface problems with dynamical
  immersed spaces on unfitted meshes: Fully discrete analysis}, SIAM J. Numer.
  Anal., 2 (2021), pp.~797--828.

\bibitem{2020GuoLin}
{\sc R.~Guo and T.~Lin}, {\em An immersed finite element method for elliptic
  interface problems in three dimensions}, J. Comput. Phys., 414 (2020).

\bibitem{2020GuoLinZou}
{\sc R.~Guo, Y.~Lin, and J.~Zou}, {\em Solving two dimensional
  {H}(curl)-elliptic interface systems with optimal convergence on unfitted
  meshes}, arXiv:2011.11905,  (2020).

\bibitem{2020GuoLinLinZhuang}
{\sc R.~Guo, Y.~L. T.~Lin, and Q.~Zhuang}, {\em Error analysis of symmetric
  linear/bilinear partially penalized immersed finite element methods for
  helmholtz interface problems}, J. Comput. Appl. Math.,  (2021).

\bibitem{2014HansboLarsonZahedi}
{\sc P.~Hansbo, M.~G. Larson, and S.~Zahedi}, {\em A cut finite element method
  for a stokes interface problem}, Appl. Numer. Math., 85 (2014), pp.~90--114.

\bibitem{2010HiptmairLiZou}
{\sc R.~Hiptmair, J.~LI, and J.~Zou}, {\em Convergence analysis of finite
  element methods for ${H}(div;{\Omega})$-elliptic interface problems}, J.
  Numer. Math., 18 (2010), pp.~187--218.

\bibitem{2012HiptmairLiZou}
{\sc R.~Hiptmair, J.~Li, and J.~Zou}, {\em Convergence analysis of finite
  element methods for {H}(curl; {$\Omega$})-elliptic interface problems},
  Numer. Math., 122 (2012), pp.~557--578.

\bibitem{2011HIPTMAIRMOIOLAPERUGIA}
{\sc R.~Hiptmair, A.~Moiola, and I.~Perugia}, {\em Stability results for the
  time-harmonic {M}axwell equations with impedance boundary conditions}, Math.
  Models Methods Appl. Sci., 21 (2011), pp.~2263--2287.

\bibitem{2007HiptmairXu}
{\sc R.~Hiptmair and J.~Xu}, {\em Nodal auxiliary space preconditioning in
  {H}(curl) and {H}(div) spaces}, SIAM J. Numer. Anal., 45 (2007),
  pp.~2483--2509.

\bibitem{2013HouSongWangZhao}
{\sc S.~Hou, P.~Song, L.~Wang, and H.~Zhao}, {\em A weak formulation for
  solving elliptic interface problems without body fitted grid}, J. Comput.
  Phys., 249 (2013), pp.~80 -- 95.

\bibitem{2005HoustonPerugiaSchneebeli}
{\sc P.~Houston, I.~Perugia, A.~Schneebeli, and D.~Sch{\"o}tzau}, {\em Interior
  penalty method for the indefinite time-harmonic {M}axwell equations}, Numer.
  Math., 100 (2005), pp.~485--518.

\bibitem{2004HoustonPerugiaSchotzau}
{\sc P.~Houston, I.~Perugia, and D.~Schotzau}, {\em Mixed discontinuous
  galerkin approximation of the {M}axwell operator}, SIAM J. Numer. Anal., 42
  (2004), pp.~434--459.

\bibitem{2021HuWang}
{\sc J.~Hu and H.~Wang}, {\em An optimal multigrid algorithm for the combining
  $p_1$-$q_1$ finite element approximations of interface problems based on
  local anisotropic fitting meshes}, J. Sci. Comput., 88 (2021).

\bibitem{2008HuShuZou}
{\sc Q.~Hu, S.~Shu, and J.~Zou}, {\em A mortar edge element method with nearly
  optimal convergence for three-dimensional {M}axwell's equations}, Math.
  Comp., 77 (2008), pp.~1333--1353.

\bibitem{2005KafafyLinLinWang}
{\sc R.~Kafafy, T.~Lin, Y.~Lin, and J.~Wang}, {\em Three-dimensional immersed
  finite element methods for electric field simulation in composite materials},
  Internat. J. Numer. Methods Engrg., 64 (2005), pp.~940--972.

\bibitem{2010LiMelenkWohlmuthZou}
{\sc J.~Li, J.~M. Melenk, B.~Wohlmuth, and J.~Zou}, {\em Optimal a priori
  estimates for higher order finite elements for elliptic interface problems},
  Appl. Numer. Math., 60 (2010), pp.~19--37.

\bibitem{2015LinLinZhang}
{\sc T.~Lin, Y.~Lin, and X.~Zhang}, {\em Partially penalized immersed finite
  element methods for elliptic interface problems}, SIAM J. Numer. Anal., 53
  (2015), pp.~1121--1144.

\bibitem{2020LiuZhangZhangZheng}
{\sc H.~Liu, L.~Zhang, X.~Zhang, and W.~Zheng}, {\em Interface-penalty finite
  element methods for interface problems in {$H^1$}, {H}(curl), and {H}(div)},
  Comput. Methods Appl. Mech. Engrg., 367 (2020).

\bibitem{1986MakinsonShah}
{\sc G.~Makinson and A.~Shah}, {\em An iterative solution method for solving
  sparse nonsymmetric linear systems}, Journal of Computational and Applied
  Mathematics, 15 (1986), pp.~339--352.

\bibitem{2003Monk}
{\sc P.~Monk}, {\em Finite Element Methods for {M}axwell's Equations}, Oxford
  University Press, 2003.

\bibitem{1980Nedelec}
{\sc J.~C. Nedelec}, {\em Mixed finite elements in $\mathbb{R}^3$}, Numer.
  Math., 35 (1980), pp.~315--341.

\bibitem{1986Nedelec}
{\sc J.~C. N{\'e}d{\'e}lec}, {\em A new family of mixed finite elements in
  $\mathbb{R}^3$}, Numer. Math., 50 (1986), pp.~57--81.

\bibitem{2001OsherFedkiw}
{\sc S.~Osher and R.~P. Fedkiw}, {\em Level set methods: An overview and some
  recent results}, J. Comput. Phys., 169 (2001), pp.~463--502.

\bibitem{2004PerssonStrang}
{\sc P.-O. Persson and G.~Strang}, {\em A simple mesh generator in {M}atlab},
  SIAM Review, 46 (2004), pp.~329--345.

\bibitem{1977RaviartThomas}
{\sc P.~A. Raviart and J.~M. Thomas}, {\em A mixed finite element method for
  2-nd order elliptic problems}, in Mathematical Aspects of Finite Element
  Methods, I.~Galligani and E.~Magenes, eds., Berlin, Heidelberg, 1977,
  Springer Berlin Heidelberg, pp.~292--315.

\bibitem{2020RoppertSchoderToth}
{\sc K.~Roppert, S.~Schoder, F.~Toth, and M.~Kaltenbacher}, {\em Non-conforming
  nitsche interfaces for edge elements in curl--curl-type problems}, IEEE
  Trans. Magn., 56 (2020).

\bibitem{1999Sethian}
{\sc J.~A. Sethian}, {\em Level set methods and fast marching methods}, vol.~3
  of Cambridge Monographs on Applied and Computational Mathematics, Cambridge
  University Press, Cambridge, second~ed., 1999.

\bibitem{2016WangXiaoXu}
{\sc F.~Wang, Y.~Xiao, and J.~Xu}, {\em High-order extended finite element
  methods for solving interface problems}, Comput. Methods Appl. Mech. Engrg.,
  364 (2020).

\bibitem{2003WarburtonHesthaven}
{\sc T.~Warburton and J.~S. Hesthaven}, {\em On the constants in {$hp$}-finite
  element trace inverse inequalities}, Comput. Methods Appl. Mech. Engrg., 192
  (2003), pp.~2765--2773.

\bibitem{2008XuZhu}
{\sc J.~XU and Y.~ZHU}, {\em Uniform convergent multigrid methods for elliptic
  problems with strongly discontinuous coefficients}, Mathematical Models and
  Methods in Applied Sciences, 18 (2008), pp.~77--105.

\bibitem{1998XuZou}
{\sc J.~Xu and J.~Zou}, {\em Some nonoverlapping domain decomposition methods},
  SIAM Review, 40 (1998), pp.~857--914.

\bibitem{1994CaiLazarovManteuffelMcCormick}
{\sc R.~L. Z.~Cai, T.~A. Manteuffel, and S.~F. McCormick}, {\em First-order
  system least squares for second-order partial differential equations: Part
  i}, SIAM J. Numer. Anal., 31 (1994).

\bibitem{2004ZhaoWei}
{\sc S.~Zhao and G.~W. Wei}, {\em High-order {FDTD} methods via derivative
  matching for {M}axwell's equations with material interfaces}, J. Comput.
  Phys., 200 (2004), pp.~60--103.

\end{thebibliography}

\end{document}